\documentclass[aos,seceqn,nameyear,dvips]{arximspdf}
\usepackage{accents}

\doi{10.1214/10-AOS810}
\volume{38}
\issue{6}
\pubyear{2010}
\firstpage{3245}
\lastpage{3299}

\makeatletter

\newtheorem{Lem}{Lemma}[section]
\newtheorem{Prop}{Proposition}[section]
\newproclaim{assumption}{Assumption}

\newcommand{\ort}{^{\bot}}
\newcommand{\R}{\mathbb R}
\newcommand{\rr}{\mathbb R}
\newcommand{\N}{\mathbb N}
\newcommand{\sirc}{{\circ}}
\newcommand{\Xb}{\mathbf{X}}
\newcommand{\Sb}{\mathbf{S}}

\newcommand{\Vb}{\mathbf{V}}

\newcommand{\hats}[1]{\hspace*{1.5pt}\hat{\hspace*{-1.5pt}#1}}

\newcommand{\Bb}{\mathbf{B}}
\newcommand{\Hb}{\mathbf{H}}
\newcommand{\Ab}{\mathbf{A}}
\newcommand{\Ub}{\mathbf{U}}
\newcommand{\Mb}{\mathbf{M}}
\newcommand{\Yb}{\mathbf{Y}}
\newcommand{\Gb}{\mathbf{G}}
\newcommand{\wb}{\mathbf{w}}

\newcommand{\I}{{I}}
\newcommand{\II}{\mathit{II}}
\newcommand{\III}{\mathit{III}}
\newcommand{\IV}{\mathit{IV}}

\newcommand{\lam}{\lambda}

\newcommand{\Omegab}{{\bolds\Omega}}
\newcommand{\iotab}{{\bolds\iota}}
\newcommand{\Upsib}{{\bolds\Upsilon}}
\newcommand{\tetb}{{\bolds\theta}}
\newcommand{\varthetab}{{\bolds\vartheta}}
\newcommand{\Thetab}{{\bolds\Theta}}
\newcommand{\Lamb}{{\bolds\Lambda}}
\newcommand{\Sigb}{{\bolds\Sigma}}
\newcommand{\Deltab}{{\bolds\Delta}}
\newcommand{\taub}{{\bolds\tau}}
\newcommand{\etab}{{\bolds\eta}}
\newcommand{\betab}{{\bolds\beta}}
\newcommand{\Xib}{{\bolds\Xi}}
\newcommand{\Gamb}{{\bolds\Gamma}}

\newcommand{\ikf}{\mathcal{I}_k(f_1)}

\newcommand{\vechop}{\operatorname{vech}}
\newcommand{\dvecrond}{\operatorname{\mathrm{d}\hspace*{-6pt}\accentset{\circ}{\phantom{a}\mathrm{vec}}}}
\newcommand{\dvecronds}{\mathrm{d}\hspace*{-4pt}\accentset{\circ}{\phantom{a}\mathrm{vec}}}
\newcommand{\dvec}{\operatorname{dvec} }
\newcommand{\vecop}{\operatorname{vec}}
\newcommand{\pr}{^{\prime}}
\newcommand{\npr}{^{(n)\prime}}

\newcommand{\utphi}{\mathop{\phi}\limits_{\widetilde{}}}
\newcommand{\utQ}{\mathop{Q}\limits_{\widetilde{}}}
\newcommand{\utT}{\mathop{T}\limits_{\widetilde{}}}
\newcommand{\utDelta}{\mathop{\bolds\Delta}\limits_{\widetilde{}}}
\newcommand{\ubDelta}{{{\bolds\Delta}}}
\newcommand{\utSb}{\mathop{\mathbf{S}}\limits_{\hspace*{-0.5pt}\widetilde{}}}
\newcommand{\ny}{n\rightarrow\infty}

\newcommand{\btilde}{{\tilde{\phantom u}}\hspace{-1.8mm}b}
\newcommand{\bbar}{{\bar{\phantom u}}\hspace{-1.8mm}b}
\newcommand{\lbar}{{\bar{\phantom u}}\hspace{-1.8mm}l}
\newcommand{\dbar}{\hspace*{1pt}{\bar{\phantom u}}\hspace{-2.4mm}d}
\newcommand{\dbarr}{{\bar{\phantom u}}\hspace{-1.8mm}d}
\makeatother

\begin{document}
\begin{frontmatter}

\title{Optimal rank-based testing for principal components}
\runtitle{Rank-based tests for PCA}

\begin{aug}
\author[A]{\fnms{Marc} \snm{Hallin}\corref{}\thanksref{t2,t1,t5}\ead[label=e1]{mhallin@ulb.ac.be}},
\author[B]{\fnms{Davy} \snm{Paindaveine}\thanksref{t4,t5}\ead[label=e2]{dpaindav@ulb.ac.be}} and
\author[C]{\fnms{Thomas} \snm{Verdebout}\ead[label=e3]{thomas.verdebout@univ-lille3.fr}\ead[label=u1,url]{http://homepages.ulb.ac.be/\textasciitilde dpaindav}}
\runauthor{M. Hallin, D. Paindaveine and T. Verdebout}
\affiliation{Universit\'e Libre de Bruxelles and Princeton  University, Universit\'e Libre
de~Bruxelles~and~Universit\' e Lille Nord de France}
\address[A]{M. Hallin\\
Institut de Recherche en Statistique\\
E.C.A.R.E.S., CP114\\
Universit\'{e} Libre de Bruxelles\\
50, Avenue F. D. Roosevelt\\
B-1050 Bruxelles\\
Belgium\\
\printead{e1}}
\address[B]{D. Paindaveine\\
Institut de Recherche en Statistique\\
E.C.A.R.E.S., CP114\\
\quad and D\'{e}partement de Math\'{e}matique\\
Universit\'{e} Libre de Bruxelles\\
50, Avenue F. D. Roosevelt \\
B-1050 Bruxelles\\
Belgium\\
\printead{e2} \\
\printead{u1}}
\address[C]{T. Verdebout\\
EQUIPPE-GREMARS \\
Universit\' e Lille III\\
Domaine Universitaire du Pont de Bois, BP 60149 \\
F-59653 Villeneuve d'Ascq Cedex\\
France\\
\printead{e3}}
\end{aug}

\thankstext{t2}{Supported by the Sonderforschungsbereich ``Statistical
modelling of nonlinear dynamic processes'' (SFB 823) of the
Deutsche Forschungsgemeinschaft, and by a Discovery grant of
the Australian Research Council.}
\thankstext{t4}{Supported by a contract of the National Bank of Belgium and
a Mandat d'Impulsion Scientifique of the
Fonds National de la Recherche Scientifique,
Communaut\' e fran\c caise de Belgique.}
\thankstext{t1}{Member
of Acad\'emie Royale de Belgique and
extra-muros Fellow of CenTER, Tilburg University.}
\thankstext{t5}{Member
of ECORE, the association between CORE and ECARES.}

% HISTORY:
\received{\smonth{4} \syear{2009}}
\revised{\smonth{1} \syear{2010}}

% ABSTRACT
%
\begin{abstract}
This paper provides parametric and rank-based optimal tests for
eigenvectors and eigenvalues of covariance or scatter matrices in
elliptical families. The parametric tests extend the Gaussian
likelihood ratio tests of \citet{A63} and their pseudo-Gaussian
robustifications by \citet{D77} and Tyler (\citeyear{T81}, \citeyear{T83}). The rank-based
tests address a much broader class of problems, where covariance
matrices need not exist and principal components are associated with
more general scatter matrices. The proposed tests are shown to
outperform daily practice both from the point of view of validity as
from the point of view of efficiency. This is achieved by utilizing the
Le Cam theory of locally asymptotically normal experiments, in the
nonstandard context, however, of a \textit{curved} parametrization. The
results we derive for curved experiments are of independent interest,
and likely to apply in other contexts.
\end{abstract}

% KEYWORDS
%

\end{frontmatter}

%s1 ###
\section{Introduction}\label{intro}

This fairly detailed introduction aims at providing a comprehensive and
nontechnical overview of the paper, including its asymptotic theory
aspects, and a rough description of some of the rank-based test
statistics to be derived. It is expected to be accessible
to a broad readership. It should be sufficiently informative for the
reader not interested in the technical aspects of asymptotic theory, to
proceed to Sections~\ref{gausscase} (Gaussian and pseudo-Gaussian tests)
and~\ref{ranktests} (rank-based tests), where the proposed testing procedures are
described, and for the reader mainly interested in asymptotics, to
decide whether he/she is interested in the treatment of a LAN family
with curved parametrization developed in Sections~\ref{LANsection} and~\ref{paramtests}.

%s1.1 ###
\subsection{Hypothesis testing for principal components}

Principal components
are probably the most popular and widely used device in the traditional
multivariate analysis toolkit.
Introduced by \citet{P1901}, principal component analysis (PCA) was
rediscovered by \citet{H33}, and ever since has been an essential
part of daily statistical practice, basically in all domains of application.

The general objective of PCA is to reduce the dimension of some
observed $k$-dimensional random vector $\Xb$ while preserving most of
its total variability. This is achieved by considering an adequate
number $q$ of linear combinations of the form $\betab_1\pr\Xb,\ldots
, \betab_q\pr\Xb$, where $\betab_j$, $j=1,\ldots, k$, are the
eigenvectors associated with the eigenvalues $\lambda_1,\ldots,
\lambda_k$ of $\Xb$'s covariance matrix $\Sigb_{\mathrm{cov}}$, ranked
in decreasing order of magnitude. Writing $\betab$ for the orthogonal
$k\times k$ matrix with columns $\betab_1 ,\ldots, \betab_k$ and $
\Lamb_{\Sigb_{\mathrm{cov}}}$ for the diagonal matrix of eigenvalues
$\lambda_1,\ldots, \lambda_k$, the matrix $\Sigb_{\mathrm{cov}}$ thus
factorizes into $\Sigb_{\mathrm{cov}}=\betab\Lamb_{\Sigb_{\mathrm{cov}}}
\betab\pr$. The random variable $\betab_j\pr\Xb$, with
variance $\lambda_j$, is known as $\Xb$'s \textit{$j$th principal
component}.

Chapters on inference for eigenvectors
and eigenvalues can be found in most textbooks on multivariate
analysis, and mainly cover Gaussian maximum likelihood estimation (MLE)
and the corresponding Wald and Gaussian likelihood ratio tests (LRT).
The MLEs of $\betab$ and $\Lamb_{\Sigb_{\mathrm{cov}}}$ are the
eigenvectors and eigenvalues of the empirical covariance matrix
\[
\mathbf{S}^{(n)}:=\frac{1}{n}\sum_{i=1}^n \bigl(\mathbf{X}_i-\bar{\mathbf
X}^{(n)}\bigr)\bigl(\mathbf{X}_i-\bar{\mathbf X}^{(n)}\bigr)\pr\qquad \mbox{with }
\bar{\mathbf X}^{(n)}:=\frac{1}{n}\sum_{i=1}^n \mathbf{X}_i,
\]
while testing problems classically include testing for sphericity
(equality of eigenvalues), testing for \textit{subsphericity} (equality
among some given subset of eigenva\-lues---typically, the last $k-q$
ones), testing that the $\ell$th eigenvector has some specified
direction, or that the proportion of variance accounted for by the last
$k-q$ principal components is larger than some fixed proportion of the
total variance: see, for instance, \citet{A03} or \citet{J86}.

Gaussian MLEs and the corresponding tests (Wald or likelihood ratio
tests---since they are asymptotically equivalent, in the sequel we
indistinctly refer to LRTs) for covariance matrices and functions
thereof are notoriously sensitive to violations of Gaussian
assumptions; see \citet{MW80} for a classical
discussion of this fact, or \citet{YTM05} for a more recent overview.
The problems just mentioned about the eigenvectors and eigenvalues of
$\Sigb_{\mathrm{cov}}$ are no exception to that rule, although belonging,
in Muirhead and Waternaux's terminology, to the class of ``easily
robustifiable'' ones. For such problems, \textit{adjusted} LRTs remaining
valid under the whole class of elliptical distributions with finite
fourth-order moments can be obtained via a correction factor involving
estimated kurtosis coefficients [see \citet{SB87} for a
general result on the ``easy'' cases, and \citet{HP08c}
for the ``harder'' ones]. Such adjusted LRTs were obtained by
Tyler (\citeyear{T81}, \citeyear{T83}) for eigenvector problems and
by \citet{D77} for eigenvalues.

Tyler actually constructs tests for the \textit{scatter matrix} $\Sigb$
characterizing the density contours [of the form $(\mathbf{x}-\tetb)\pr
\Sigb^{-1}(\mathbf{x}-\tetb)=$ constant] of an elliptical family. His
tests are the Wald tests associated with any available
estimator $\hat{\Sigb}$ of ${\Sigb}$ such that $n^{1/2} \vecop
(\hat{\Sigb}- \Sigb)$ is asymptotically normal, with mean zero and
covariance matrix ${\bolds\Psi} _{f}$, say, under ${f}\in
\mathcal{F}$, where $\mathcal F$ denotes some class of elliptical
densities and ${\bolds\Psi} _{f}$ either is known or (still, under
${f}\in\mathcal{F}$) can be estimated consistently. The resulting
tests then are valid under the class $\mathcal F$. When the estimator
$\hat{\Sigb}$ is the empirical covariance matrix $\mathbf{S}^{(n)}$, these
tests under Gaussian densities are asymptotically equivalent to
Gaussian LRTs. Unlike the latter, however, they remain
(asymptotically) valid under the class $\mathcal{F}^4$ of all
elliptical distributions with finite moments of order four, and hence
qualify as \textit{pseudo-Gaussian} versions of the Gaussian LRTs.

Due to their importance for applications, throughout this paper, we
concentrate on the following two problems:

(a) testing the null hypothesis ${\mathcal H}_{0}^{\betab}$ that the
first principal direction $ \betab_{1}$ coincides (up to the sign)
with some specified unit vector $\betab_{}^{0}$ (the choice of the
\textit{first} principal direction here is completely arbitrary, and made
for the simplicity of exposition only), and

(b) testing the null hypothesis ${\mathcal H}_{0}^{\Lamb
}$ that ${\sum_{j=q+1}^{k} \lambda_{j%; \Vb
}}/{\sum_{j=1}^{k}\lambda_{j%; \Vb
}}=p$ against the one-sided alternative under which ${\sum_{j=q+1}^{k}
\lambda_{j%; \Vb
}}/{\sum_{j=1}^{k}\lambda_{j%; \Vb
}}<p$, $p \in(0,1)$ given.

The Gaussian LRT for (a) was introduced in a seminal paper by
\citet{A63}. Denoting by ${\lambda}_{j; \Sb}$ and $\betab_{j;
\Sb}$, $j=1,\ldots,k$, respectively, the eigenvalues and eigenvectors
of $\mathbf{S}^{(n)}$, this test---denote it by $\phi^{(n)}_{\betab
;\mathrm{Anderson}}$---rejects ${\mathcal H}_{0}^{\betab}$ (at
asymptotic level $\alpha$) as soon as
%
%e1.1 ###
%
\begin{eqnarray}\label{AndTesta}
Q_{\mathrm{Anderson}}^{(n)}
:\!&=&
n \bigl[ {\lambda}_{1; \Sb} \betab{}^{0\prime} \bigl(\mathbf{S}^{(n)}\bigr)^{-1}
\betab_{}^{0} + {\lambda}_{1; \Sb}^{-1} \betab^{0\prime} \mathbf
{S}^{(n)}\betab_{}^{0}-2 \bigr] \nonumber\\[-8pt]\\[-8pt] %[1mm]
& = & \frac{n}{\lambda_{1; {\Sb}} }\sum_{j=2}^{k} \frac{(\lambda
_{j; {\Sb}}- \lambda_{1; {\Sb}} )^{2}}{\lambda_{j; {\Sb}}^{3}}
\bigl( \betab_{j; \Sb}\pr\mathbf{S}^{(n)}\betab^{0}
\bigr)^{2}\nonumber
\end{eqnarray}
exceeds the $\alpha$ upper-quantile of the chi-square distribution
with $(k-1)$ degrees of freedom. The behavior of this test being
particularly poor under non-Gaussian densities, Tyler (\citeyear{T81}, \citeyear{T83})
proposed a pseudo-Gaussian version $\phi^{(n)}_{\betab;\mathrm
{Tyler}}$, which he obtains via an empirical kurtosis correction
%
%e1.2 ###
%
\begin{equation}\label{TylTesta}
Q^{(n)}_{\mathrm{Tyler}}:= \bigl(1+ \hat{\kappa}^{(n)}\bigr)^{-1} Q_{\mathrm
{Anderson}}^{(n)}
\end{equation}
of~(\ref{AndTesta}) (same asymptotic distribution), where $\hat
{\kappa}^{(n)}$ is some consistent estimator of the underlying
kurtosis parameter $\kappa_k$; see Section~\ref{pseudovec} for a definition.

A related test of \citet{S91} addresses the same problem where
however $\betab_{1}$ is the first eigenvector of the \textit
{correlation} matrix.

The traditional Gaussian test for problem (b) was introduced in the
same paper by \citet{A63}. For any $k\times k$ diagonal
matrix $\Lamb$ with diagonal entries $\lambda_1,\ldots,\lambda_k$,
let ${a}_{p,q}({\Lamb}):=2 ( p^{2} \sum_{j=1}^{q} {\lambda
}_{j}^{2}+ (1-p)^{2} \sum_{j=q+1}^{k} {\lambda}_{j}^{2} )$.
Defining\vspace*{-2pt} $\betab_\Sb:=(\betab_{1; \Sb},\ldots,\betab_{k; \Sb})$
and $\mathbf{c}_{p,q}:=
( -p \mathbf{1}_{q}\pr{\,}\vdots{\,}(1-p) \mathbf{1}_{k-q}\pr)\pr
$, with\vspace*{1pt} $\mathbf{1}_\ell:= (1,\ldots, 1)\pr\in\rr^\ell$, and
denoting by $\dvec(\Ab)$ the vector
obtained by stacking the diagonal elements of a square matrix $\Ab$,
Anderson's test, $\phi^{(n)}_{\Lamb; \mathrm{Anderson}}$, say,
rejects the null hypothesis at asymptotic level $\alpha$ whenever
%
%e1.3 ###
%
\begin{eqnarray}\label{TAnd}
T_{\mathrm{Anderson}}^{(n)}
:\!&=&
n^{1/2}
({a}_{p,q}({\Lamb}_{\mathbf{S}}))^{-1/2}
\mathbf{c}_{p,q}\pr\dvec\bigl({\betab}_\Sb\pr{\Sb}^{(n)}{\betab}_\Sb\bigr)
\nonumber\\[-8pt]\\[-8pt]
&=&
n^{1/2}
({a}_{p,q}({\Lamb}_{\mathbf{S}}))^{-1/2}
\Biggl( (1-p) \sum_{j=q+1}^{k} {\lambda}_{j;\Sb}- p\sum_{j=1}^{q}
{\lambda}_{j;\Sb} \Biggr)\nonumber
\end{eqnarray}
is less than the standard normal $\alpha$-quantile. Although he does
not provide any explicit form, \citet{D77} briefly explains how to
derive the pseudo-Gaussian version
%
%e1.4 ###
%
\begin{equation}\label{TDav}
T_{\mathrm{Davis}}^{(n)}:= \bigl(1+ \hat{\kappa}{}^{(n)}\bigr)^{-1/2}
T_{\mathrm{Anderson}}^{(n)}
\end{equation}
of~(\ref{TAnd}), where $\hat{\kappa}^{(n)}$ again is any consistent
estimator of the underlying kurtosis parameter $\kappa_k$. The
resulting test (same asymptotic standard normal distribution) will be
denoted as $\phi^{(n)}_{\Lamb; \mathrm{Davis}}$.

Being based on empirical covariances, though, the pseudo-Gaussian tests
based on~(\ref{TylTesta}) and~(\ref{TDav}) unfortunately remain
poorly robust. They still are very sensitive to the presence of
outliers---an issue which we do not touch here; see, for example, \citet{CH00},
\citet{SAW06}, and the
references therein. Moreover, they do require finite moments of order
four---hence lose their validity under heavy tails, and only address
the traditional covariance-based concept of principal components.

This limitation is quite regrettable, as principal components,
irrespective of any moment conditions, clearly depend on the elliptical
geometry of underlying distributions only. Recall that an elliptical
density over ${\R}^{k}$ is determined by a \textit{location vector}
${\bolds\theta}\in\R^{k}$, a~\textit{scale} parameter $\sigma\in\R
_{0}^{+}$ (where $\sigma^2 $ is not necessarily a variance), a~real-valued $k \times k$ symmetric and positive definite matrix ${\Vb
}$ called the \textit{shape} matrix, and a
\textit{standardized radial density} $f_{1}$ (whenever the elliptical
density has finite second-order moments, the shape and covariance
matrices $\mathbf{V}$ and $\Sigb_{\mathrm{cov}}$ are proportional, hence
share the same collection of eigenvectors and, up to a positive factor,
the same collection of eigenvalues). Although traditionally described
in terms of the covariance matrix $\Sigb_{\mathrm{cov}}$, most
inference problems in multivariate analysis naturally extend to
arbitrary elliptical models, with the shape matrix $\Vb$ or the
\textit{scatter matrix} $\Sigb:= \sigma^2\Vb$ playing the role of $\Sigb
_{\mathrm{cov}}$. Principal components are no exception; in
particular, problems (a) and (b) indifferently can be formulated in
terms of shape or covariance eigenvectors and eigenvalues. Below,
$\Lamb_{\Vb}:=\operatorname{diag}(\lambda_{1;\Vb},\ldots,\lambda_{k;\Vb
})$ and $\betab:=(\betab_1,\ldots,\betab_k)$ collect the
eigenvalues and eigenvectors of the shape matrix $\Vb$.

Our objective in this paper is to provide a class of signed-rank
tests which remain valid under arbitrary elliptical densities, in the
absence of \textit{any} moment assumption, and hence are not limited to
the traditional covariance-based concept of principal components. Of
particular interest within that class are the \textit{van der
Waerden}---that is, \textit{normal-score}---tests, which are
asymptotically equivalent, under Gaussian densities, to the
corresponding Gaussian LRTs (the asymptotic optimality of which we
moreover establish in Section~\ref{gausscase}, along with local
powers). Under non-Gaussian conditions, however, these van der Waerden
tests uniformly dominate, in the Pitman sense,
the pseudo-Gaussian tests based on~(\ref{TylTesta}) and~(\ref{TDav})
above, which, as a result, turn out to be nonadmissible (see
Section~\ref{secare}).

Our tests are based on the multivariate signs and ranks previously
considered by Hallin and Paindaveine
(\citeyear{HP06a}, \citeyear{HP08a}) and \citet{HOP06}. Denote by $\mathbf{X}_1, \ldots, \mathbf{X}_n$
an observed $n$-tuple of $k$-dimensional elliptical vectors with
location ${\bolds\theta}$ and shape $\Vb$. Let $\mathbf{Z}_{i}:=
\mathbf{V}^{-1/2}(\mathbf{X}_i - {\bolds\theta})$ denote the
\textit{sphericized} version of $\mathbf{X}_i$ (throughout
$\mathbf{A}^{1/2}$, for a symmetric and positive definite matrix
$\mathbf{A}$, stands for the symmetric and positive definite root of
$\mathbf A$): the corresponding multivariate signs are defined as the unit
vectors $\mathbf{U}_i= \mathbf{U}_i({\bolds\theta}, \Vb):=
\mathbf{Z}_i/\Vert\mathbf{Z}_i\Vert$, while the ranks
$R^{(n)}_i=R^{(n)}_i({\bolds\theta}, \Vb)$ are those of the norms
$\Vert\mathbf{Z}_i\Vert$, $i=1,\ldots, n$. Our rank tests are based on
signed-rank covariance matrices of the form
\[
{\utSb}{}_{K}^{(n)}
:= \frac{1}{n} \sum_{i=1}^{n} K \biggl(\frac{ R^{(n)}_{i}}{n+1}
\biggr) \mathbf{U}_{i} \mathbf{U}_{i}\pr,
\]
where $K\dvtx(0,1)\to\R$ stands for some \textit{score function}, and
$\mathbf{U}_{i}=\mathbf{U}_i(\hat{\bolds\theta}, \hats{\Vb})$ and
$R^{(n)}_{i} =
R^{(n)}_i(\hat{\bolds\theta}, \hats{\Vb})$ are computed\vspace*{1pt} from
appropriate estimators
$\hat{\bolds\theta}$ and $ \hats{\Vb}$ of ${\bolds\theta}$ and~$\Vb$.
More precisely,\vspace*{-3pt}
for the testing problem (a), the rank-based test ${\utphi
}{}^{(n)}_{\betab;
K}$ rejects the null hypothesis ${\mathcal H}_{0}^{\betab}$ (at
asymptotic level $\alpha$) whenever
\[
{\utQ}{}_{K}^{(n)}
:=
\frac{nk(k+2)}{\mathcal{J}_k(K)} \sum_{j=2}^{k} \bigl(\tilde{\betab
}_{j}\pr{\utSb}{}_{K}^{(n)}\betab^{0} \bigr)^{2}
\]
exceeds the $\alpha$ upper-quantile of the chi-square distribution
with $(k-1)$ degrees of freedom; here, $\mathcal{J}_k(K)$ is a
standardizing constant and $\tilde{\betab}_j$ stands for a
constrained estimator of $\Vb$'s $j$th eigenvector; see (\ref
{choiceprelim}) for details.
As for problem (b), our rank tests ${\utphi}{}^{(n)}_{\Lamb; K}$ are
based on statistics of the form
\[
{\utT}{}_{K}^{(n)}:=
\biggl(\frac{nk(k+2)}{\mathcal{J}_k(K)} \biggr)^{1/2}
(a_{p,q}(\tilde{\Lamb}_\Vb))^{-1/2}
\mathbf{c}_{p,q}\pr\dvec\bigl(\tilde{\Lamb}_{\Vb}^{1/2}\hat
{\betab}\pr{\utSb}{}_{K}^{(n)}\hat{\betab}\tilde{\Lamb}_{\Vb
}^{1/2}\bigr),
\]
where $\tilde{\Lamb}_{\Vb}$ and $\hat{\betab}$ are adequate
estimators of $\Lamb_{\Vb}$ and $\betab$, respectively.
The null hypothesis ${\mathcal H}_{0}^{\Lamb
}$ is to be rejected at asymptotic level $\alpha$ whenever
${\utT}{}_{K}^{(n)}$ is smaller than the standard normal $\alpha$-quantile.

These tests are not just validity-robust, they also are efficient. For
any smooth radial density $f_1$, indeed, the score function $K=K_{f_1}$
(see Section~\ref{scores}) provides a signed-rank test which is \textit
{locally and asymptotically optimal} (\textit{locally and asymptotically
most stringent}, in the Le Cam sense) under radial density $f_1$. In
particular, when based on \textit{normal or van der Waerden} scores
$K=K_{\phi_1}:=\Psi_k^{-1}$, where \label{defphik} $\Psi_k$ denotes
the chi-square distribution function with $k$ degrees of freedom, our
rank tests achieve the same asymptotic performances as the optimal
Gaussian ones at the multinormal, while enjoying maximal validity
robustness, since no assumption is required on the underlying density
beyond ellipticity. Moreover, the asymptotic relative efficiencies
(AREs) under non-Gaussian densities of these van der Waerden tests are
uniformly larger than one with respect to their pseudo-Gaussian
parametric competitors; see Section~\ref{secare}. On all counts,
validity, robustness, and efficiency, our van der Waerden tests thus
perform uniformly better than the daily practice Anderson tests and
their pseudo-Gaussian extensions.

%s1.2 ###
\subsection{Local asymptotic normality for principal components}

The methodological tool we are using throughout is Le Cam's theory of
\textit{locally asymptotically normal} (LAN) \textit{experiments} [for
background reading on LAN, we refer to \citet{Lcam86}, \citet{CY00} or \citet{V98};
see also \citet{S85} or
\citet{R94}]. Although this powerful method has been used quite
successfully in inference problems for elliptical families
[Hallin and Paindaveine (\citeyear{HP02}, \citeyear{HP04}, \citeyear{HP05}, \citeyear{HP06a}),
\citet{HOP06} and \citet{HP08a} for location,
VARMA dependence, linear models, shape and scatter, resp.], it
has not been considered so far in problems involving eigenvectors and
eigenvalues, and, as a result, little is\vadjust{\goodbreak} known about optimality issues
in that context. The main reason, probably, is that the eigenvectors
$\betab$ and eigenvalues $\Lamb$ are complicated functions of the
covariance or scatter matrix $\Sigb$, with unpleasant identification
problems at possibly multiple eigenvalues. These special features of
eigenvectors and eigenvalues, as we shall see, make the LAN approach
more involved than in standard cases.

LAN (actually, ULAN) has been established, under appropriate regularity
assumptions on radial densities, in \citet{HP06a}, for
elliptical families when parametrized by a location vector ${\bolds
\theta}$ and a scatter matrix $\Sigb$ [more precisely, the vector
$\operatorname{vech}(\Sigb)$ resulting from stacking the upper diagonal elements of
$\Sigb$]. Recall, however, that LAN or ULAN are properties of the
parametrization of a family of distributions, not of the family itself.
Now, due to the complicated relation between $({\bolds\theta},
\operatorname{vech} \Sigb)$ and the quantities of interest $\Lamb$
and $\betab$, the $({\bolds\theta}, \operatorname{vech}(\Sigb
))$-parametrization is not convenient in the present context. Another
parametrization, involving location, scale, and shape eigenvalues and
eigenvectors is much preferable, as the hypotheses to be tested then
take simple forms. Therefore, we show (Lemma~\ref{LElemme}) how the ULAN result
of \citet{HP06a} carries over to this new
parametrization where, moreover, the information matrix, very
conveniently, happens to be block-diagonal---a structure that greatly
simplifies inference in the presence of nuisance parameters.
Unfortunately, this new parametrization, where $\betab$ ranges over
the set ${\mathcal SO}_k$ of $k\times k$ real orthogonal matrices with
determinant one, raises problems of another nature. The
subparameter $\operatorname{vec}(\betab)$ indeed ranges over $\operatorname{vec}({\mathcal SO}_k)$, a~nonlinear manifold of $\R^{k^2} $, yielding a \textit{curved} ULAN experiment.
% where
By a \textit{curved experiment}, we mean a parametric model indexed by a
$\ell$-dimensional parameter ranging over some nonlinear manifold of
$\R^\ell$, such as in \textit{curved} exponential families, for instance.
Under a $\operatorname{vec}(\betab)$-parametrization, the local experiments are not
the traditional \textit{Gaussian shifts} anymore, but \textit{curved}
Gaussian location ones, that is, Gaussian location models under which
the mean of a multinormal observation with specified covariance
structure ranges over a nonlinear manifold of $\R^\ell$, so that the
simple local asymptotic optimality results associated with local
Gaussian shifts no longer hold.
To the best of our knowledge, such experiments never have been
considered in the LAN literature.

A third parametrization, however, can be constructed from the fact that
$\betab$ is in ${\mathcal SO}_k$ if it can be expressed as the
exponential of a $k\times k$ skew-symmetric matrix~$\iotab$. Denoting
by $\operatorname{vech}^+(\iotab)$ the vector resulting from stacking the upper
off-diagonal elements of $\iotab$, this yields a parametrization
involving location, scale, shape eigenvalues and $\operatorname{vech}^+(\iotab)$; the
latter subparameter ranges freely over $\R^{k(k-1)/2}$, yielding a
well-behaved ULAN parametrization where local experiments converge to
the classical Gaussian shifts, thereby allowing for the classical
construction [\citet{Lcam86}, Section 11.9] of locally asymptotically
optimal tests. The trouble is that translating null hypotheses (a) and
(b) into the $\iotab$-space in practice seems unfeasible.

Three distinct ULAN structures are thus coexisting on the same families
of distributions:

(ULAN1) proved in \citet{HP06a} for the $({\bolds
\theta}, \operatorname{vech}(\Sigb))$-parametrization, serving as
the mother of all subsequent ones;

(ULAN2) for the location-scale-eigenvalues--eigenvectors
parametrization, where the null hypotheses of interest take simple
forms, but the local experiments happen to be \textit{curved} ones;

(ULAN3) for the location-scale-eigenvalues--skew symmetric matrix param\-
etrization, where everything is fine from a decision-theoretical point
of view, with, however, the major inconvenience that explicit solutions
cannot be obtained in terms of original parameters.

The main challenge of this paper was the delicate interplay between
these three structures. Basically, we are showing (Lemma~\ref{LElemme}) how ULAN
can be imported from the first parametrization, and (Section \ref
{curved}) optimality results from the third parametrization, both to
the second one. These results then are used in order to derive locally
asymptotically optimal Gaussian, pseudo-Gaussian and rank-based tests
for eigenvectors and eigenvalues of shape. This treatment we are giving
of curved ULAN experiments, to the best of our knowledge, is original, and
likely to apply in a variety of other contexts.

%s1.3 ###
\subsection{Outline of the paper} \label{outlinesec}
Section~\ref{sec2} contains, for easy reference, some basic notation and
fundamental assumptions to be used later on. The main ULAN result, of a
nonstandard curved nature, is established in Section~\ref{LANsection},
and its consequences for testing developed in Section \ref
{paramtests}. As explained in the \hyperref[intro]{Introduction}, optimality is imported
from an untractable parametrization involving skew-symmetric matrices.
This is elaborated, in some detail, in Section~\ref{curved}, where a
general result is derived, and in (Section~\ref{howeigenvectors}),
where that result is applied to the particular case of eigenvectors and
(Section~\ref{eigenvaluesprob}) eigenvalues of shape, under arbitrary
radial density $f_1$. Special attention is given, in Sections \ref
{Gaussbetab} and~\ref{Gausslamb}, to the Gaussian case ($f_1=\phi
_1$); in Sections~\ref{pseudovec} and~\ref{pseudoval}, those Gaussian
tests are extended to a pseudo-Gaussian context with finite
fourth-order moments. Then, in Section~\ref{ranktests}, rank-based
procedures, which do not require any moment assumptions, are
constructed: Section~\ref{rankHajek} provides a general asymptotic
representation result [Proposition~\ref{Hajek}(i)] in the H\' ajek
style; asymptotic normality, under the null as well as under local
alternatives, follows as a corollary [Proposition~\ref{Hajek}(ii)].
Based on these results, Sections~\ref{gsjdlr} and~\ref{rankeigvalues}
provide optimal rank-based tests for the eigenvector and eigenvalue
problems considered throughout; Sections~\ref{secare} and
\ref{secsimu} conclude with asymptotic relative efficiencies and
simulations. Technical proofs are concentrated in the \hyperref[app]{Appendix}.

The reader interested in inferential results and principal components
only (the form of the tests, their optimality properties and local
powers) may skip Sections~\ref{LANsection} and~\ref{paramtests},
which are devoted to curved LAN experiments, and concentrate on
Section~\ref{gausscase} for the ``parametric'' procedures, on
Section~\ref{ranktests} for the rank-based ones, on Sections~\ref
{secare} and~\ref{secsimu} for their asymptotic and finite-sample performances.

%s1.4 ###
\subsection{Notation} \label{notasec}
The following notation will be used throughout. For any $k \times
k$ matrix $\Ab=(A_{ij})$, write $\operatorname{vec}(\Ab)$ for the
$k^2$-dimensional vector obtained by stacking the columns of $\Ab$,
$\operatorname{vech}(\Ab)$ for the $[k(k+1)/2]$-dimensional vector
obtained by stacking the upper diagonal elements of those columns,
$\operatorname{vech}^+(\Ab)$ for the $[k(k-1)/2]$-dimensional vector
obtained by stacking\vspace*{1pt} the upper off-diagonal elements of the same, and
$\dvec(\Ab)=:(A_{11}, (\dvecrond(\Ab))\pr)\pr$
for the $k$-dimensional vector obtained by stacking the diagonal
elements of $\Ab;\dvecrond(\Ab)$ thus is $\dvec(\Ab)$ deprived
of its first component. Let $\Hb_{k}$ be the $k \times
k^{2}$ matrix such that $\Hb_{k} \vecop({\Ab})= \dvec(\Ab)$. Note
that we then have that $\Hb_{k}\pr\dvec(\Ab)= \vecop(\Ab)$ for
any $k\times k$ diagonal matrix $\Ab$, which implies that $\Hb_{k}\Hb
_{k}\pr=\mathbf{I}_k$.
Write $\operatorname{diag}(\mathbf{B}_{1}, \ldots, \mathbf{B}_{m})$ for the
block-diagonal matrix with blocks ${\Bb}_{1}, \ldots, {\Bb}_m$ and
$\mathbf{A}^{\otimes2}$ for the Kronecker
product $\mathbf{A}\otimes\mathbf{A}$. Finally, denoting by
$\mathbf{e}_\ell$ the $\ell$th vector in the canonical basis
of $\R^k$, write $\mathbf{K}_{k}:= \sum_{i,j =1}^{k}(\mathbf{e}_i\mathbf
{e}_j\pr)\otimes(\mathbf{e}_j\mathbf{e}_i\pr)$ for the $k^2\times k^2$
\textit{commutation matrix}.

%s2 ###
\section{Main assumptions}\label{sec2}
%s2.1 ###
\subsection{Elliptical densities}
\label{defelliptttt}
We throughout assume that the observations are elliptically
symmetric. More precisely, defining
\[
\mathcal{F} := \{ h\dvtx\R^+_0\to\R^+ \dvtx\mu_{k-1;h}<\infty\} ,
\]
where $\mu_{\ell;h}:=\int_0^{\infty} r^{\ell} h(r) \,dr$, and
\[
\mathcal{F}_1 := \biggl\{ h_1
\in\mathcal{F} \dvtx(\mu_{k-1;
h_1})^{-1}\int_0^{1}r^{k-1}h_1(r) \,d r =1/2 \biggr\},
\]
we denote by $\Xb_{1}^{(n)}, \ldots, \Xb^{(n)}_n$ an observed
$n$-tuple of mutually independent $k$-dimensional random vectors with
probability density function of the form
%
%e2.1 ###
%
\begin{equation} \label{density}
f(\mathbf{x}):=c_{k,f_1} \vert{\bolds\Sigma} \vert^{-1/2}
f_{1} \bigl(
\bigl(
(\mathbf{x} - {\bolds\theta})\pr{\bolds\Sigma}^{-1} (\mathbf{x} -
{\bolds\theta})
\bigr)^{1/2} \bigr),\qquad \mathbf{x}\in\R^{k},
\end{equation}
for some $k$-dimensional vector ${\bolds\theta}$ (\textit{location}),
some symmetric and positive definite $(k\times k)$ \textit{scatter}
matrix $\Sigb$, and some $f_1$ in the class $\mathcal{F}_1 $ of \textit
{standardized radial densities}; throughout, $\vert\Ab\vert$ stands
for the determinant of the square matrix $\Ab$.

Define the \textit{elliptical coordinates} of $\Xb_{i}^{(n)}$ as
%
%e2.2 ###
%
\begin{eqnarray} \label{Ud}
\mathbf{U}_{i}^{(n)}({\bolds\theta}, {\Sigb}) &:=& \frac{\Sigb^{-1/2}
(\Xb_{i}^{(n)} -{\bolds\theta})}{\| \Sigb^{-1/2} (\Xb
_{i}^{(n)} -{\bolds\theta})\|}\quad\mbox{and}\nonumber\\[-8pt]\\[-8pt]
d_{i}^{(n)}({\bolds\theta}, {\Sigb}) &:=& \bigl\| \Sigb^{-1/2}
\bigl(\Xb_{i}^{(n)} -{\bolds\theta}\bigr)\bigr\|.\nonumber
\end{eqnarray}
Under the assumption of ellipticity, the \textit{multivariate signs}
$\mathbf{U}_{i}^{(n)}({\bolds\theta}, {\Sigb})$,\break
$ i=1,\ldots, n $, are i.i.d. uniform over the unit sphere
in $\R^k$, and independent of the \textit{standardized elliptical
distances} $d^{(n)}_{i} ({\bolds\theta}, {\Sigb})$. Imposing~that
$f_1\in\mathcal{F}_1$\break
implies that the $d^{(n)}_{i}({\bolds\theta}, {\Sigb})$'s, which have common density $\tilde
{f}_{1k}(r) :=\break (\mu_{k-1;f_1})^{-1} r^{k-1} f_1 ( {r} )
I_{[r>0]}$, with distribution function $\tilde{F}_{1k}$, have median
one [$\tilde{F}_{1k}(1)=1/2$]---a constraint which identifies $\Sigb$
without requiring any moment assumptions [see \citet{HP06a}
for a discussion]. Under finite second-order
moments, the scatter matrix $\Sigb$ is proportional to the traditional
covariance
matrix $\Sigb_{\mathrm{cov}}$.% of the $\Xb_{i}\n$'s.

Special instances %of elliptical densities
are the $k$-variate
multinormal distribution, with radial density
$f_1(r)=\phi_1(r):=\exp(- a_k r^2/2)\label{ak}$, the $k$-variate
Student distributions, with radial densities (for $\nu\in\R^+_0$
degrees of freedom) $f_1(r)=f_{1,\nu}^t(r):=(1 + a_{k,\nu}
r^2/\nu)^{-(k+\nu)/2}$, and the $k$-variate power-exponential
distributions, with radial densities of the form $f_1(r)=
f_{1,\eta}^e (r):=\exp(- b_{k,\eta} r^{2\eta})$, $\eta\in\R^+_0$;
the positive constants $a_k$, \label{defnak} $a_{k,\nu}$, and
$b_{k,\eta}$ are such that $f_1\in\mathcal{F}_1$.

The derivation of locally and asymptotically optimal tests at standardized
radial density $f_1$ will be based on the \textit{uniform local and
asymptotic normality} (ULAN) of the model \textit{at given $f_1$}.
This ULAN property---the statement of which requires some further
preparation and is delayed to Section~\ref{LANsection}---only holds
under some further mild regularity conditions on $f_1$. More
precisely, we require $f_1$
to belong to the collection $\mathcal{F}_a$ of all absolutely
continuous densities in $\mathcal{F}_1$ for which, denoting by $\dot
f_1$ the
a.e. derivative of $f_1$ and letting
$\varphi_{f_1}:= -{\dot f_1}/f_1$, the integrals
%
%e2.3 ###
%
\begin{equation}\label{definfo}\qquad
\mathcal{I}_k(f_1) := \int_{0}^1 \varphi_{f_1}^2(r) \tilde
{f}_{1k}(r) \,dr
%$$
\quad\mbox{and}\quad
\mathcal{J}_k(f_1):=\int_{0}^1
r^2 \varphi_{f_1}^2(r)\tilde{f}_{1k}(r) \,dr
\end{equation}
are finite. The quantities $\mathcal{I}_k(f_1)$ and $\mathcal
{J}_k(f_1)$ play the roles of \textit{radial Fisher information for
location} and \textit{radial Fisher information for shape/scale},
respectively. Slightly less stringent assumptions, involving
derivatives in the sense of distributions, can be found in
\citet{HP06a}, where we refer to for details. The intersection of
$\mathcal{F}_a$ and $\mathcal{F}^4_1:=\{f_1\!\in\!\mathcal{F}_1 \dvtx\int
_0^\infty r^4\tilde{f}_{1k}(r) \,dr<\infty\}$ will be denoted as
$\mathcal{F}_a^4$.

%s2.2 ###
\subsection{Score functions} \label{scores}
The various \textit{score functions} $K$ appearing in the rank-based
statistics to be introduced in Section~\ref{ranktests}
%(\ref{rankcovv})
will be assumed to satisfy a few regularity assumptions which we are
listing here for convenience.
\renewcommand{\theassumption}{$(\mathrm{S})$}
\begin{assumption}\label{assuS}
The score function $K \dvtx(0,1)
\rightarrow\rr$
(S1) is continuous and square-integrable,
(S2) can be expressed as the difference of two monotone increasing
functions, and
(S3) satisfies $\int_0^1 K(u) \,du=k$.
\end{assumption}

Assumption (S3) is a normalization constraint that is
automatically satisfied by the score functions
$K(u)=K_{f_1}(u):= \varphi_{f_1} ( {\tilde
F}{}^{-1}_{1k}(u)) {\tilde F}{}^{-1}_{1k}(u)$\vadjust{\goodbreak} associated with any radial
density $f_1\in\mathcal{F}_a$ (at which ULAN
holds); see Section~\ref{LANsection}. For score functions $K,K_1,K_2$
satisfying Assumption~\ref{assuS}, let [throughout, $U$ stands for a random
variable uniformly
distributed over $(0,1)$]
%
%e2.4 ###
%
\begin{equation}\label{infoK}
\mathcal{J}_k(K_1,K_2) :=\mathrm{E}[K_1(U)K_2(U)],\qquad
\mathcal{J}_k(K):=\mathcal{J}_k(K,K)
\end{equation}
and
%
%e2.5 ###
%
\begin{equation}\label{infoKf}
\mathcal{J}_k(K,f_1):= \mathcal{J}_k(K,K_{f_1});
\end{equation}
with this notation, $\mathcal{J}_k(f_1) = \mathcal
{J}_k(K_{f_1},K_{f_1})$.

The \textit{power} score functions $K_a(u) := k(a+1) u^a$ ($a\geq0$), with
$\mathcal{J}_k(K_a)=k^2 (a+1)^2/(2a+1)$, provide
some traditional score functions satisfying Assumption~\ref{assuS}: the sign,
Wilcoxon, and Spearman scores
are obtained for $a=0$, $a=1$ and $a=2$, respectively.
As for the score functions of the form $K_{f_1}$, an important
particular case is that of van der Waerden or \textit{normal} scores,
obtained for $f_1=\phi_1$. Then
%
%e2.6 ###
%
\begin{equation}\label{normalscores}
K_{\phi_1}(u)=\Psi_k^{-1}(u) \quad\mbox{and}\quad \mathcal
{J}_k(\phi_1)=k(k+2),
\end{equation}
where $\Psi_k$ was defined in page \pageref{defphik}.
Similarly, Student densities $f_1=f^t_{1,\nu}$ (with $\nu$ degrees of
freedom) yield the scores
\[
K_{f_{1,\nu}^t}(u)=\frac{k(k+\nu)G_{k,\nu}^{-1}(u)}{\nu+kG_{k,\nu
}^{-1}(u)}
\]
and
\[
\mathcal{J}_k(f_{1,\nu}^t )=\frac
{k(k+2)(k+\nu)}{k+\nu+2},
\]
where $G_{k,\nu}$ stands for the Fisher--Snedecor distribution
function with $k$ and $\nu$ degrees of freedom.

%s3 ###
\section{Uniform local asymptotic normality (ULAN) and curved Gaussian
location local experiments}\label{LANsection}

%s3.1 ###
\subsection{Semiparametric modeling of elliptical families}
Consider an i.i.d. $n$-tuple $\mathbf{X}_1^{(n)},\ldots, \mathbf
{X}_n^{(n)}$ with elliptical density~(\ref{density})
characterized by
$\tetb$, $\Sigb$, and $f_1\dvtx({\bolds\theta},\Sigb)$ or, if a
vector is to be preferred, $({\bolds\theta}\pr, (\operatorname
{vech} \Sigb)\pr)\pr$, provides a perfectly valid parametrization
of the elliptical family with standardized radial density $f_1$.
However, in the problems we are considering in this paper, it will be
convenient to have eigenvalues and eigenvectors appearing explicitly in
the vector of parameters. Decompose therefore the scatter matrix $\Sigb
$ into $\Sigb=\sigma^2\mathbf{V}= \betab\Lamb_{\Sigb} \betab\pr=\betab
\sigma^2 \Lamb_{\Vb}
\betab\pr$, where $\sigma\in\rr^+_0$ is a
\textit{scale parameter} (equivariant under multiplication by a
positive constant), and $\mathbf V$ a \textit{shape matrix} (invariant
under multiplication by a positive constant) with eigenvalues $\Lamb
_{\Vb}=\operatorname{diag}(\lambda_{1;\Vb},\ldots, \lambda_{k;\Vb
})=\sigma^{-2}\operatorname{diag}(\lambda_{1;\Sigb},\ldots,
\lambda_{k;\Sigb}) = \sigma^{-2}\Lamb_{\Sigb}$; $\betab$ is an
element of the so-called \textit{special orthogonal group}\vadjust{\goodbreak}
${\mathcal SO}_{k}:= \{ \mathbf{O} | \mathbf{O}\pr\mathbf{O}= \mathbf{I}_k$,
$|\mathbf{O}|=1\}$
diagonalizing both $\Sigb$ and $\Vb$, the columns $\betab_1,\ldots,
\betab_k$ of which are the eigenvectors (common to $\Sigb$ and $\Vb
$) we are interested in.

Such decomposition of scatter into scale and shape can be achieved in
various ways. Here, we adopt the determinant-based definition of scale
\[
\sigma:=\vert\Sigb\vert^{1/2k}=
\prod_{j=1}^{k} \lambda_{j;\Sigb} ^{1/2k}\qquad \mbox{hence }
\Vb:= \Sigb/ \sigma^2 = \Sigb/ \vert\Sigb\vert^{1/k},
\]
which implies that $\vert\Vb\vert=\prod_{j=1}^{k} \lambda_{j;\Vb
}=1$. As shown by \citet{P08}, this choice indeed is the only
one for which the information matrix for scale and shape is
block-diagonal, which greatly simplifies inference. The parametric
families of elliptical distributions with specified standardized radial
density $f_1$ then are indexed by the $L= k(k+2)$-dimensional parameter
\[
\varthetab:= ({\bolds\theta}^{ \prime}, \sigma^{2},
(\dvecrond
\Lamb_{\Vb})\pr, (\vecop\betab)\pr)\pr=:(\varthetab_\I\pr
,\vartheta_\II,\varthetab_\III\pr,\varthetab_\IV\pr)\pr,
\]
where $\dvecrond(\Lamb_{\Vb})=(\lambda_{2;\Vb}, \ldots,
\lambda_{k;\Vb})\pr$ since $\lambda_{1;\Vb} =\prod_{j=2}^{k}
\lambda_{j;\Vb}^{-1}$.\vspace*{2pt}

This $\varthetab$-parametrization however requires a fully identified
$k$-tuple of eigenvectors, which places the following restriction on
the eigenvalues $\Lamb_{\Vb}$.
\renewcommand{\theassumption}{$(\mathrm{A})$}
\begin{assumption}\label{assuA}
The eigenvalues $\lam_{j;\Vb}$ of the shape
matrix $\Vb$ are all distinct, that is, since $\Sigb$ (hence also
$\Vb$) is positive definite, $\lam_{1;\Vb}
> \lam_{2; \Vb}
> \cdots> \lam_{k;\Vb}>0$.
\end{assumption}

Denote by $\mathrm{P}^{(n)}_{\varthetab; f_1}$
the joint
distribution of $\mathbf{X}_1^{(n)},\ldots, \mathbf{X}_n^{(n)}$ under parameter
value $\varthetab$ and standardized radial density $f_1\in\mathcal
{F}_1$; the parameter space [the definition of which includes
Assumption~\ref{assuA}] then
is
\[
\Thetab: =\rr^k\times\rr^+_0\times\mathcal{C}^{k-1}\times
\operatorname{vec}({\mathcal SO}_{k}),
\]
where $\mathcal{C}^{k-1}$ is the open cone of $(\rr^+_0)^{k-1}$ with
strictly ordered (from largest to smallest) coordinates.

Since $\operatorname{vec}({\mathcal SO}_{k})$ is a nonlinear manifold of $\rr^{k^2}
$: the $\operatorname{vec}(\betab)$-parametri\-zed experiments are \textit{curved
experiments}, in which the standard methods [see Section 11.9 of
\citet{Lcam86}] for constructing locally asymptotically optimal tests do
not apply. It is well known, however [see, e.g., \citet{KG89}],
that any element $\betab$ of ${\mathcal SO}_{k}$ can be expressed as
the exponential $\exp(\iotab)$ of a $k\times k$ \textit{skew-symmetric
matrix}~$\iotab$, itself characterized by the $k(k-1)/2$-vector
$\operatorname{vech}^+(\iotab)$ of its upper off-diagonal elements.
The differentiable mapping $\hbar\dvtx\operatorname{vech}^+(\iotab)
\mapsto\hbar(\operatorname{vech}^+(\iotab)):= \operatorname
{vec}(\exp(\iotab))$ from $\rr^{k(k-1)/2}$ to $ {\mathcal SO}_{k}$
is one-to-one, so that $\operatorname{vech}^+(\iotab)\in\rr
^{k(k-1)/2}$ also can be used as a parametrization instead of
$\operatorname{vec}(\betab)\in\operatorname{vec}({\mathcal SO}_{k})$. Both
parametrizations yield uniform local asymptotic normality (ULAN).
Unlike the $\operatorname{vec}(\betab)$-parametrized one, the $\operatorname
{vech}^+(\iotab)$-parametrized experiment is not curved, as
$\operatorname{vech}^+(\iotab)$ freely\vadjust{\goodbreak} ranges over $\rr^{k(k-1)/2}$,
so that the standard methods for constructing locally asymptotically
optimal tests apply---which is not the case with curved experiments. On
the other hand, neither the $\operatorname{vech}^+(\iotab)$-part of
the central sequence, nor the image in the $\operatorname
{vech}^+(\iotab)$-space of the null hypothesis $\mathcal{H}_0^{\betab
}$ yield tractable %explicit
forms. Therefore, we rather state ULAN for the curved
$\operatorname{vec}(\betab)$-parametrization. Then (Section~\ref{curved}), we develop a general
theory of locally asymptotically optimal tests for differentiable
hypotheses in curved ULAN experiments.

Without Assumption~\ref{assuA}, the $\varthetab$-parametrization is not valid,
and cannot enjoy LAN nor ULAN; optimality properties (of a local and
asymptotic nature) then cannot be obtained. As far as validity issues
(irrespective of optimality properties) are considered, however, this
assumption can be weakened. If the null hypothesis $\mathcal{H}^\betab
_0$ is to make any sense, the first eigenvector $\betab_1$ clearly
should be identifiable, but not necessarily the remaining ones. The
following assumption on the $\lambda_{j;\Vb}$'s, under which $\betab
_2,\ldots,\betab_k$ need not be identified, is thus minimal in that case.
\renewcommand{\theassumption}{$(\mathrm{A}\pr_1)$}
\begin{assumption}\label{assuApr1}
The eigenvalues of the shape
matrix $\Vb$ are such that
$\lam_{1;\Vb}
> \lam_{2;\Vb}
\geq\cdots\geq\lam_{k;\Vb}>0$.
\end{assumption}

Under Assumption~\ref{assuApr1}, $\Thetab$ is broadened into a larger
parameter space $\Thetab\pr_1$, which does not provide a valid
parametrization anymore, and for which the ULAN property of
Proposition~\ref{LAN} below no longer holds. As we shall see, all the
tests we are proposing for $\mathcal{H}^\betab_0$ nevertheless
remains valid under the extended null hypothesis $ {\mathcal
H}_{0;1}^{\betab\prime}$ resulting from weakening~\ref{assuA} into
\ref{assuApr1}. Note that, in case the null hypothesis is dealing with $\betab
_q$ instead of $\betab_1$, the appropriate weakening of Assumption~\ref{assuA}
is the following.
\renewcommand{\theassumption}{$(\mathrm{A}\pr_q)$}
\begin{assumption}\label{assuAprq}
The eigenvalues of the shape
matrix $\Vb$ are such that
$\lam_{1;\Vb}\geq\cdots\geq\lambda_{q-1;\Vb}>
\lam_{q;\Vb}>\lambda_{q+1;\Vb}
\geq\cdots\geq\lam_{k;\Vb}>0$.
\end{assumption}

This yields enlarged parameter space $\Thetab\pr_q$ and null
hypothesis ${\mathcal H}_{0;q}^{\betab\prime}$.

Similarly, the null hypothesis $ {\mathcal H}_{0}^{\Lamb}$ requires
the identifiability of the groups of $q$ largest (hence $k-q$ smallest)
eigenvalues; within each group, however, eigenvalues may coincide,
yielding the following assumption.
\renewcommand{\theassumption}{$(\mathrm{A}^{\prime\prime}_q)$}
\begin{assumption}\label{assuAq}
The eigenvalues of the shape
matrix $\Vb$ are such that
$\lam_{1;\Vb}\geq\cdots\geq\lambda_{q-1;\Vb}\geq
\lam_{q;\Vb}>\lambda_{q+1;\Vb}
\geq\cdots\geq\lam_{k;\Vb}>0$.
\end{assumption}

This yields enlarged parameter space
$\Thetab^{\prime\prime}_q$ and null hypothesis $ {\mathcal
H}_{0;q}^{\Lamb\prime\prime}$, say. As we shall see, the tests we
are proposing for $ {\mathcal H}_{0}^{\Lamb}$ remain valid under $
{\mathcal H}_{0;q}^{\Lamb\prime\prime}$.

%s3.2 ###
\subsection{Curved ULAN experiments}
\label{curvLAN}
Uniform local asymptotic normality\break
(ULAN) for the parametric families or \textit{experiments}
$\mathcal{P}^{(n)}_{f_1}:=\{\mathrm{P}^{(n)}_{\varthetab; f_1} \dvtx
\varthetab\in\Thetab
\}
$, with classical root-$n$ rate, is the main technical tool of this\vadjust{\goodbreak}
paper. For any $\varthetab:= ({\bolds\theta}\pr, \sigma^{2}$,
$(\dvecrond\Lamb_{\Vb})\pr, (\vecop\betab)\pr)\pr$ $\in
\Thetab$, a~\textit{local alternative} is a sequence $\varthetab
^{(n)}\in\Thetab$ such that $(\varthetab^{(n)}- \varthetab) $ is $
O(n^{-1/2})$. For any such $\varthetab^{(n)}$, consider a further
sequence $\varthetab^{(n)} +n^{-1/2} \taub^{(n)}$, with
$\taub^{(n)} =
%& =&
( (\taub^{\I(n)})\pr, \tau^{\II(n)} , (\taub^{\III(n)})\pr,\break
(\taub^{\IV(n)})\pr)\pr$ such that $\sup_{n}
\taub\npr\taub^{(n)} < \infty$ and $\varthetab^{(n)} +n^{-1/2}
\taub^{(n)} \in\Thetab$ for all $n$. Note that such $\taub^{(n)}$
exist: $\taub^{\I(n)}$ can be any bounded sequence of $\rr^k$, $\tau
^{\II(n)}$ any bounded sequence with $\tau^{\II(n)}>-n^{1/2}\sigma
^{2(n)}$, $\taub^{\III(n)}$ any bounded sequence of real
$(k-1)$-tuples $(\tau^{\III(n)}_1,\ldots, \tau^{\III(n)}_{k-1})$
such that
\begin{eqnarray*}
0 &<& \lambda^{(n)}_{k;\Vb} + n^{-1/2}\tau^{\III(n)}_{k-1} < \cdots\\
&<&\lambda^{(n)}_{3;\Vb} + n^{-1/2} \tau^{\III(n)}_2
< \lambda^{(n)}_{2;\Vb} + n^{-1/2}\tau^{\III(n)}_1 \\
&<&
\prod_{j=2}^k\bigl(\lambda^{(n)}_{j;\Vb} + n^{-1/2}\tau^{\III
(n)}_{j-1}\bigr)^{-1} ,
\end{eqnarray*}
which ensures that the perturbed eigenvalues $\lambda_{j;\Vb
}^{(n)}+n^{-1/2}\ell^{(n)}_j$, %$i=1,\ldots, k$,
with
%
%e3.1 ###
%
\begin{eqnarray}\label{ell1}
\ell^{(n)}_1 :\!&=& n^{1/2} \Biggl(\prod_{j=2}^k\bigl(\lambda^{(n)}_{j;\Vb}
+ n^{-1/2}\tau^{\III(n)}_{j-1}\bigr)^{-1} - \lambda^{(n)}_{1;\Vb}
\Biggr) \nonumber\\[-8pt]\\[-8pt]
&=&
-\lambda^{(n)}_{1;{\Vb}} \sum_{j=2}^k \bigl(\lambda^{(n)}_{j;\Vb
}\bigr)^{-1}\tau^{\III(n)}_{j-1} + O(n^{-1/2})\nonumber
\end{eqnarray}
and $(\ell^{(n)}_2,\ldots,\ell^{(n)}_k):=\taub^{\III(n)}$, still
satisfy Assumption~\ref{assuA} and yield determinant value one.
Writing
${\bolds\ell}^{(n)}$ for the diagonal $k\times k$ matrix with
diagonal elements $\ell^{(n)}_1,\ldots, \ell^{(n)}_k$, we then have
\begin{eqnarray*}
\operatorname{tr}\bigl(\bigl(\Lamb_{\Vb}^{(n)}\bigr)^{-1} {\bolds\ell}^{(n)}\bigr)&=&
\bigl(\lambda^{(n)}_{1;{\Vb}}\bigr)^{-1} \Biggl[
-\lambda^{(n)}_{1;{\Vb}} \sum_{j=2}^k \bigl(\lambda^{(n)}_{j;{\Vb
}}\bigr)^{-1}\tau^{\III(n)}_{j-1} + O(n^{-1/2})
\Biggr] \\
&&{}
+ \sum_{j=2}^k \bigl(\lambda^{(n)}_{j;{\Vb}}\bigr)^{-1} \tau^{\III(n)}_{j-1}
\\
&=& O(n^{-1/2}).
\end{eqnarray*}

Finally, denote by $\mathbf{M}_{k}\pr(\lambda_2,\ldots,\lambda_k)=
( -\lambda_1(\lambda_2^{-1},\ldots, \lambda_k^{-1})\pr
%_i^{-1}
{\,}\vdots{\,}\mathbf{I}_{k-1} )\pr$ the value at $(\lambda_2,\ldots
, \lambda_k)$ of the Jacobian matrix of
\[
( \lambda_2, \ldots, \lambda_k)
\mapsto\Biggl(\lambda_1:=\prod_{j=2}^k\lambda_j^{-1},\lambda_2,\ldots,
\lambda_k\Biggr)%\in\R^{k}$
.
\]
Letting $\Lamb:= \operatorname{diag}(\lambda_1,\lambda_2,\ldots,
\lambda_k)$, we have $\mathbf{M}_{k} \pr(\lambda_2,\ldots,\lambda
_k)\dvecrond(\mathbf{l})= \dvec(\mathbf{l})$ for any $k \times k$ real
matrix $\mathbf{l}$ such that $\operatorname{tr}(\Lamb^{-1} \mathbf
{l})=0$. Indeed,
\begin{eqnarray*}%\label{Jacobienne}
\mathbf{M}_{k} \pr(\lambda_2,\ldots,\lambda_k) \dvecrond(\mathbf{l})
&=& \bigl( -\lambda_1(\lambda_2^{-1},\ldots, \lambda_k^{-1})\pr
{\,}\vdots{\,}\mathbf{I}_{k-1} \bigr)\pr(\dvecrond\mathbf{l}) \\
&=& \bigl( - \lambda_{1 } \bigl(\operatorname{tr}(\Lamb^{-1} \mathbf{l})- (\lambda
_{1})^{-1} \mathbf{l}_{11}\bigr) {\,}\vdots{\,} (\dvecrond\mathbf{l})'
\bigr)'\\
&=& \dvec(\mathbf{l}),
\end{eqnarray*}
an identity that will be used later on for $\mathbf{M}_{k}^{{\bolds
\Lambda}_\Vb}:=\mathbf{M}_{k}(\dvecrond({\bolds\Lambda}_{\Vb}))$.

The problem is slightly more delicate for $\taub^{\IV(n)}$, which
must be such that
$\vecop({\betab}^{(n)}) + n^{-1/2} \taub^{\IV(n)}$ remains in
$\vecop({\mathcal SO}_{k})$. That is, $\taub^{\IV(n)}$ must be of
the form $ \taub^{\IV(n)}=\vecop(\mathbf{b}^{(n)})$, with
%
%e3.2 ###
%
\begin{eqnarray}\label{antisym}
\mathbf{0} &=& \bigl({\betab}^{(n)}+ n^{-1/2} \mathbf{b}^{(n)} \bigr)\pr
\bigl({\betab}^{(n)}+ n^{-1/2} \mathbf{b}^{(n)} \bigr)-\mathbf{I}_k
\nonumber\\[-8pt]\\[-8pt]
&=& n^{-1/2} \bigl( \betab^{(n)\prime} \mathbf{b}^{(n)}+ \mathbf
{b}^{(n)\prime} \betab^{(n)} \bigr) +
n^{-1}\mathbf{b}^{(n)\prime}\mathbf{b}^{(n)}.\nonumber
\end{eqnarray}
That is, $\betab^{(n)\prime} \mathbf{b}^{(n)}+ n^{-1/2}\mathbf
{b}^{(n)\prime}\mathbf{b}^{(n)}/2$ should be skew-symmetric. Such local
perturbations admit an intuitive interpretation: we have indeed
\[
{\betab}^{(n)}+ n^{-1/2} \mathbf{b}^{(n)}= \betab^{(n)}\betab
^{(n)\prime} \bigl({\betab}^{(n)}+ n^{-1/2} \mathbf{b}^{(n)}\bigr) = \betab
^{(n)}\bigl(\mathbf{I}_k + n^{-1/2} \betab^{(n)\prime} \mathbf{b}^{(n)}\bigr)
\]
an expression in which $\mathbf{I}_k + n^{-1/2}\betab^{(n)\prime}
\mathbf{b}^{(n)}$, up to a $O(n^{-1})$ quantity, coincides with the first-order
approximation of the exponential of a skew-symmetric matrix, and
therefore can be interpreted as an infinitesimal rotation.
Identity~(\ref{antisym}) provides a characterization of ${\mathcal
SO}_k$ in the vicinity of $\betab^{(n)}$. The tangent space [in $\rr
^{k^2} $, at $\operatorname{vec}(\betab)$] to $\operatorname{vec}({\mathcal SO}_k)$ is obtained by
linearizing~(\ref{antisym}). More precisely, this tangent space is of the
form
%
%e3.3 ###
%
\begin{eqnarray} \label{tangentbetab}
&&\{
\operatorname{vec}(\betab+ \mathbf{b})
\vert
\operatorname{vec}(\mathbf{b})\in\rr^{k^2}
\mbox{ and }
\betab\pr\mathbf{b} + \mathbf{b}\pr\betab= \mathbf{0}
\}
\nonumber\\[-8pt]\\[-8pt]
&&\qquad
= \{
\operatorname{vec}(\betab+ \mathbf{b})
\vert
\operatorname{vec}(\mathbf{b})\in\rr^{k^2}
\mbox{ and }
\betab\pr\mathbf{b} \mbox{ skew-symmetric}
\}.\nonumber
\end{eqnarray}

We then have the following result (see the \hyperref[app]{Appendix} for the proof).
\begin{Prop}\label{LAN}
The experiment $\mathcal{ P}^{(n)}_{f_1} := \{
\mathrm{P}^{(n)}_{\varthetab;f_1} \vert
{\bolds\vartheta}\in\Thetab\}$ is ULAN, with central sequence
$\Deltab^{(n)}_{{\bolds\vartheta};f_1}:= ( %\begin{array}{c}
\Deltab^{\I\prime} _{{\bolds\vartheta};f_1},
\Delta^{\II} _{{\bolds\vartheta};f_1} ,
\Deltab^{\III\prime}_{{\bolds\vartheta};f_1},
\Deltab^{\IV\prime} _{{\bolds\vartheta};f_1}
)\pr$,
where [with $d_i:= d_i^{(n)}({\bolds\theta}, \Vb)$ and $\mathbf
{U}_i:=\mathbf{U}^{(n)}_i({\bolds\theta}, \Vb)$ as defined in (\ref
{Ud}), and letting $\mathbf{M}_{k}^{{\bolds\Lambda}_\Vb}:=\mathbf
{M}_{k}(\dvecrond{\bolds\Lambda}_\Vb)$],
\begin{eqnarray*}
\Deltab^{\I} _{{\bolds\vartheta};f_1} &:=&
\frac{1}{\sqrt{n}\sigma}
\sum_{i=1}^{n} \varphi_{f_1} \biggl(\frac{d_{i}}{\sigma} \biggr)
{\Vb}^{-1/2}\mathbf{U}_{i},
\\
\Delta^{\II} _{{\bolds\vartheta};f_1} &:=&
\frac{1}{2\sqrt{n} \sigma^{2}}
\sum_{i=1}^{n} \biggl( \varphi_{f_1}
\biggl(\frac{d_{i}}{\sigma} \biggr) \frac{d_{i}}{\sigma}
-k \biggr),
\\
\Deltab^{\III} _{{\bolds\vartheta};f_1} &:=&
\frac{1}{2\sqrt{n}} \mathbf{M}_{k}^{{\bolds\Lambda}_\Vb} \mathbf{H}_{k} (
\Lamb_{\Vb}^{-1/2}\betab\pr)^{\otimes2}
\sum_{i=1}^{n} \vecop\biggl( \varphi_{f_1}
\biggl(\frac{d_{i}}{\sigma} \biggr) \frac{d_{i}}{\sigma} \mathbf{U}_{i}\mathbf{U}\pr
_{i} \biggr)
\end{eqnarray*}
and
\[
\Deltab^{\IV} _{{\bolds\vartheta};f_1} :=
\frac{1}{2\sqrt{n}}
\mathbf{G}_{k}^{\betab}\mathbf{L}^{\betab, \Lamb_{\Vb}}_{k} ( {\Vb
}^{\otimes2} )^{ -1/2}
\sum_{i=1}^{n} \vecop\biggl( \varphi_{f_1}
\biggl(\frac{d_{i}}{\sigma} \biggr) \frac{d_{i}}{\sigma} \mathbf{U}_{i}\mathbf{U}\pr
_{i} \biggr),
\]
with
${\Gb}_{k}^{\betab}:=({\Gb}_{k;12}^{\betab} {\Gb
}_{k;13}^{\betab} \cdots
{\Gb}_{k;(k-1)k}^{\betab}), {\Gb}_{k;jh}^{\betab
}:=\mathbf{e}_{j} \otimes
{\betab}_{h}-\mathbf{e}_{h} \otimes{\betab}_{j}$
and
\[
\mathbf{L}^{\betab, \Lamb_{\Vb}}_{k}:=\bigl(\mathbf{L}^{\betab, \Lamb_{\Vb
}}_{k;12} \mathbf{L}^{\betab, \Lamb_{\Vb}}_{k;13} \cdots\mathbf
{L}^{\betab, \Lamb_{\Vb}}_{k;(k-1)k}\bigr)^{\prime},\qquad
\mathbf{L}^{\betab, \Lamb_{\Vb}}_{k;jh}:=(\lambda_{h;\Vb}-\lambda
_{j;\Vb})(\betab_{h} \otimes
\betab_{j}),
\]
and with block-diagonal information matrix
%
%e3.4 ###
%
\begin{equation}\label{Gamb}
\Gamb_{{\bolds\vartheta};f_1}=
\operatorname{diag}(
\Gamb^{\I}_{{\bolds\vartheta};f_1},
\Gamma^{\II}_{{\bolds
\vartheta};f_1},\Gamb^{\III}_{{\bolds\vartheta};f_1},\Gamb^{\IV
}_{{\bolds\vartheta};f_1}),
\end{equation}
where, defining $\mathbf{D}_k(\Lamb_\Vb):= \frac{1}{4} \mathbf
{M}_{k}^{{\bolds\Lambda}_\Vb} \mathbf{H}_{k} [
\mathbf{I}_{k^2} + \mathbf{K}_k ]
(\Lamb_{\Vb}^{-1})^{\otimes2} \mathbf{H}_{k}\pr
(\mathbf{M}_{k}^{{\bolds\Lambda}_\Vb})\pr$,
\[
\Gamb^{\I}_{{\bolds\vartheta};f_1}= \frac{\ikf}{k \sigma
^{2}}{\Vb}^{-1},\qquad %$$
%
%$$
\Gamma^{\II}_{{\bolds\vartheta};f_1}=\frac{\mathcal
{J}_k(f_1)-k^{2}}{4 \sigma^{4}},\qquad
\Gamb^{\III}_{{\bolds\vartheta};f_1}=\frac{\mathcal
{J}_k(f_1)}{k(k+2)}\mathbf{D}_k(\Lamb_\Vb)
\]
and
\[
\Gamb^{\IV}_{{\bolds\vartheta};f_1} := \frac{1}{4}\frac
{\mathcal{J}_k(f_1)}{k(k+2)} \mathbf{G}_{k}^{\betab} \operatorname
{diag}\bigl(\nu_{12}^{-1}
,\nu_{13}^{-1}, \ldots, \nu_{(k-1)k}^{-1}\bigr) (\mathbf{G}_{k}^{\betab
})\pr,
\]
where
$\nu_{jh}:= \lambda_{j;\Vb}
\lambda_{h;\Vb}/(\lambda_{j;\Vb}-\lambda_{h;\Vb})^{2}$. More
precisely, for any local alternative $\varthetab^{(n)}$
and any bounded sequence $\taub^{(n)}$ such that $\varthetab^{(n)}+
n^{-1/2}\taub^{(n)}\in\Thetab$, we have, under $\mathrm{P}^{(n)}
_{ \varthetab^{(n)};f_1}$,
\begin{eqnarray*}
\Lambda^{(n)}_{\varthetab^{(n)}+ n^{-1/2}\taub^{(n)}/ \varthetab^{(n)};f_1}
:\!&=&\log\bigl(d\mathrm{P}^{(n)}_{\varthetab^{(n)}+
n^{-1/2}\taub^{(n)};f_1 }/ d\mathrm{P}^{(n)}_{\varthetab^{(n)};f_1 } \bigr)
\\
&=& \bigl(\taub^{(n)}\bigr)\pr\Deltab^{(n)}_{\varthetab^{(n)};f_1 }
-\tfrac{1}{2}\bigl(\taub^{(n)}\bigr)\pr\Gamb_{\varthetab;f_1}\taub^{(n)}+
o_{\mathrm P}(1)
\end{eqnarray*}
and $ \Deltab_{\varthetab^{(n)};f_1 } \stackrel{\mathcal
{L}}{\longrightarrow} \mathcal{N}( \mathbf{0}, \Gamb_{\varthetab;f_1} )
$, as $\ny$.
\end{Prop}

The block-diagonal structure of the information matrix ${\Gamb
}_{\varthetab; f_{1}}$ implies that inference on $\betab$ (resp.,
$\Lamb_\Vb$) can be conducted under unspecified ${\bolds\theta}$,
$\sigma$ and ${\bolds\Lambda_{\mathbf{V}}}$ (resp., $\betab$) as if the latter
were known, at no asymptotic cost. The orthogonality between the
eigenvalue and eigenvector parts of the central sequence is structural,
while that between the eigenvalue and eigenvector parts on one hand and
the scale parameter part on the other hand is entirely due to the
determinant-based parametrization of scale [see \citet{HP06b} or \citet{P08}].
Note that ${\Gamb_{\varthetab
; f_{1}}^{\IV}}$, with rank $k(k-1)/2 <k^2$, is not invertible.

%s3.3 ###
\subsection{Locally asymptotically optimal tests for differentiable
hypotheses in curved ULAN experiments}\label{curved}
Before addressing testing problems involving eigenvalues and
eigenvectors, we need a general theory for locally asymptotically
optimal tests in curved ULAN experiments, which we are developing in
this section.

Consider a ULAN sequence of experiments $\{\mathrm{P}^{(n)}_{\bolds\xi}
\dvtx{\bolds\xi} \in\Xib\}$, where $ \Xib$ is an open subset of $
\mathbb{R}^m$, with central sequence $\Deltab_{\bolds\xi}$ and
information $\Gamb_{\bolds\xi}$. For the simplicity of exposition,
assume that $\Gamb_{\bolds\xi}$ for any ${\bolds\xi}$ has full
rank $m$. Let $\hbar\dvtx\Xib\to\mathbb{R}^p$, $p\geq m$, be a
continuously differentiable mapping such that the Jacobian matrix
$D\hbar({\bolds\xi})$ has full rank $m$ for all ${\bolds\xi}$, and
consider the experiments $\{\mathrm{P}^{(n)}_\varthetab\dvtx\varthetab
\in
\Thetab:=\hbar(\Xib)\}$, where, with a slight abuse of
notation, $\mathrm{P}^{(n)}_{\varthetab}:=\mathrm{P}^{(n)}_{\bolds\xi}$
for $\varthetab=\hbar({\bolds\xi})$. This sequence also is ULAN,
with central sequence $\Deltab_\varthetab$ and information matrix
$\Gamb_\varthetab$ such that [see~(\ref{LAQvartheta2}) and the proof
of Lemma~\ref{LElemme}], at $\varthetab=\hbar({\bolds\xi})$, and
[up to $o^{(n)}_{\mathrm{P}_\varthetab}(1)$'s which, for simplicity, we
omit here] $\Deltab_{\bolds\xi}=D\hbar\pr({\bolds\xi})\Deltab
_\varthetab$ and $\Gamb_{\bolds\xi}=D\hbar\pr({\bolds\xi
})\Gamb_\varthetab D\hbar({\bolds\xi})$---throughout, we write
$D\hbar\pr(\cdot)$, $D\bbar\pr(\cdot)$, etc., instead of $(D\hbar(\cdot))\pr
$, $(D\bbar(\cdot))\pr$, etc. In general, $\Thetab$ is a nonlinear
manifold of $\rr^p$; the experiment parametrized by $\Thetab$ then is
a \textit{curved} experiment.

Next, denoting by $C$ an $r$-dimensional manifold in $\mathbb{R}^p$,
$r<p$, consider the null hypothesis
$\mathcal{H}_0\dvtx\varthetab\in C \cap\Thetab$---in general, a~nonlinear restriction of the parameter space $\Thetab$.
The same hypothesis can be expressed in the ${\bolds\xi
}$-parametrization as $\mathcal{H}_0\dvtx{\bolds\xi} \in\Xib_0$,
where $\Xib_0:=\hbar^{-1}(C \cap\Thetab)$ is a ($\ell
$-dimensional, say) submanifold of~$\Xib$. Fix ${\bolds\xi}_0 =\hbar
^{-1}(\varthetab_0)\in\Xib_0$, and let $\lbar\dvtx B\subset
\mathbb{R}^\ell\to\Xib$ be a local (at ${\bolds\xi}_0$) chart for
this manifold.

Define ${\bolds\alpha}_0:=\lbar^{-1}({\bolds\xi}_0)$. At
${\bolds\xi}_0$, $\mathcal{H}_0$ is linearized into
$\mathcal{H}_{{\bolds\xi}_0}\dvtx{\bolds\xi} \in{\bolds\xi}_0 +
\mathcal{M}(D\lbar({\bolds\alpha}_0))$, where $D\lbar({\bolds
\alpha}_0)$ is the Jacobian matrix of $\lbar$ (with rank $\ell$)
computed at ${\bolds\alpha}_0$ and $\mathcal{M}(\mathbf{A})$ denotes
the vector space spanned by the columns of a matrix $\mathbf A$. At
${\bolds\alpha}_0$, a~locally asymptotically most stringent test
statistic (at ${\bolds\xi}_0$) for $\mathcal{H}_{{\bolds\xi}_0}$ is
%
%e3.5 ###
%
\begin{equation}\label{lams}
Q_{{\bolds\xi}_0}:=
\Deltab_{{\bolds\xi}_0}\pr
\bigl( \Gamb_{{\bolds\xi}_0}^{-1} - D\lbar({\bolds\alpha}_0) ( D\lbar
\pr({\bolds\alpha}_0)\Gamb_{{\bolds\xi}_0}D\lbar({\bolds\alpha
}_0))^{-1} D\lbar\pr({\bolds\alpha}_0) \bigr)
\Deltab_{{\bolds\xi}_0}
\end{equation}
[see Section 11.9 of \citet{Lcam86}]. This test statistic is nothing else
but the squared Euclidean norm of the orthogonal projection, onto the
linear space orthogonal to $ \Gamb_{{\bolds\xi}_0}^{1/2}D\lbar
({\bolds\alpha}_0)$, of the standardized central sequence $ \Gamb
_{{\bolds\xi}_0}^{-1/2}\Deltab_{{\bolds\xi}_0}$.
In view of ULAN, the asymptotic behavior of $\Deltab_{{\bolds\xi
}_0}$ is the same under local alternatives in $ \Xib_0$ as under local
alternatives in $ {\bolds\xi}_0 + \mathcal{M}(D\lbar({\bolds\alpha
}_0))$, so that the same test statistic~$
Q_{{\bolds\xi}_0}$, which (at ${\bolds\xi}_0$) is locally
asymptotically most stringent for $\mathcal{H}_{{\bolds\xi}_0}$, is
also locally asymptotically most stringent for $\mathcal{H}_{0}$.

In many cases, however, it is highly desirable to express the most
stringent statistic in the curved $\Thetab$-parametrization, which, as
is the case for the eigenvalues/eigenvectors problems considered in
this work, is the natural parametrization. This is the objective of the
following result (see the \hyperref[app]{Appendix} for the proof).
\begin{Prop}\label{Optitest}
With the same notation as above, a~locally asymptotically most
stringent statistic (at $\varthetab_0$) for testing
$\mathcal{H}_{0}\dvtx\varthetab\in C\cap\Thetab$ is
%
%e3.6 ###
%
\begin{equation} \label{Qxi}
Q_{{\bolds\xi}_{0}}
=
Q_{\varthetab_{0}}:=\Deltab_{\varthetab_0}\pr
\bigl( \Gamb_{\varthetab_0}^{-} - D\btilde(\etab_0) ( D\btilde\pr
(\etab_0)\Gamb_{\varthetab_0}D\btilde(\etab_0))^{-} D\btilde\pr
(\etab_0) \bigr)
\Deltab_{\varthetab_0},\hspace*{-34pt}
\end{equation}
where $\btilde\dvtx A\subset\mathbb{R}^\ell\to\mathbb{R}^p$ is a
local (at $\varthetab_0$) chart for the tangent (still at $\varthetab
_0$) to the manifold $C\cap\Thetab$, $\etab_{0}:= \bbar
^{-1}(\varthetab_{0})$, and $\Ab^-$ denotes the Moore--Penrose inverse
of~$\Ab$.
\end{Prop}

Hence, a~locally asymptotically most stringent (at ${\bolds\xi}_0$ or
$\varthetab_0$, depending on the para\-metrization) test for $\mathcal
{H}_0$ can be based on either of the two quadratic forms $Q_{{\bolds
\xi}_0}$ or $Q_{\varthetab_0}$, which coincide, and are
asymptotically chi-square $[(m-\ell)$ degrees of freedom] under $\mathrm
{P}_{{\bolds\xi}_0}^{(n)}= \mathrm{P}_{\varthetab_0}^{(n)}$, for
$\varthetab_0 = \hbar({\bolds\xi}_0)$. For practical
implementation, of course, an adequately discretized root-$n$
consistent estimator has to be substituted for the unknown $\varthetab
_0 $ or ${\bolds\xi}_0$---which asymptotically does not affect the
test statistic.

Provided that $\Xib$ remains an open subset of $\rr^m$, the
assumption of a full-rank information matrix $\Gamb_{\bolds\xi}$ is
not required. Hallin and Puri [(\citeyear{HP94}), Lemma~5.12] indeed have shown, in
the case of ARMA experiments, that~(\ref{lams}) remains locally
asymptotically most stringent provided that generalized inverses (not
necessarily Moore--Penrose ones) are substituted for the inverses of
noninvertible matrices, yielding
\[
Q_{{\bolds\xi}_0}:=
\Deltab_{{\bolds\xi}_0}\pr
\bigl( \Gamb_{{\bolds\xi}_0}^{-} - D\lbar({\bolds\alpha}_0) ( D\lbar
\pr({\bolds\alpha}_0)\Gamb_{{\bolds\xi}_0}D\lbar({\bolds\alpha
}_0))^{-} D\lbar\pr({\bolds\alpha}_0) \bigr)
\Deltab_{{\bolds\xi}_0}.
\]
The same reasoning as in the proof of Proposition~\ref{Optitest} then
applies, mutatis mutandis, when ``translating'' $Q_{{\bolds\xi}_0}$
into $Q_{\varthetab_0}$ (with appropriate degrees of freedom).

%s4 ###
\section{Parametrically optimal tests for principal components}
\label{paramtests}

%s4.1 ###
\subsection{Optimal parametric tests for eigenvectors}\label{howeigenvectors}
Testing the hypothesis $\mathcal{H}^{\betab}_0$ on eigenvectors is a
particular case of the problem considered in the previous section. The
$\operatorname{vech}^+(\iotab)$ parametrization [$\iotab$ an arbitrary
skew-symmetric $(k\times k)$ matrix] yields a standard ULAN experiment,
with parameter
\[
{\bolds\xi}:= ({\bolds\theta}\pr, \sigma^2, (\dvecrond(\Lamb
_{\Vb}))\pr, (\operatorname{vech}^+(\iotab))\pr
)\pr\in\rr^k\times\rr^+\times\mathcal{C}^{k-1}\times\rr
^{k(k-1)/2}=:\Xib,
\]
hence $m= k(k+3)/2$, while Proposition~\ref{LAN} provides the curved
ULAN experiment, with parameter $\varthetab\in\Thetab\subset\rr
^{p}$ and $p=k(k+2)$. ULAN for the ${\bolds\xi}$-experiment readily
follows from the fact that the mapping $ \operatorname{vech}^+(\iotab
)\mapsto\operatorname{vec}(\betab)= \operatorname{vec}(\exp(\iotab
))$ is continuously differentiable.

As explained before, the block-diagonal structure of the information
matrix~(\ref{Gamb}) implies that locally asymptotically optimal
inference about $\betab$ can be based on $\Deltab^{\IV}_{\varthetab
;f_1}$ only, as if ${\bolds\theta}$, $\sigma^2$ and $\dvecrond
({\bolds\Lambda_\Vb})$ were specified. Since this also allows for
simpler exposition and lighter notation, let us assume that these
parameters take on specified values ${\bolds\theta}$, $\sigma^2$
and $(\lambda_{2; \Vb},\ldots, \lambda_{k; \Vb})$, respectively.
The resulting experiment then is parametrized either by $\vecop
\betab\in\vecop({\mathcal SO}_k)\subset\rr^{k^2}$ (playing the
role of $\varthetab\in\Thetab\subset\rr^{p}$ in the notation of
Proposition~\ref{Optitest}) or by $\operatorname{vech}^+(\iotab)\in
\rr^{k(k-1)/2}$ (playing the role of ${\bolds\xi}$).

In this experiment, the null hypothesis $\mathcal{H}^{\betab}_0$
consists in the intersection of the linear manifold
$C:= ( \betab^{0\prime},\mathbf{0}_{1\times(k-1)k}%,\ldots, \mathbf{0}_k
)\pr+ \mathcal{M}(\Upsib% ^\betab_0
)$, where
$\Upsib%^\betab_0
:=
(
\mathbf{0}_{k(k-1)\times k},\break
\mathbf{I}_{k(k-1)}
)\pr,
$
with the nonlinear manifold $\operatorname{vec}({\mathcal SO}_k)$. Let $\betab
_{0}:=(\betab^{0}, \betab_{2} ,\ldots, \betab_{k})$ be such that
$\operatorname{vec}(\betab_0)$ belongs to that intersection.
In view of Proposition~\ref{Optitest}, a~most stringent test statistic
[at $\operatorname{vec}(\betab_{0})$] for $\mathcal{H}^{\betab}_0$
requires a chart for the tangent to $C\cap\vecop({\mathcal SO}_k)$ at
$\operatorname{vec}(\betab_0)$. It follows from~(\ref{tangentbetab}) that this
tangent space reduces to
\[
\{ \vecop( \betab_{0} + \mathbf{b}) | \mathbf{b}:=(\mathbf{0} , \mathbf
{b}_{2}, \ldots, \mathbf{b}_{k}) \mbox{ such that } \betab_{0}\pr
\mathbf{b} + {\mathbf{b}}\pr\betab_{0}= \mathbf{0} \}.
\]

Solving for $\operatorname{vec}(\mathbf{b})=(\mathbf{0}\pr, \mathbf{b}_{2}\pr, \ldots,
\mathbf{b}_{k}\pr)\pr$ the system of constraints $\betab_{0}\pr\mathbf
{b} +
{\mathbf{b}}\pr\betab_{0}= \mathbf{0}$ yields $\operatorname{vec}(\mathbf{b}) \in
\mathcal
{M}(\mathbf{P}_{k}^{\betab_{0}})$, where
%
%e4.1 ###
%
\begin{equation}\label{defP}\quad
\mathbf{P}_{k}^{\betab_{0}}:=
\pmatrix{
\mathbf{0}_{k \times k(k-1)}\vspace*{2pt}\cr
\mathbf{I}_{k-1}\otimes[ \mathbf{I}_k- \betab^{0}\betab^{0\prime
} ]
-
\displaystyle\sum_{i,j=1}^{k-1}
[\mathbf{e}_{i;k-1}\mathbf{e}_{j;k-1}\pr
\otimes
\betab_{j+1}\betab_{i+1}\pr]}
\end{equation}
(with $\mathbf{e}_{i;k-1}$ denoting the $i$th vector of the canonical
basis of $\R^{k-1}$). A~local chart for the tangent space of interest is then simply $\btilde
\dvtx{\bolds\eta}\in\rr^{(k-1)k}\mapsto\btilde({\bolds\eta}):=
\vecop(\betab_0) + \mathbf{P}_{k}^{\betab_{0}}{\bolds\eta}$, with
$\etab_0=\btilde^{-1}(\vecop(\betab_{0}))=\mathbf{0}_{(k-1)k}$ and
$D\btilde(\etab_0)=\mathbf{P}_{k}^{\betab_{0}}$.
Letting $\varthetab_0:=({\bolds\theta}' , \sigma^2, (\dvecrond
\Lamb_{\Vb})' ,(\vecop\betab_0)')'$, the test statistic (\ref
{Qxi}) takes the form
%
%e4.2 ###
%
\begin{eqnarray}\label{Qparam}\quad
Q_{\varthetab_0 ; f_{1}}^{(n)}
&=&
\Deltab_{\varthetab_0;f_1}^{\IV\prime}
[(\Gamb^{\IV}_{\varthetab_0 ;f_1})^{-}- \mathbf{P}_{k}^{\betab
_{0}} ((\mathbf{P}_{k}^{\betab_{0}}) \pr\Gamb^{\IV
}_{\varthetab_0 ;f_1} \mathbf{P}_{k}^{\betab_{0}} )^{-} (\mathbf
{P}_{k}^{\betab_{0}}) \pr]
\Deltab_{\varthetab_0 ;f_1}^{\IV} \nonumber\\[-8pt]\\[-8pt]
&=&\frac{nk(k+2)}{ \mathcal{J}_k(f_1)} \sum_{j=2}^{k} \bigl(
\betab_j\pr\mathbf{S}_{\varthetab_0;f_1}^{(n)}\betab^0\bigr)^2,\nonumber
\end{eqnarray}
with
%
%e4.3 ###
%
\begin{equation}\label{Sparam}
\mathbf{S}_{\varthetab;f_1}^{(n)}:=
\frac{1}{n}
\sum_{i=1}^n \varphi_{f_1}
\biggl(\frac{d_{i}({\bolds\theta},\Vb)}{\sigma} \biggr) \frac
{d_{i}({\bolds\theta},\Vb)}{\sigma} \Ub_{i}({\bolds\theta},\Vb
) \Ub_{i}\pr({\bolds\theta},\Vb) ,
\end{equation}
where $\Vb$ denotes the unique shape value associated with the
parameter $\varthetab$.

After simple algebra, we obtain
%
%e4.4 ###
%
\begin{eqnarray} \label{degreefreedom}
&&\Gamb_{\varthetab_{0};f_1}^{\IV} [(\Gamb^{\IV}_{\varthetab_0
;f_1})^{-}- \mathbf{P}_{k}^{\betab_{0}} ((\mathbf{P}_{k}^{\betab
_{0}}) \pr\Gamb^{\IV}_{\varthetab_0 ;f_1} \mathbf{P}_{k}^{\betab
_{0}} )^{-} (\mathbf{P}_{k}^{\betab_{0}}) \pr]
\nonumber\\[-8pt]\\[-8pt]
&&\qquad
= \tfrac{1}{2}\mathbf{G}_{k}^{\betab_0}
\operatorname{diag}\bigl(\mathbf{I}_{k-1},\mathbf{0}_{(k-2)(k-1)/2 \times
(k-2)(k-1)/2}\bigr)(\mathbf{G}_{k}^{\betab_0})\pr, \nonumber
\end{eqnarray}
which is idempotent with rank $(k-1)$. Since, moreover, $\Deltab
_{\varthetab_{0};f_1 }^{\IV}$, under $\mathrm{P}^{(n)}_{\varthetab
_{0};f_{1}}$, is asymptotically $\mathcal{N}( \mathbf{0}, \Gamb
_{\varthetab_{0};f_1}^{\IV})$, Theorem 9.2.1 in \citet{RM71}
then shows that $ Q_{\varthetab_0 ; f_{1}}^{(n)}$, still under $\mathrm
{P}^{(n)}_{\varthetab_{0};f_{1}}$, is asymptotically chi-square with
$(k-1)$ degrees of freedom.

The resulting test, which rejects $\mathcal{H}_0^{\betab}$ at
asymptotic level $\alpha$ whenever $ Q_{\varthetab_0 ; f_{1}}^{(n)}$
exceeds the $\alpha$-upper quantile $\chi^2_{k-1,1-\alpha}$ of the
$\chi^2_{k-1}$ distribution, will be denoted as $\phi^{(n)}_{\betab;
f_1}$.\vspace*{1pt} It is locally asymptotically most stringent, at $\varthetab_0$
and under correctly specified standardized radial density $f_1$ (an
unrealistic assumption). Of course, even if $f_{1}$ were supposed to be
known, $Q_{\varthetab_0 ; f_{1}}^{(n)}$ still depends on the
unspecified ${\bolds\theta}, \sigma^2, \Lamb_{\Vb}$ and $\betab
_2,\ldots,\betab_k$. In order to obtain a genuine test statistic,
providing a locally asymptotically most stringent test at \textit{any}
$\varthetab_0\in\mathcal{H}_0^\betab$ (with an obvious abuse
of notation), we would need replacing those nuisance parameters with
adequate estimates. We will not pursue any further with this problem
here, as it is of little practical interest for arbitrary density
$f_1$. The same problem will be considered in Section~\ref{gausscase}
for the Gaussian and pseudo-Gaussian versions of~(\ref{Qparam}), then
in Section~\ref{ranktests} for the rank-based ones.

%s4.2 ###
\subsection{Optimal parametric tests for eigenvalues}\label{eigenvaluesprob}

We now turn to the problem of testing the null hypothesis ${\mathcal
H}_{0}^{\Lamb}\dvtx\sum_{j=q+1}^{k} \lambda_{j; \Vb}- p \sum
_{j=1}^{k}\lambda_{j; \Vb}=0$ against alternatives of the form
${\mathcal H}_{1}^{\Lamb}\dvtx\sum_{j=q+1}^{k} \lambda_{j; \Vb}- p
\sum_{j=1}^{k}\lambda_{j; \Vb}<0$, for given $p \in(0,1)$.
Letting
\[
h\dvtx(\lambda_{2 }, \lambda_{3 }, \ldots,\lambda_{k })\pr\in
\mathcal{C}^{k-1}
%(\subset(\rr^+_0)^{k-1})
\mapsto
\sum_{j=q+1}^{k} \lambda_{j }- p
\Biggl(
\prod_{j=2}^{k}\lambda_{j }^{-1}+ \sum_{j=2}^{k}\lambda_{j }
\Biggr)
\]
and recalling that $ \prod_{j=1}^{k} \lambda_{j;\Vb}=1$, ${\mathcal
H}_{0}^{\Lamb}$ rewrites, in terms of $\dvecrond(\Lamb_{\Vb})$, as
${\mathcal H}_{0}^{\Lamb}\dvtx h( \dvecrond(\Lamb_{\Vb})) =0$, a~highly
nonlinear but smooth constraint on $\dvecrond(\Lamb_{\Vb})$.
It is easy to check that, when computed at $\dvecrond(\Lamb_{\Vb})$,
the gradient of $h$ is %given by
\begin{eqnarray*}
&&
\operatorname{grad} h(\dvecrond(\Lamb_{\Vb}))
\\
&&\qquad= \bigl(
p ( \lambda_{1;\Vb}\lambda_{2;\Vb}^{-1}-1),
\ldots,
p ( \lambda_{1;\Vb}\lambda_{q;\Vb}^{-1}-1),
\\
&&\hspace*{36.74pt}
1+p ( \lambda_{1;\Vb} \lambda_{q+1; \Vb}^{-1}-1 ), \ldots, 1+p (
\lambda_{1;\Vb} \lambda_{k; \Vb}^{-1}-1 )
\bigr)\pr.
\end{eqnarray*}

Here again, in view of the block-diagonal form of the information
matrix, we may restrict our attention to the $\dvecrond(\Lamb_{\Vb
})$-part $\Deltab_{\varthetab;f_1}^{\III}$ of the central sequence
as if ${\bolds\theta}$, $\sigma^2$ and $\betab$ were known; the
parameter space then reduces to the $(k-1)$-dimensional open
cone $\mathcal{C}^{k-1}$. Testing a nonlinear constraint on a
parameter ranging over an open subset of $\rr^{k-1}$ is much easier
however than the corresponding problem involving a curved experiment,
irrespective of the possible noninvertibility of the information\vadjust{\goodbreak}
matrix. In the noncurved experiment, indeed, a~linearized version
${\mathcal H}_{0,\mathrm{lin}}^{\Lamb} \dvtx\dvecrond(\Lamb_{\Vb})\in
\dvecrond(\Lamb_0) +
{\mathcal M}\ort(\operatorname{grad} h(\dvecrond\Lamb_0))
$ of ${\mathcal H}_{0}^{\Lamb}$ in the vicinity of $\dvecrond({\bolds
\Lambda}_0)$ satisfying $h( \dvecrond\Lamb_0 )=0$ makes sense
[${\mathcal M}\ort(\mathbf{A})$ denotes the orthogonal complement
of ${\mathcal M}(\mathbf{A})$]. And, as mentioned in Section \ref
{curved}, under ULAN, the asymptotic behavior of $\Deltab^{\III}
_{{\bolds\vartheta}_0;f_1}$, with $\varthetab_0=({\bolds\theta}\pr
,\sigma^2,(\dvecrond\Lamb_0)\pr,(\vecop\betab)\pr)\pr$, is
locally the same under ${\mathcal H}_{0,\mathrm{lin}}^{\Lamb}$
as under ${\mathcal H}_{0}^{\Lamb}$. As for the ``linearized
alternative'' ${\mathcal H}_{1,\mathrm{lin}}^{\Lamb} $ consisting of
all $\dvecrond\Lamb$ values such that $(\dvecrond\Lamb- \dvecrond
\Lamb_0)\pr\operatorname{grad} h(\dvecrond\Lamb_0)<0$, it locally
and asymptotically coincides with ${\mathcal H}_{1}^{\Lamb}$: indeed,
although the symmetric difference ${\mathcal H}_{1}^{\Lamb} \Delta
{\mathcal H}_{1,\mathrm{lin}}^{\Lamb} $, for fixed $n$, is not empty,
any $\dvecrond\Lamb_0 + n^{-1/2}\taub^{\III}\in{\mathcal
H}_{1,\mathrm{lin}}^{\Lamb} $ eventually belongs to ${\mathcal
H}_{1}^{\Lamb} $, and conversely. Therefore, a~locally (at
$\dvecrond\Lamb_0$) asymptotically optimal test for ${\mathcal
H}_{0,\mathrm{lin}}^{\Lamb}$ against ${\mathcal H}_{1,\mathrm
{lin}}^{\Lamb} $ is also locally asymptotically optimal for ${\mathcal
H}_{0 }^{\Lamb}$ against ${\mathcal H}_{1 }^{\Lamb} $, and
conversely, whatever the local asymptotic optimality concept adopted.
Now, in the problem of testing ${\mathcal H}_{0,\mathrm{lin}}^{\Lamb
}$ against ${\mathcal H}_{1,\mathrm{lin}}^{\Lamb} $ the null
hypothesis is (locally) a hyperplane of $\rr^{k-1}$, with an
alternative consisting of the halfspace lying ``below'' that hyperplane.
For such one-sided problems (locally and asymptotically, still
at $\varthetab_0$) uniformly most powerful tests exist; a \textit{most
powerful} test statistic is [\citet{Lcam86}, Section 11.9]
%
%e4.5 ###
%
\begin{eqnarray}\label{Tparam1}
T_{\varthetab_0; f_{1}}^{(n)}&:=& (\operatorname{grad}\pr
h(\dvecrond\Lamb_{0})
(\Gamb^{\III}_{{\bolds\vartheta}_0;f_1})^{-1}
\operatorname{grad} h(\dvecrond\Lamb_{0}) )^{-1/2}
\nonumber\\[-8pt]\\[-8pt]
& &{}
\times\operatorname{grad}\pr h(\dvecrond\Lamb_{0}) (\Gamb
^{\III}_{{\bolds\vartheta}_0;f_1})^{-1} \Deltab^{\III} _{{\bolds
\vartheta}_0;f_1},\nonumber
\end{eqnarray}
which, under $\mathrm{P}^{(n)}_{\varthetab_0 ; f_1}$, is asymptotically
standard normal.
An explicit form of $T_{\varthetab_0; f_{1}}^{(n)}$ requires a closed
form expression of the inverse of $\Gamb^{\III}_{{\bolds\vartheta
};f_1} =
({\mathcal{J}_k(f_1)}/\break{k(k+2)})\times\mathbf{D}_k(\Lamb_\Vb) $. The following
lemma provides such an expression for the inverse of $\mathbf
{D}_{k}(\Lamb
_{\Vb})$ (see the \hyperref[app]{Appendix} for the proof).
\begin{Lem} \label{infoinverse}
Let $\mathbf{P}_{k}^{\Lamb_\Vb}:=\mathbf{I}_{k^{2}}- \frac{1}{k} \Lamb
_{\Vb}^{\otimes2} \vecop(\Lamb_{\Vb}^{-1})(\vecop(\Lamb_{\Vb
}^{-1}))\pr$ and $\mathbf{N}_k:=(\mathbf{0}_{(k-1)\times1},\mathbf
{I}_{k-1})$. Then,
$(\mathbf{D}_{k}(\Lamb_{\Vb}))^{-1}=\mathbf{N}_{k}\mathbf{H}_{k} \mathbf
{P}_{k}^{\Lamb_\Vb} (\mathbf{I}_{k^{2}}+ \mathbf{K}_{k}) \Lamb_{\Vb
}^{\otimes2} (\mathbf{P}_{k}^{\Lamb_\Vb} )\pr\mathbf{H}_{k}\pr
\mathbf{N}_{k}\pr$.
\end{Lem}

Using this lemma, it follows after some algebra that, for any
$\varthetab_0\in{\mathcal H}_{0}^{\Lamb}$,
\begin{eqnarray*}
&&
\operatorname{grad}\pr h (\dvecrond\Lamb_{0})(\mathbf{D}_{k}(\Lamb
_{0}))^{-1}\operatorname{grad} h(\dvecrond\Lamb_{0})\\
&&\qquad
=
2 \Biggl\{ p^{2} \sum_{j=1}^{q} \lambda_{j;0}^{2} + (1-p)^{2} \sum
_{j=q+1}^{k} \lambda_{j;0}^{2} \Biggr\}=a_{p,q}(\Lamb_0)
\end{eqnarray*}
[where ${\bolds\Lambda}\mapsto a_{p,q}({\bolds\Lambda})$ is the
mapping defined in~(\ref{TAnd})],
and
\[
\operatorname{grad}\pr h (\dvecrond\Lamb_{0}) (\mathbf{D}_{k}(\Lamb
_{0}))^{-1} \mathbf{M}_{k}^{\Lamb_{0}}\mathbf{H}_{k}({\Lamb
_{0}^{-1/2}})^{\otimes2}=\mathbf{c}_{p,q}\pr\mathbf{H}_{k}({\Lamb
_{0}^{1/2}})^{\otimes2}.
\]
This and the definition of $\mathbf{H}_{k}$ yields
%
%e4.6 ###
%
\begin{equation}\label{Tparam}\qquad
T_{\varthetab_0; f_{1}}^{(n)}= \biggl( \frac{nk(k+2)}{\mathcal
{J}_k(f_1)} \biggr)^{1/2}
(a_{p,q}(\Lamb_{0} ))^{-1/2}
\mathbf{c}_{p,q}\pr\dvec\bigl(\Lamb_{0}^{1/2} \betab\pr\mathbf
{S}_{\varthetab_0;f_{1}} ^{(n)}\betab\Lamb_{0}^{1/2}\bigr)\hspace*{-20pt}
\end{equation}
with $\mathbf{S}_{\varthetab; f_{1}}^{(n)}$ defined in~(\ref{Sparam}).
The corresponding test, which rejects ${\mathcal H}_{0}^{\Lamb}$ for
small values of $T_{\varthetab_0; f_{1}}^{(n)}$, will be denoted as
$\phi^{(n)}_{\Lamb; f_1}$.

%s4.3 ###
\subsection{Estimation of nuisance parameters}\label{estimnuis}

The tests $\phi^{(n)}_{\betab; f_1}$ and $\phi^{(n)}_{\Lamb; f_1}$
derived in Sections~\ref{howeigenvectors} and~\ref{eigenvaluesprob}
typically are valid under standardized radial density $f_1$ only; they
mainly settle the optimality bounds at given density~$f_1$, and are of
little practical value. Due to its central role in multivariate
analysis, the Gaussian case ($f_1=\phi_1$) is an exception. In this
subsection devoted to the treatment of nuisance parameters, we
therefore concentrate on the Gaussian tests $\phi^{(n)}_{\betab; \phi
_1}$ and $\phi^{(n)}_{\Lamb; \phi_1}$, to be considered in more
detail in Section~\ref{gausscase}.

The test statistics derived in Sections~\ref{howeigenvectors} and \ref
{eigenvaluesprob} indeed still involve nuisance parameters which in
practice have to be replaced with estimators. The traditional way of
handling this substitution in ULAN families consists in assuming, for a
null hypothesis of the form ${\varthetab} \in{\mathcal H}_{0}$, the
existence of a sequence $\hat{\varthetab}^{(n)}$ of estimators of
${\varthetab}$ satisfying all or part of the following assumptions (in
the notation of this paper).
\renewcommand{\theassumption}{$(\mathrm{B})$}
\begin{assumption}\label{assuB}
We say that a sequence of estimators
($\hat{\varthetab}^{(n)}, n\in\N$) satisfies Assumption~\ref{assuB}
for the null $\mathcal{H}_0$ and the density $f_1$ if
$\hat{\varthetab}^{(n)}$ is:

\begin{enumerate}[(B3)]
\item[(B1)]\hypertarget{B1} \mbox{\textit{constrained}}: $\mathrm{P}^{(n)} _{\varthetab;
f_{1}} [\hat{\varthetab}^{(n)} \in{\mathcal H}_{0} ] =1$
for all $n$ and all $\varthetab\in
{\mathcal H}_{0}$;

\item[(B2)]\hypertarget{B2} \mbox{\textit{$n^{1/2}$-consistent}}: for all $\varthetab\in
{\mathcal H}_{0} $, $n^{1/2}(\hat{\varthetab}^{(n)}-\varthetab
)=O_{\mathrm{P}}(1)$ under $\mathrm{P}^{(n)}_{\varthetab; f_1}$,
as $n\rightarrow\infty$;

\item[(B3)]\hypertarget{B3} \mbox{\textit{locally asymptotically discrete}}: for all $\varthetab\in
{\mathcal H}_{0}$ and all $c>0$, there exists $M=M(c)>0$ such that the number
of possible values of $\hat{\varthetab}^{(n)}$ in balls of the form
$\{
\mathbf{t}%\in\R^{L}
\dvtx n^{1/2} \| (\mathbf{t}-{\varthetab}) \| \leq c\}$ is bounded
by $M$, uniformly as $n\rightarrow\infty$.
\end{enumerate}
\end{assumption}

These assumptions will be used later on. In the Gaussian or
pseudo-Gaussian context we are considering here, however, Assumption
\hyperlink{B3}{(B3)} can be dispensed with under arbitrary densities with finite
fourth-order moments. The following asymptotic linearity result
characterizes the asymptotic impact,\vspace*{1pt} on $\Deltab_{\varthetab; \phi
_1}^{\III}$ and $\Deltab_{\varthetab; \phi_1}^{\IV}$, under any
elliptical density $g_1$ with finite fourth-order moments, of
estimating $\varthetab$ (see the \hyperref[app]{Appendix} for the proof).
\begin{Lem}\label{parametricasymplin}
Let Assumption~\ref{assuA} hold, fix $\varthetab\in\Thetab$ and $g_{1} \in
{\mathcal F}_{1}^{4}$, and write $D_k(g_1):=\mu_{k+1;g_1}/\mu
_{k-1;g_1}$. Then, for any
root-$n$ consistent estimator
$\hat{\varthetab}:= ( \hat{\varthetab}^{\I\prime} , \hat
{\vartheta}^\II, \hat{\varthetab}^{\III\prime}, \hat{\varthetab
}^{\IV\prime} )\pr$
of $\varthetab$ under $\mathrm{P}^{(n)}_{\varthetab; g_{1}}$, both\vadjust{\goodbreak}
$\Deltab_{\hat{\varthetab}; \phi_1}^{\III}- \Deltab_{\varthetab;
\phi_1}^{\III}+ a_{k} ( D_{k}(g_{1})/k) \times\Gamb
_{\varthetab; \phi_1}^{\III}
n^{1/2}(\hat{\varthetab}^{\III}-{\varthetab}^{\III})$ and\vspace*{-2pt}
$\Deltab_{\hat{\varthetab}; \phi_1}^{\IV}- \Deltab_{\varthetab;
\phi_1}^{\IV}+ a_{k} ( D_{k}(g_{1})/k) \Gamb
_{\varthetab; \phi_1}^{\IV} n^{1/2}(\hat{\varthetab}^{\IV
}-{\varthetab}^{\IV}) $
are $o_\mathrm{P}(1)$ under $\mathrm{P}_{\varthetab; g_{1}}^{(n)}$, as
$\ny
$, where $a_k$ was defined in Section~\ref{defelliptttt}.
\end{Lem}

%s5 ###
\section{\mbox{Optimal Gaussian and pseudo-Gaussian tests for principal
components}}\label{gausscase}

%s5.1 ###
\subsection{Optimal Gaussian tests for eigenvectors}\label{Gaussbetab}
For $f_1=\phi_1$, the test statistic~(\ref{Qparam}) takes the form
%
%e5.1 ###
%
\begin{equation} \label{betagauss}\quad
Q_{\varthetab_0 ; \phi_1 }^{(n)}= n \sum_{j=2}^{k} \bigl({\betab
}_{j}\pr\mathbf{S}^{(n)}_{\varthetab_0;\phi_1 } \betab_{}^{0}
\bigr)^{2} =
n \betab^{0\prime} \mathbf{S}^{(n)}_{\varthetab_0;\phi_1 } (\mathbf
{I}_k-\betab_{}^{0}\betab^{0\prime}) \mathbf{S}^{(n)}_{\varthetab
_0;\phi_1 } \betab_{}^{0},
\end{equation}
with ${\Sb}_{\varthetab;\phi_1 }^{(n)}:= \frac{a_{k}}{n\sigma^2}
\sum_{i=1}^{n} \Vb^{-1/2}(\Xb_{i}- {\bolds\theta})(\Xb_{i}-
{\bolds\theta})\pr\Vb^{-1/2}$. This statistic still depends on
nuisance parameters, to be replaced with estimators. Letting $\Sb^{(n)}
=
\frac{1}{n} \sum_{i=1}^{n} ({\Xb}_{i}-\bar{\mathbf X})({\Xb}_{i}-\bar
{\mathbf{X}})\pr$, a~natural choice for such estimators would be
$\hat{\bolds\theta}=\bar{\mathbf X}:=\frac{1}{n}\sum_{i=1}^n\Xb_i$
and
\[
\Sb^{(n)}
=:
\bigl| \Sb^{(n)}\bigr|^{1/k} \hats{\Vb}
=
\biggl( \frac{| \Sb^{(n)}|^{1/k}}{\hat\sigma^2} \biggr) \hat
\sigma^2 \hats{\Vb}
=:
\biggl( \frac{| \Sb^{(n)}|^{1/k}}{\hat\sigma^2} \biggr) \hat
\sigma^2
\hat\betab_\mathbf{V} {\hat\Lamb}_{\Vb}\hat\betab_\mathbf{V}\pr,
\]
where ${\hat\Lamb}_{\Vb}$ is the diagonal matrix collecting the
eigenvalues of $\hats{\Vb}$ (ranked in decreasing order), $\hat\betab
_\mathbf{V}
:=
( \hat\betab_{1;\mathbf{V}}, \ldots, \hat\betab_{k;\mathbf{V}})$ is the
corresponding matrix of eigenvectors, and $\hat\sigma^2$ stands for
the empirical median of $d_{i}^{2}(\bar{\Xb}, \hats{\Vb})$, $i=
1,\ldots,n$.
%[Link with the Gaussian MLEs].
For $\betab$, however, we need a constrained estimator ${\tilde\betab
}$ satisfying Assumption~\ref{assuB} for $\mathcal{H}^\betab_0$ ($\hat\betab_\mathbf{V}$ in
general does not). Thus, we rather propose estimating $\varthetab$ by
%
%e5.2 ###
%
\begin{equation} \label{choiceprelim}
\hat{\varthetab}:=
(\bar{\mathbf X}\pr, \hat{\sigma}^{2},
(\dvecrond\hat\Lamb_{\Vb})\pr,
(\vecop{\tilde\betab}_{0})\pr)\pr,
\end{equation}
where
${\tilde\betab}_{0}:=( \betab_{}^{0} , \tilde{\betab}_{2}, \ldots
, \tilde{\betab}_{k})$ can be obtained from $(\hat\betab_{1;\Vb},
\ldots,\hat\betab_{k;\Vb})$ via the following Gram--Schmidt
technique. Let
${\tilde\betab}_{2}:=(\mathbf{I}_k- \betab_{}^{0}\betab_{}^{0 \prime
}) \hat{\betab}_{2;\Vb}/ \|(\mathbf{I}_k-\break \betab_{}^{0}\betab_{}^{0
\prime}) \hat{\betab}_{2;\Vb} \|$. By construction, ${\tilde\betab
}_{2}$ is the unit-length vector proportional to the projection of the
second eigenvector of $\mathbf{S}^{(n)}$ onto the space which is
orthogonal to $\betab_{}^{0}$. Iterating this procedure, define
\[
{\tilde\betab}_{j}=\frac{(\mathbf{I}_k- \betab_{}^{0}\betab_{}^{0
\prime}-\sum_{h= 2 }^{j-1} \tilde{\betab}_{h}\tilde{\betab
}_{h}\pr) \hat{\betab}_{j;\Vb}}{\| (\mathbf{I}_k- \betab_{}^{0}\betab
_{}^{0 \prime}-\sum_{h= 2}^{j-1} \tilde{\betab}_{h}\tilde{\betab
}_{h}\pr) \hat{\betab}_{j;\Vb}\|},\qquad j=3, \ldots, k.
\]
The corresponding (constrained) estimator of the scatter $\Sigb$ is
$\tilde\Sigb
:=
\hat\sigma^2 \tilde\Vb
:=
\hat\sigma^2 {\tilde\betab}_{0} {\hat\Lamb}_{\Vb} {\tilde
\betab}{}'_{0}$.

It is easy to see that ${\tilde\betab}_{0}$, under $\varthetab_0 \in
{\mathcal H}_{0}^{\betab}$, inherits $\hat\betab_\mathbf{V}$'s root-$n$
consistency, which holds under any elliptical density $g_1$ with finite
fourth-order moments.
Lemma~\ref{parametricasymplin} thus applies. Combining Lemma \ref
{parametricasymplin} with~(\ref{degreefreedom}) and the fact that
\[
\mathbf{G}_{k}^{\betab_0}
\operatorname{diag}\bigl(\mathbf{I}_{k-1},\mathbf{0}_{(k-2)(k-1)/2 \times
(k-2)(k-1)/2}\bigr)(\mathbf{G}_{k}^{\betab_0})\pr\vecop(\tilde{\betab
}_{0}- \betab_0)=\mathbf{0}
\]
%
%under $\mathrm{P}^{(n)}_{\varthetab_0;g_1}$
(where $\betab_0$ is the
matrix of eigenvectors associated with $\varthetab_0$), one easily
obtains that substituting $\hat\varthetab$ for $\varthetab_0$ in
(\ref{betagauss}) has no asymptotic impact on $Q_{\varthetab_0 ; \phi
_1 }^{(n)}$---more precisely,
$Q_{\hat\varthetab; \phi_1 }^{(n)}- Q_{\varthetab_0 ; \phi_1
}^{(n)}$ is $o_\mathrm{P}(1)$ as $\ny$ under $\mathrm
{P}^{(n)}_{\varthetab
_0;g_1}$, with $g_1 \in{\mathcal F}_{1}^{4} $. It follows that
$Q_{\hat\varthetab; \phi_1 }^{(n)}$ shares the same asymptotic
optimality properties as $ Q_{\varthetab_0 ; \phi_1 }^{(n)}$,
irrespective of the value of $\varthetab_0\in\mathcal{H}^\betab_0$.
Thus, a~locally and asymptotically \textit{most stringent} Gaussian test
of $\mathcal{H}^\betab_0$---denote it by $\phi^{(n)}_{\betab
;\mathcal{N}}$---can be based on the asymptotic chi-square
distribution [with $(k-1)$ degrees of freedom] of
%
%e5.3 ###
%
\begin{eqnarray} \label{QN}
Q_{\hat\varthetab; \phi_1 }^{(n)}
&=&
\frac{n a_{k}^{2}}{\hat{\sigma}^4}
\sum_{j=2}^k \bigl( \tilde\betab_j\pr\tilde\Vb^{-1/2}\mathbf{S}^{(n)}\tilde
\Vb^{-1/2}\betab^0
\bigr)^2
\nonumber\\
&=&
\frac{n a_{k}^{2}}{\hat{\sigma}^4 \hat\lambda_{1; \Vb}}
\sum_{j=2}^k \hat\lambda_{j; \Vb}^{-1} \bigl( \tilde\betab
_j\pr\mathbf{S}^{(n)}\betab^0
\bigr)^2 \\
&=&
\frac{n a_{k}^{2} | \Sb^{(n)}|^{2/k}}{\hat{\sigma}^4 \lambda_{1;
\Sb}}
\sum_{j=2}^k \lambda_{j; \Sb}^{-1} \bigl( \tilde\betab_j\pr
\mathbf{S}^{(n)}\betab^0
\bigr)^2
=:Q_{{\mathcal N}}^{(n)}.\nonumber
\end{eqnarray}
Since $\hat{\sigma}^{2}/| \Sb^{(n)}|^{1/k}$ converges to $a_{k}$ as
$\ny$ under the null ${\mathcal H}_{0}^{\betab}$ and Gaussian
densities, one can equivalently use the statistic
\[
\bar{Q}_{\mathcal{N}}^{(n)}:=\frac{ n}{\lambda_{1;\mathbf{S}}} \sum
_{j=2}^{k} { \lambda}_{j;\mathbf{S}} ^{-1} \bigl(\tilde{\betab
}_{j}\pr\mathbf{S}^{(n)}\betab_{}^{0} \bigr)^{2},
\]
which, of course, is still a locally and asymptotically \textit{most
stringent} Gaussian test statistic. For results on local powers, we
refer to Proposition~\ref{pseudogausstestbeta}.

This test is valid under Gaussian densities only (more precisely, under
radial densities with Gaussian kurtosis). On the other hand, it remains
valid in case Assumption~\ref{assuA} is weakened [as in \citet{A63} and
Tyler (\citeyear{T81}, \citeyear{T83})] into Assumption~\ref{assuApr1}. Indeed, the consistency
of $ \tilde{\Sigb}$ remains unaffected under the null, and $ \betab
_{}^{0}$ still is an eigenvector for both $ \tilde{\Sigb}^{-1/2}$
and $\Sigb$, so that $[\mathbf{I}_k-\betab_{}^{0}\betab^{0\prime} ]
\tilde{\Sigb}^{-1/2}\Sigb\tilde{\Sigb}^{-1/2} \betab_{}^{0} =\mathbf
{0}$. Hence,
\begin{eqnarray*}
Q_{\mathcal{N}}^{(n)}
&=&
n a_{k}^{2}
\sum_{j=2}^k \bigl( \tilde\betab_j\pr\tilde{\Sigb}^{-1/2}\mathbf
{S}^{(n)}\tilde{\Sigb}^{-1/2}\betab^0
\bigr)^2 \\
&=&n a_{k}^{2} \betab_{}^{0\prime} \tilde{\Sigb}^{-1/2}\mathbf{S}^{(n)}
\tilde{\Sigb}^{-1/2} \Biggl(\sum_{j=2}^k\tilde\betab
_j\tilde\betab_j\pr\Biggr)
\tilde{\Sigb}^{-1/2}\mathbf{S}^{(n)} \tilde{\Sigb}^{-1/2} \betab
_{}^{0} \\
&=&n a_{k}^{2} \betab_{}^{0\prime} \tilde{\Sigb}^{-1/2}\mathbf{S}^{(n)}
\tilde{\Sigb}^{-1/2} [\mathbf{I}_k-\betab_{}^{0}\betab
^{0\prime} ] \tilde{\Sigb}^{-1/2}\mathbf{S}^{(n)} \tilde{\Sigb
}^{-1/2} \betab_{}^{0} \\
&=&n a_{k}^{2} \betab_{}^{0\prime} \tilde{\Sigb}^{-1/2}\bigl(\mathbf
{S}^{(n)}- a_{k}^{-1} \Sigb\bigr) \tilde{\Sigb}^{-1/2} \\
& &{}\times
[\mathbf{I}_k-\betab_{}^{0}\betab^{0\prime} ] \tilde{\Sigb
}^{-1/2}\bigl(\mathbf{S}^{(n)}- a_{k}^{-1}\Sigb\bigr) \tilde{\Sigb}^{-1/2} \betab
_{}^{0} \\
&=&n a_{k}^{2} \betab_{}^{0\prime}{\Sigb}^{-1/2}\bigl(\mathbf{S}^{(n)}-
a_{k}^{-1}\Sigb\bigr){\Sigb}^{-1/2} \\
& &{}\times
[\mathbf{I}_k-\betab_{}^{0}\betab^{0\prime}]{\Sigb}^{-1/2}\bigl(\mathbf
{S}^{(n)}- a_{k}^{-1} \Sigb\bigr){\Sigb}^{-1/2} \betab_{}^{0} +o_\mathrm
{P}(1),
\end{eqnarray*}
as $n\rightarrow\infty$ under ${\mathcal H}_{0;1}^{\betab\prime}$.
Since $n^{1/2} a_{k} {\Sigb}^{-1/2}(\mathbf{S}^{(n)}- a_{k}^{-1} \Sigb
){\Sigb}^{-1/2} \betab^{0}$ is asymptotically $\mathcal
{N}(\mathbf{0},\mathbf{I}_{k} + \betab^{0} \betab^{0 \prime} )$ as $\ny$
under ${\mathcal H}_{0;1}^{\betab\prime}$ and Gaussian densities,
this idempotent quadratic form remains asymptotically chi-square, with
$(k-1)$ degrees of freedom, even when~\ref{assuA} is weakened into~\ref{assuApr1},
as was to be shown.\looseness=-1

This test is also invariant under the group of
transformations $\mathcal{G}_{\mathrm{rot},\circ}$ mapping $(\Xb
_1,\ldots
,\Xb_n)$ onto $(\mathbf{O}\Xb_1+\mathbf{t},\ldots,\mathbf{O}\Xb
_n+\mathbf{t})$,
where $\mathbf{t}$ is an arbitrary $k$-vector and $\mathbf{O} \in
{\mathcal
SO}_{k}^{\betab^{0}}:= \{ \mathbf{O} \in{\mathcal SO}_{k} \vert
\mathbf{O} \betab^{0}= \betab^{0} \}$, provided that the estimator of
$\varthetab_0$ used is equivariant under the same group---which the
esti\-mator $\hat\varthetab$ proposed in~(\ref{choiceprelim}) is. Indeed,
denoting by $Q_{\mathcal N}^{(n)}(\mathbf{O},\mathbf{t}) $, $\hat
\varthetab
(\mathbf{O},\mathbf{t})$, $\Lamb_{\Sb}(\mathbf{O},\mathbf{t})$, $\tilde
\betab
(\mathbf{O},\mathbf{t})$, $\tilde\Sigb(\mathbf{O},\mathbf{t})$ and $\Sb
^{(n)}(\mathbf{O},\mathbf{t})$ the statistics $Q_{\mathcal N}^{(n)}$,
$\hat
\varthetab$, $\Lamb_{\Sb}$, $\tilde\betab$, $\tilde\Sigb$ and
$\Sb^{(n)}$
computed from the transformed sample $(\mathbf{O} \Xb_{1}^{(n)}+\mathbf{t},
\ldots, \mathbf{O} \Xb_{n}^{(n)}+\mathbf{t})$, one easily checks that, for
any $\mathbf{O}\in{\mathcal SO}_{k}^{\betab^0}$, $\Lamb_{\Sb}(\mathbf
{O},\mathbf{t})=\Lamb_{\Sb}$, $\tilde\betab(\mathbf{O},\mathbf{t}) =
\mathbf{O} \tilde\betab$, $\tilde\Sigb(\mathbf{O},\mathbf{t}) = \mathbf
{O}{\tilde
\Sigb} \mathbf{O}\pr$ and $\Sb^{(n)}(\mathbf{O},\mathbf{t}) = \mathbf
{O}\Sb
^{(n)}\mathbf{O}\pr$, so that (noting that $\mathbf{O}\pr\betab^0=\betab^0$)
\[
Q_{\mathcal N}^{(n)}(\mathbf{O},\mathbf{t})
=
n a_{k}^{2} \sum_{j=2}^k \bigl(\tilde\betab_j\pr\mathbf{O}\pr\mathbf{O}\tilde
\Sigb^{-1/2} \mathbf{O}\pr\mathbf{O}\mathbf{S}^{(n)}\mathbf{O}\pr\mathbf
{O}\tilde\Sigb^{-1/2} \mathbf{O}\pr\betab^0
\bigr)^2
=Q_{\mathcal N}^{(n)}.
\]

Finally, let us show that $Q_{\mathrm{Anderson}}^{(n)}$ and
$Q_{\mathcal{N}}^{(n)}$ asymptotically coincide, under ${\mathcal
H}_{0;1}^{\betab\prime}$ and Gaussian densities, hence also under
contiguous alternatives. This asymptotic equivalence indeed is not a
straightforward consequence of the definitions~(\ref{AndTesta}) and
(\ref{QN}). Since $\sum_{j=2}^{k} {\lambda}_{j; \mathbf{S}}^{-1}
(\betab_{j; \Sb}\betab_{j; \mathbf{S}}\pr- \tilde\betab_j\tilde
\betab_j\pr)$ is $o_\mathrm{P}(1)$ and $n^{1/2}(\mathbf{S}^{(n)}- \tilde
\betab_0{ \Lamb}_\mathbf{S}\tilde\betab_0\pr)$ is $O_\mathrm{P}(1)$ as
$n\to\infty$, under ${\mathcal H}_{0;1}^{\betab\prime}$
and Gaussian densities [with ${\Lamb}_\mathbf{S}:=\operatorname
{diag}({\lambda
}_{1; \mathbf{S}},\ldots,{\lambda}_{k; \mathbf{S}})$], it follows from
Slutsky's lemma that
\begin{eqnarray*}
Q_{\mathrm{Anderson}}^{(n)}
&=&
\frac{n}{{\lambda}_{1; \mathbf{S}}} \sum_{j=2}^{k} {\lambda}_{j;
\mathbf{S}}^{-1} [ ({\lambda}_{j; \mathbf{S}}-{\lambda}_{1;\mathbf{S}})
{\betab}_{j; \mathbf{S}}\pr\betab_{}^{0} ]^{2} \\
&=&
\frac{n}{{\lambda}_{1; \mathbf{S}}} \sum_{j=2}^{k} {\lambda}_{j;
\mathbf{S}}^{-1} \bigl[ {\betab}_{j; \mathbf{S}}\pr\bigl(\mathbf{S}^{(n)}-
\tilde\betab_0{ \Lamb}_\mathbf{S}\tilde\betab_0\pr\bigr) \betab
_{}^{0} \bigr]^{2}
\\
&=&
\frac{n}{{\lambda}_{1; \mathbf{S}}} \sum_{j=2}^{k} {\lambda}_{j;
\mathbf{S}}^{-1} \bigl[ \tilde{\betab}_{j}\pr\bigl(\mathbf{S}^{(n)}- \tilde
\betab_0{ \Lamb}_\mathbf{S}\tilde\betab_0\pr\bigr) \betab_{}^{0}
\bigr]^{2}+o_\mathrm{P}(1) \\
&=&
\frac{n}{{\lambda}_{1; \mathbf{S}}} \sum_{j=2}^{k} {\lambda}_{j;
\mathbf{S}}^{-1} \bigl( \tilde{\betab}_{j}\pr\mathbf{S}^{(n)}\betab
_{}^{0} \bigr)^{2}+o_\mathrm{P}(1) \\
&=& \bar{Q}_{\mathcal{N}}^{(n)}+ o_\mathrm{P}(1)
\end{eqnarray*}
as $\ny$, still under ${\mathcal H}_{0;1}^{\betab\prime}$
and Gaussian densities. The equivalence between $Q_{\mathrm
{Anderson}}^{(n)}$ and $Q_{\mathcal{N}}^{(n)}$ in the Gaussian case
then follows since $\bar{Q}_{\mathcal{N}}^{(n)}={Q}_{\mathcal
{N}}^{(n)}+ o_\mathrm{P}(1)$ as $\ny$, under ${\mathcal H}_{0;1}^{\betab
\prime}$ and Gaussian densities.

%s5.2 ###
\subsection{Optimal Gaussian tests for eigenvalues}\label{Gausslamb}

Turning to ${\mathcal H}_{0}^{\Lamb}$, we now consider the Gaussian
version of the test statistic $T_{\varthetab_0 ; f_{1}}^{(n)}$
obtained in Section~\ref{eigenvaluesprob}. In view of~(\ref{Tparam}),
we have
%
%e5.4 ###
%
\begin{equation}\label{Tgauss}
T_{\varthetab_0; \phi_1 }^{(n)}
= n^{1/2}( a_{p,q}(\Lamb_{0} ))^{-1/2} \mathbf{c}_{p,q}\pr\dvec
\bigl(\Lamb_{0}^{1/2}\betab\pr\Sb_{\varthetab_0;\phi_1 }^{(n)}\betab
\Lamb_{0}^{1/2}\bigr)
\end{equation}
%
%{eqnarray}
[recall that $\mathcal{J}_k(\phi_1 ) =k(k+2)$; see (\ref
{normalscores})]. Here also we have to estimate $\varthetab_0$ in
order to obtain a genuine test statistic. By using the fact
that $\betab\Lamb_{0}\betab\pr=\Vb_0$ (where all parameter values
refer to those in $\varthetab_0$), we obtain that,
in~(\ref{Tgauss}),
%
%e5.5 ###
%
\begin{eqnarray}
&&
n^{1/2}\mathbf{c}_{p,q}\pr\dvec\bigl(\Lamb_{0}^{1/2}\betab\pr\Sb
_{\varthetab_0;\phi_1
}^{(n)}\betab\Lamb_{0}^{1/2}\bigr)\nonumber\\[-8pt]\\[-8pt]
%%[2mm]
&&\qquad= \frac{n^{1/2} a_{k}}{\sigma^{2}} \mathbf{c}_{p,q}\pr\dvec
\Biggl(\betab\pr\frac{1}{n}\sum_{i=1}^n(\Xb_i-{\bolds\theta})(\Xb
_i-{\bolds\theta})\pr
\betab\Biggr),\nonumber
\end{eqnarray}
a $O_\mathrm{P}(1)$ expression which does not depend on $\Lamb_{0}$. In
view of Lemma~\ref{parametricasymplin} and the block-diagonal form of
the information matrix, estimation of ${\bolds\theta}$, $\sigma^2$
and $\betab$ has no asymptotic impact on the eigenvalue part $\Deltab
_{\varthetab; \phi_1 }^{\III}$ of the central sequence, hence on~$
T_{\varthetab_{0}; \phi_1 }^{(n)}$.\vspace*{2pt} As for $a_{p,q}(\Lamb_0 )$, it
is a continuous function of $\Lamb_0$, so that, in view of Slutsky's
lemma, plain consistency of the estimator of $\Lamb_0$ is sufficient.
Consequently, we safely can use here the unconstrained estimator
%
%e5.6 ###
%
\begin{equation}\label{thetachapeaulamb}
\hat{\varthetab}:=
(\bar{\mathbf X}\pr, \hat{\sigma}^{2},
(\dvecrond\hat\Lamb_{\Vb})',
(\vecop{\hat\betab_{\Vb}})\pr)\pr;
\end{equation}
see the beginning of Section~\ref{Gaussbetab}. Using again the fact
that, under Gaussian densities, $\hat{\sigma}^{2}/| \Sb
^{(n)}|^{1/k}$ converges to $a_{k}$ as $\ny$, a~locally and asymptotically \textit{most powerful} Gaussian test
statistic therefore is given by
%
%e5.7 ###
%
\begin{eqnarray}\label{commeAnderson}
T_{\mathcal{N}}^{(n)}
:\!&=&
\frac{n^{1/2} a_{k}}{\hat{\sigma}^{2}} ({a}_{p,q}(\hat{\Lamb
}_{\mathbf{V}}))^{-1/2} \mathbf{c}_{p,q}\pr\dvec\bigl(\hat\betab_{\Vb
}\pr{\Sb}^{(n)}\hat\betab_{\Vb}\bigr)\nonumber\\
&=&
\frac{n^{1/2} a_{k}|\Sb^{(n)}|^{1/k}}{\hat{\sigma}^{2}}
({a}_{p,q}({\Lamb}_{\mathbf{S}}))^{-1/2} \mathbf{c}_{p,q}\pr\dvec
\bigl(\hat\betab_{\Vb}\pr{\Sb}^{(n)}\hat\betab_{\Vb}\bigr)\nonumber\\[-8pt]\\[-8pt]
&=&n^{1/2} (a_{p,q}({\Lamb}_\mathbf{S}))^{-1/2}
\Biggl(
(1-p) \sum_{j=q+1}^k \lambda_{j;\Sb}-p \sum_{j=1}^q \lambda
_{j;\Sb}
\Biggr)\nonumber\\
&&{} +o_\mathrm{P}(1),\nonumber
\end{eqnarray}
under Gaussian densities as $\ny$.
The corresponding test, $\phi^{(n)}_{\Lamb; \mathcal{N}}$ say,
rejects $\mathcal{H}^\Lamb_0$ whenever $T_{\mathcal{N}}^{(n)}$ is
smaller than the standard normal $\alpha$-quantile; (\ref
{commeAnderson}) shows that $T_{\mathcal{N}}^{(n)}$ coincides [up to
$o_\mathrm{P}(1)$] with $T_{\mathrm{Anderson}}^{(n)}$ given in (\ref
{TAnd}), which entails that (i) $\phi^{(n)}_{\Lamb; \mathrm
{Anderson}}$ is also locally and asymptotically most powerful under
Gaussian densities, and that (ii) the validity of $\phi^{(n)}_{\Lamb;
\mathcal{N}}$ extends to $\mathcal{H}^{\Lamb\prime\prime} _{0;q}$
(since the validity of $\phi^{(n)}_{\Lamb; \mathrm{Anderson}}$ does).

%%%%%%%%%%%%%%%%%%%%%%%%%%%%%%%%%%%%%%%%%%%%%%%%%%%%%%%%%%%%%%%%%%%%%%%%%%%%%%%%%%%%%%%%%%%%%%%%%%%%%%%%%%%%%%%%%%%%%%%%%%%%%%%%%%%%%%%%%%%%%%%%%%%%%%%%%%%%
%%%%%%%%%%%%%%%%%%%%%%%%%%%%%%%%%%%%%%%%%%%%%%%%%%%%%%%%%%%%%%%%%%%%%%%%%%%%%%%%%%%%%%%%%%%%%%%%%%%%%%%%%%%%%%%%%%%%%%%%%%%%%%%%%%%%%%%%%%%%%%%%%%%%%%%%%%%%

%s5.3 ###
\subsection{Optimal pseudo-Gaussian tests for eigenvectors}
\label{pseudovec}

The Gaussian tests $\phi^{(n)}_{\betab; \mathcal{N}}$ and $\phi
^{(n)}_{\Lamb; \mathcal{N}}$ of Sections~\ref{Gaussbetab} and
\ref{Gausslamb}
unfortunately are valid under multinormal densities only (more
precisely, as we shall see, under densities with Gaussian kurtosis). It
is not difficult, however, to extend their validity to the whole class
of elliptical populations with finite fourth-order moments, while
maintaining their optimality properties at the multinormal.

Let us first introduce the following notation. For any $g_{1}\in
\mathcal{F}_1^{4}$, let (as in Lem\-ma~\ref{parametricasymplin}) $D_k(g_{1}):=\mu
_{k+1;g_1}/\mu_{k-1;g_1}
=\sigma^{-2}
\mathrm{E}_{\varthetab;g_{1}} [d^2_{i}(\tetb,\Vb)]=\int_0^1({\tilde
G}{}^{-1}_{1k}(u))^2 \,du$
and
$E_k(g_{1})
:= \sigma^{-4}
\mathrm{E}_{\varthetab;g_{1}} [d^4_{i}(\tetb,\Vb)]=\int_0^1({\tilde
G}{}^{-1}_{1k}(u))^4 \,du$,
where ${\tilde G}_{1k}(r):= \break(\mu_{k-1;g_1})^{-1}
\int_0^r s^{k-1}
g_1(s) \,ds$; see Section~\ref{defelliptttt}. Then
\[
\kappa_k(g_{1}):= \frac{k}{k+2} \frac{E_k(g_{1})}{D_k^{2}(g_{1})}-1
\]
is the \textit{kurtosis} of the elliptic population with radial density
$g_{1}$ [see, e.g., page 54 of \citet{A03}]. For Gaussian densities,
$E_k(\phi_1 ) = k(k+2)/a_{k}^{2}$, $D_k(\phi_{1})=k/a_{k}$ and
$\kappa_k(\phi_1 )=0$.

Since the asymptotic covariance matrix of $\Deltab_{{\varthetab} ;
\phi_1 }^{\IV}$ under $\mathrm{P}^{(n)}_{\varthetab;g_{1}}$ (with
$\varthetab\in\mathcal{H}^\betab_0$ and $g_{1} \in\mathcal
{F}_{1}^{4}$) is $ (a_{k}^{2}E_{k}(g_{1})/k(k+2)) \Gamb^{\IV
}_{{\varthetab} ;\phi_1}$,
it is natural to base our pseudo-Gaussian tests on statistics of the
form [compare with the $f_1=\phi_1$ version of~(\ref{betagauss})]
\begin{eqnarray*}
Q_{\varthetab_0, \mathcal{N}*}^{(n)}
:\!&=&
\frac{k(k+2)}{a_{k}^{2}E_{k}(g_{1})} \Deltab_{{\varthetab_0} ;
\phi_1 }^{\IV\prime}
[ (\Gamb^{\IV}_{\varthetab_0 ;\phi_1})^{-} -
\mathbf{P}_{k}^{\betab_{0}} ((\mathbf{P}_{k}^{\betab_{0}}) \pr
\Gamb^{\IV}_{\varthetab_0 ; \phi_1} \mathbf{P}_{k}^{\betab_{0}}
)^{-} (\mathbf{P}_{k}^{\betab_{0}}) \pr]
\Deltab_{{\varthetab_0} ; \phi_1 }^{\IV}
\\
& = &
\bigl(1+ \kappa_{k}(g_{1})\bigr)^{-1} \frac{k^{2}}{D_{k}^{2}(g_{1})
a_{k}^{2}} Q_{\varthetab_0, \phi_{1}}^{(n)}
=:\bigl(1+ \kappa_{k}(g_{1})\bigr)^{-1} Q_{\varthetab_{0}, \mathcal{N}}(g_{1}).
\end{eqnarray*}
As in the Gaussian case, and with the same $\hat\varthetab$ as
in~(\ref{choiceprelim}), Lemma~\ref{parametricasymplin} entails that
$
Q_{\hat\varthetab, \phi_{1}}^{(n)}= Q_{\varthetab_0, \phi
_{1}}^{(n)}+ o_\mathrm{P}(1),
$
as $\ny$ under $\mathrm{P}^{(n)}_{\varthetab_{0}; g_{1}}$, with
$\varthetab_{0} \in{\mathcal H}_{0}^{\betab}$ and $g_{1} \in
{\mathcal F}_{1}^{4}$.
Since $\hat{\sigma}^{2}/|\mathbf{S}^{(n)}|^{1/k}$ consistently estimates
$k/D_k(g_1)$ under $\mathrm{P}^{(n)}_{\varthetab_0; g_{1}}$,
with $\varthetab_0 \in{\mathcal H}_{0}^{\betab}$ and $g_{1} \in
{\mathcal F}_{1}^{4}$,
it follows from Slutsky's lemma that
\[
\hat{Q}_{\hat{\varthetab}, \mathcal{N}} :=
\frac{\hat{\sigma}^{2}}{|\mathbf{S}^{(n)}|^{1/k} a_{k}^{2}} Q_{\hat
\varthetab, \phi_{1}}^{(n)}
\]
satisfies $\bar{Q}_{\mathcal{N}}^{(n)}= \hat{Q}_{\hat{\varthetab},
\mathcal{N}}= Q_{\varthetab_0, \mathcal{N}}(g_{1})+o_\mathrm{P}(1)$ as
$\ny$, still under $ {\mathcal H}_{0}^{\betab}$, $g_{1} \in{\mathcal
F}_{1}^{4}$. The pseudo-Gaussian test $\phi^{(n)}_{\betab;\mathcal
{N}*}$ we propose is based on
%
%e5.8 ###
%
\begin{equation}\label{QN*}
Q_{\mathcal{N}*}^{(n)} := (1+\hat\kappa_{k})^{-1} \bar{Q}_{\mathcal
{N}}^{(n)} ,
\end{equation}
where $\hat{\kappa}_k :=
(kn^{-1}\sum_{i=1}^n \hat{d}_{i}^4)/((k+2)(n^{-1}\sum_{i=1}^n \hat
{d}_{i}^2)^2)-1
$, with $\hat{d}_{i}:={d}_{i}(\bar{\Xb}, \mathbf{S}^{(n)})$.
The statistic $Q_{\mathcal{N}*}^{(n)}$ indeed remains asymptotically
chi-square [$(k-1)$ degrees of freedom] under $\mathcal{H}^{\betab
\prime} _{0;1}$ for any $g_{1} \in\mathcal{F}_{1}^{4}$.
Note that\vspace*{1pt} $\phi^{(n)}_{\betab;\mathcal{N}*}$ is obtained from $\phi
^{(n)}_{\betab;\mathcal{N}}$ by means of the standard kurtosis
correction of \citet{SB87}, and asymptotically coincides
with $\phi^{(n)}_{\betab;\mathrm{Tyler}}$; see~(\ref{TylTesta}).

Local powers for $\phi^{(n)}_{\betab;\mathcal{N}*}$ classically
follow from applying Le Cam's third lemma.
Let $\taub^{(n)}:=((\taub^{\I(n)})\pr, \tau^{\II(n)}, (\taub
^{\III(n)})\pr, (\taub^{\IV(n)})\pr)\pr, $
with $\taub^{(n)\prime}\taub^{(n)}$ uniformly bounded, where
$ \taub^{\IV(n)}
=\vecop(\mathbf{b}^{(n)})$ is a perturbation of $\operatorname{vec}(\betab_0)=\vecop
(\betab^{0},\betab_2,\ldots, \betab_k)$ such that $\betab_0\pr
\mathbf{b}$, with $\mathbf{b}= ( \mathbf{b}_{1}^{\prime}, \ldots,
\mathbf{b}_{k}^{\prime} )\pr:=\lim_{\ny}\mathbf{b}^{(n)}$, is
skew-symmetric; see~(\ref{antisym}) and~(\ref{tangentbetab}). Assume
furthermore that the corresponding perturbed value of $\varthetab_0\in
\mathcal{H}^\betab_0$ does not belong to $\mathcal{H}^\betab_0$,
that is, $\mathbf{b}_1\neq\mathbf{0}$, and define
%
%e5.9 ###
%
\begin{eqnarray}\label{rbeta}
r_{\varthetab_0 ; \taub}^{\betab}
:\!&=& \lim_{\ny} \bigl(\operatorname{vec} \mathbf{b}^{(n)}\bigr)
^\prime\mathbf{G}_{k}^{\betab_{0}} \operatorname{diag}\bigl(\nu_{12}^{-1},
\ldots
, \nu{}^{-1}_{1k}, \mathbf{0}_{1\times{(k-2)(k-1)}/{2} }\bigr)\nonumber\\
&&\hspace*{20.4pt}{}\times
(\mathbf{G}_{k}^{\betab_{0}})\pr\bigl(\operatorname{vec} \mathbf
{b}^{(n)}\bigr)\\
&=&
4 \sum_{j=2}^{k} \nu_{1j}^{-1} (\betab\pr_{j} \mathbf{b}_1)^{2}.\nonumber
\end{eqnarray}
The following result summarizes the asymptotic properties of the
pseudo-Gaussian tests $\phi^{(n)}_{\betab; \mathcal{N}*}$.
Note that optimality issues involve $\mathcal{H}^\betab_0$ [hence
require Assumption~\ref{assuA}], while validity extends to $\mathcal
{H}^{\betab\prime} _{0;1}$ [which only requires Assumption~\ref{assuApr1}].
\begin{Prop} \label{pseudogausstestbeta}
\textup{(i)} $Q^{(n)}_{\mathcal{N}*}$ is asymptotically chi-square with $(k-1)$
degrees of freedom under $\bigcup_{\varthetab\in{\mathcal
H}_{0;1}^{\betab\prime}}
\bigcup_{g_{1}\in\mathcal{F}_{1}^{4} }
\{ \mathrm{P}^{(n)}_{\varthetab;g_{1}}\}$, and asymptotically noncentral
chi-square, still with $(k-1)$ degrees of
freedom, but with noncentrality parameter
$r_{\varthetab; \taub}^{\betab}/4(1+\kappa_k(g_1)) $ under $\mathrm
{P}^{(n)}_{\varthetab+n^{-1/2}\taub^{(n)};g_{1}}$, with $\varthetab
\in\mathcal{H}^\betab_0$,
$g_{1}\in\mathcal{F}_{a}^{4}$, and $\taub^{(n)}$ as
described above;\vadjust{\goodbreak}

{\smallskipamount=0pt
\begin{longlist}[(iii)]
\item[(ii)]
the sequence of tests $\phi^{(n)}_{\betab; \mathcal{N}*}$
rejecting the null whenever $Q^{(n)}_{\mathcal{N}*}$ exceeds the
$\alpha$ upper-quantile $\chi^2_{k-1;1-\alpha}$ of the chi-square
distribution with $(k-1)$ degrees of freedom has asymptotic size
$\alpha$ under $\bigcup_{\varthetab\in
{\mathcal H}_{0;1}^{\betab\prime} }\bigcup_{g_{1}\in\mathcal
{F}_{1}^{4} } \{ \mathrm{P}^{(n)}_{\varthetab;g_{1}}\}$;

\item[(iii)]
the pseudo-Gaussian tests $\phi^{(n)}_{\betab; \mathcal{N}*}$ are
asymptotically equivalent, under $\bigcup_{\varthetab\in{\mathcal
H}_{0;1}^{\betab\prime} }\{\mathrm{P}^{(n)}_{\varthetab; \phi_1 }\}$
and contiguous alternatives, to the optimal parametric Gaussian tests
$\phi^{(n)}_{\betab; \mathcal{N}}$; hence, the sequence $\phi
^{(n)}_{\betab; \mathcal{N}*}$ is locally and asymptotically most stringent,
still at asymptotic level $\alpha$, for $\bigcup_{\varthetab\in
{\mathcal H}_{0;1}^{\betab\prime}}\bigcup_{g_{1}\in\mathcal
{F}_{1}^{4} } \{ \mathrm{P}^{(n)}_{\varthetab;g_{1}}\}$ against
alternatives of the form
$\bigcup_{\varthetab\notin{\mathcal H}_{0}^{\betab}} \{ \mathrm
{P}^{(n)}_{\varthetab;\phi_1 }\}$.
\end{longlist}}
\end{Prop}

Of course, since $\hat\kappa_{k}$ is invariant under $\mathcal
{G}_{\mathrm{rot},\circ}$, the pseudo-Gaussian test inherits the
$\mathcal
{G}_{\mathrm{rot},\circ}$-invariance features of the Gaussian one.

%s5.4 ###
\subsection{Optimal pseudo-Gaussian tests for eigenvalues}\label{pseudoval}
As in the previous section, the asymptotic null distribution of the
Gaussian test statistic $T_{\mathcal{N}}^{(n)}$ is not standard normal
anymore under radial density $g_1$ as soon as $\kappa_k(g_1)\neq
\kappa_k(\phi_1)$. The Gaussian test $\phi^{(n)}_{\Lamb; \mathcal
{N}}$ thus is\vspace*{1pt} not valid (does not have asymptotic level $\alpha$)
under such densities.
The same reasoning as before leads to a similar kurtosis correction,
yielding a pseudo-Gaussian test statistic
\[
T_{\mathcal{N}{*}}^{(n)}:=(1+\hat\kappa_k)^{-1/2} \tilde
{T}_{\mathcal{N}}^{(n)},
\]
where $\tilde{T}_{\mathcal{N}}^{(n)}:=n^{1/2}(a_{p,q}({\Lamb}_\mathbf
{S}))^{-1/2}
(
(1-p) \sum_{j=q+1}^k \lambda_{j;\Sb}-p\sum_{j=1}^q \lambda_{j;\Sb
})$ and $\hat{\kappa}_k$ is as in Section~\ref{pseudovec}. This
statistic %clearly
coincides with $T_{\mathrm{Davis}}^{(n)}$ given in~(\ref{TDav}). % of

Here also,\vspace*{1pt} local powers are readily obtained via Le Cam's third lemma.
Let $\taub^{(n)}:=((\taub^{\I(n)})^\prime, \tau^{\II(n)}, (\taub
^{\III(n)})^\prime, (\taub^{\IV(n)})^\prime)\pr$, with $\taub
^{(n)\prime}\taub^{(n)}$ uniformly
bounded, where $\taub^{\III(n)}:= \dvecrond(\mathbf{l}^{(n)})$ is such
that $\mathbf{l}:=\lim_{\ny}\mathbf{l}^{(n)}:=\lim_{\ny}
\operatorname{diag}(\ell^{(n)}_{1}, \break \ldots, \ell^{(n)}_{k})$ satisfies
$\operatorname{tr}(\Lamb_{\Vb}^{-1} \mathbf{l})=0$ [see~(\ref{ell1})
and the comments thereafter], and define
%
%e5.10 ###
%
\begin{equation}
\label{rLamb}\qquad
r_{\varthetab; \taub}^{\Lamb_{\Vb}} := \lim_{\ny}\operatorname
{grad} h(\dvecrond(\Lamb_{\Vb}))\pr\taub^{\III(n)}
=
(1-p) \sum_{j=q+1}^k \mathbf{l}_{j} -p\sum_{j=1}^q \mathbf{l}_{j}.
\end{equation}
The following proposition summarizes the asymptotic properties of the
resulting pseudo-Gaussian tests $\phi^{(n)}_{\Lamb_\Vb;\mathcal{N}*}$.
\begin{Prop} \label{pseudogausstestlambda}
\textup{(i)} $T^{(n)}_{\mathcal{N}*}$ is asymptotically normal, with mean
zero
under $\bigcup_{\varthetab\in{\mathcal H}_{0;q}^{\Lamb\prime
\prime} }
\bigcup_{g_{1}\in\mathcal{F}_{1}^{4} }
\{ \mathrm{P}^{(n)}_{\varthetab;g_{1}}\}$, mean
$
({4a_{p,q}(\Lamb_\Vb) (1+ \kappa_{k}(g_{1}))}^{-1/2} r_{\varthetab;
\taub}^{\Lamb_{\Vb}}
$
under $\mathrm{P}^{(n)}_{\varthetab+n^{-1/2}\taub^{(n)};g_{1}}$,
$\varthetab\in
{\mathcal H}_{0}^{\Lamb}$,
$g_{1}\in\mathcal{F}_a^4$ and $ \taub^{(n)}$ as described
above, and variance one under both;\vadjust{\goodbreak}

{\smallskipamount=0pt
\begin{longlist}[(iii)]
\item[(ii)] the sequence of tests $\phi^{(n)}_{\Lamb; \mathcal{N}*}$
rejecting the null whenever $T^{(n)}_{\mathcal{N}*}$ is less than the
standard normal $\alpha$-quantile $z_\alpha$ has asymptotic size
$\alpha$ under\break $\bigcup_{\varthetab\in
{\mathcal H}_{0;q}^{\Lamb\prime\prime}}\bigcup_{g_{1}\in
\mathcal{F}_{1}^{4} } \{ \mathrm{P}^{(n)}_{\varthetab;g_{1}}\}$;

\item[(iii)] the pseudo-Gaussian tests $\phi^{(n)}_{\Lamb; \mathcal{N}*}$
are asymptotically equivalent, under $\bigcup_{\varthetab\in
{\mathcal H}_{0;q}^{\Lamb\prime\prime} }\{\mathrm{P}^{(n)}_{\varthetab
; \phi_1 }\}$ and contiguous alternatives, to the optimal parametric
Gaussian tests $\phi^{(n)}_{\Lamb; \mathcal{N}}$; hence, the
sequence $\phi^{(n)}_{\Lamb; \mathcal{N}*}$ is locally and
asymptotically most powerful,
still at asymptotic level $\alpha$, for $\bigcup_{\varthetab\in
{\mathcal H}_{0;q}^{\Lamb\prime\prime} }\bigcup_{g_{1}\in\mathcal
{F}_{1}^{4} } \{ \mathrm{P}^{(n)}_{\varthetab;g_{1}}\}$ against
alternatives of the form
$\bigcup_{\varthetab\notin{\mathcal H}_{0}^{\Lamb}} \{ \mathrm
{P}^{(n)}_{\varthetab;\phi_1 }\}$.
\end{longlist}}
\end{Prop}

%s6 ###
\section{Rank-based tests for principal components} \label{ranktests}

%s6.1 ###
\subsection{Rank-based statistics: Asymptotic representation and
asymptotic normality}\label{rankHajek}
The parametric tests proposed in Section~\ref{paramtests} are valid
under specified radial densities $f_1$ only, and therefore are of
limited practical value. The importance of the Gaussian tests of
Sections~\ref{Gaussbetab} and
\ref{Gausslamb} essentially follows from the fact that they
belong to usual practice, but Gaussian assumptions are quite
unrealistic in most applications. The pseudo-Gaussian procedures of
Sections~\ref{pseudovec} and~\ref{pseudoval} are more appealing, as
they only require finite fourth-order moments. Still, moments of order
four may be infinite and, being based on empirical covariances,
pseudo-Gaussian procedures remain poorly robust. A~straightforward idea
would consist in robustifying them by substituting some robust estimate
of scatter for empirical covariance matrices. This may take care of
validity-robustness issues, but has a negative impact on powers, and
would not achieve efficiency-robustness. The picture is quite different
with the rank-based procedures we are proposing in this section. While
remaining valid under completely arbitrary radial densities, these
methods indeed also are efficiency-robust; when based on Gaussian
scores, they even uniformly outperform, in the Pitman sense, their
pseudo-Gaussian counterparts (see Section~\ref{secare}). Rank-based
inference, thus, in this problem as in many others, has much to offer,
and enjoys an extremely attractive combination of robustness and
efficiency properties.

The natural framework for principal component analysis actually is the
semiparametric context of elliptical families in which ${\bolds\theta
}$, $\dvecrond(\Lamb_\Vb)$, and $\betab$ (not $\sigma^2$) are the
parameters of interest, while the radial density $f$ [equivalently, the
couple $(\sigma^2, f_1)$] plays the role of an infinite-dimensional
nuisance. This semiparametric model enjoys the double structure
considered in \citet{HW03}, which allows for efficient
rank-based inference: the fixed-$f_1$ subexperiments, as shown in
Proposition~\ref{LAN} are ULAN, while the fixed-$({\bolds\theta}$,
$\dvecrond(\Lamb_\Vb)$, $\betab)$ subexperiments [equivalently, the
fixed-$({\bolds\theta}, \Vb)$ subexperiments] are generated by
groups of transformations acting on the observation space. Those groups
here are of the form\vadjust{\goodbreak}
$\mathcal{G}_{{\bolds\theta}, \Vb}^{(n)}, \sirc$ and consist of
the \textit{continuous monotone radial transformations}
${\mathcal{G}} ^{(n)}_h$
\begin{eqnarray*}
{\mathcal{G}}^{(n)}_h(\Xb_{1}, \ldots, \Xb
_{n})&=&{\mathcal{G}}^{(n)}_h \bigl(
{\bolds\theta}+ d_{1}({\bolds\theta}, \Vb) \Vb^{1/2} \mathbf
{U}_{1}({\bolds\theta}, \Vb), \ldots, \\
& &\hspace*{43pt}
{\bolds\theta}+ d_{n}({\bolds\theta}, \Vb) \Vb^{1/2} \mathbf
{U}_{n}({\bolds\theta}, \Vb) \bigr) \\
:\!&=& \bigl({\bolds\theta}+
h(d_{1}({\bolds\theta}, \Vb)) \Vb^{1/2} \mathbf{U}_{1}({\bolds\theta
}, \Vb), \ldots, \\
&&
\hspace*{25.4pt}{\bolds\theta}+ h(d_{n}({\bolds\theta}, \Vb)) \Vb^{1/2} \mathbf
{U}_{n}({\bolds\theta}, \Vb)\bigr),
\end{eqnarray*}
where $h\dvtx\R^{+} \rightarrow\R^{+}$ is continuous, monotone
increasing, and satisfies\break $\lim_{r \rightarrow\infty}h(r)=
\infty$ and $h(0)=0$. The group $\mathcal{G}_{{\bolds\theta},
\Vb}^{(n)}, \sirc$ generates the fixed-$({\bolds\theta}, \Vb)$ family
of distributions $ \bigcup_{\sigma^{2}} \bigcup_{f_{1}}\{\mathrm
{P}_{{\bolds\theta}, \sigma^{2}, \dvecronds(\Lamb_\Vb),\vecop (\betab) ;
f_{1}}^{(n)}\}$.\vspace*{1pt} The general results of \citet{HW03}
thus indicate that efficient inference can be based on the
corresponding maximal invariants, namely the vectors
\[
\bigl(R^{(n)}_{1}({\bolds\theta}, \Vb), \ldots,R^{(n)}_{n}({\bolds
\theta}, \Vb), \Ub_{1}({\bolds\theta}, \Vb), \ldots, \Ub
_{n}({\bolds\theta}, \Vb) \bigr)
\]
of ranks and multivariate signs, where $R^{(n)}_{i}({\bolds\theta},
\Vb)$ denotes
the rank of $d_{i}({\bolds\theta}, \Vb)$ among $d_{1}({\bolds\theta
}, \Vb), \ldots, d_{n}({\bolds\theta}, \Vb)$. Test
statistics based on such invariants automatically are distribution-free
under $\bigcup_{\sigma^{2}} \bigcup_{f_{1}}\{\mathrm{P}_{{\bolds\theta
}, \sigma^{2}, \dvecronds(\Lamb_\Vb),\vecop(\betab) ;
f_{1}}^{(n)}\}$.

Letting $R_{i}:=R_{i}({\bolds\theta}, \Vb)$ and $\Ub_{i}:=\Ub
_{i}({\bolds\theta}, \Vb)$, define
\[
{\utDelta}{}_{ \varthetab; K}^{\III}:=
\frac{1}{2\sqrt{n}}
{\Mb}_{k}^{\Lamb_{\Vb}} \mathbf{H}_{k} (
{\Lamb_{\Vb}^{-1/2}\betab\pr} )^{\otimes2}
\sum_{i=1}^{n} K \biggl(\frac{R^{(n)}_{i}}{n+1}
\biggr)\vecop( \mathbf{U}_{i}\mathbf{U}_{i}\pr)
\]
and
\[
{\utDelta}{}_{ \varthetab; K}^{\IV}:=
\frac{1}{2\sqrt{n}}
\mathbf{G}_{k}^{\betab} \mathbf{L}_{k}^{\betab,\Lamb_{\Vb}}
(\mathbf{V}^{\otimes2})^{-1/2} \sum_{i=1}^{n} K \biggl( \frac
{R^{(n)}_{i}}{n+1} \biggr) \vecop(\Ub_{i}\Ub_{i}\pr).
\]
Associated with ${\utDelta}{}_{ \varthetab; K}^{\III}$ and
${\utDelta}{}_{ \varthetab; K}^{\IV}$, let
\[
\ubDelta_{\varthetab; K, g_{1}}^{\III}:=
\frac{1}{2\sqrt{n}}
{\Mb}_{k}^{\Lamb_{\Vb}} \mathbf{H}_{k} (
{\Lamb_{\Vb}^{-1/2}\betab\pr} )^{\otimes2}
\sum_{i=1}^{n} K \biggl(\tilde{G}_{1k} \biggl(\frac
{d_{i}({\bolds\theta}, \Vb)}{\sigma} \biggr) \biggr)\vecop
(\mathbf{U}_{i}\mathbf{U}_{i}\pr)
\]
and
\[
\ubDelta_{\varthetab; K, g_{1}}^{\IV}:=
\frac{1}{2\sqrt{n}}
\mathbf{G}_{k}^{\betab} \mathbf{L}_{k}^{\betab, \Lamb_{\Vb}} (\mathbf
{V}^{\otimes2})^{-1/2} \sum_{i=1}^{n}K \biggl( \tilde{G}_{1k}
\biggl(\frac{d_{i}({\bolds\theta}, \Vb)}{\sigma} \biggr) \biggr) \vecop
(\Ub_{i}\Ub_{i}\pr),
\]
where $\tilde{G}_{1k}$ is as in Section~\ref{pseudovec}.
The following proposition provides an asymptotic representation and
asymptotic normality result for ${\utDelta}{}_{ \varthetab; K}^{\III
}$ and ${\utDelta}{}_{ \varthetab; K}^{\IV}$.

\begin{Prop} \label{Hajek}
Let Assumption~\ref{assuS} hold for the score function $K$. Then:
\begin{longlist}
\item
(\textit{asymptotic representation})\vspace*{-2pt}
$({\utDelta}{}_{ \varthetab; K}^{\III\prime},{\utDelta}{}_{
\varthetab; K}^{\IV\prime})\pr=(\ubDelta_{ \varthetab; K,
g_{1}}^{\III\prime}, \ubDelta_{ \varthetab; K, g_{1}}^{\IV\prime
})\pr+ o_{L^{2}}(1)$
as $\ny$, under $\mathrm{P}^{(n)}_{\varthetab;g_1}$, for any $\varthetab
\in\Thetab$ and $g_1\in\mathcal{F}_1$;\vspace*{2pt}

\item
(\textit{asymptotic normality}) let Assumption~\ref{assuA} hold and consider
a bound\-ed~sequence $\taub^{(n)}:=((\taub^{\I(n)})^{\prime}, \tau
^{\II(n)} , (\taub^{\III(n)})^{\prime} , (\taub^{\IV(n)})^{\prime
} )\pr$ such that both $\taub^\III:=\lim_{\ny}\taub^{\III(n)}$
and $\taub^\IV:=\lim_{\ny}\taub^{\IV(n)}$ exist. Then $(\ubDelta
_{ \varthetab; K, g_{1}}^{\III\prime}, \ubDelta_{ \varthetab; K,
g_{1}}^{\IV\prime} )\pr$
is asymptotically normal, with mean zero and mean
\[
\frac{\mathcal{J}_k(K,g_{1})}{k(k+2)}
\pmatrix{
\mathbf{D}_{k}(\Lamb_{\Vb}) \taub^{\III} \cr
\frac{1}{4} \mathbf{G}_{k}^{\betab} \operatorname{diag}\bigl(\nu_{12}^{-1},
\ldots, \nu_{(k-1)k}^{-1}\bigr) (\mathbf{G}_{k}^{\betab})\pr\taub^{\IV}}\vspace*{-2pt}
\]
[where $\mathcal{J}_k(K,g_{1})$ was defined in~(\ref{infoKf})], under
$\mathrm{P}^{(n)}_{\varthetab; g_{1}}$ (any $\varthetab\in\Thetab$ and
$g_1\in\mathcal{F}_1$) and $\mathrm{P}^{(n)}_{\varthetab+ n^{-1/2}
\taub^{(n)}; g_{1}}$ (any $\varthetab\in\Thetab$ and $g_1\in
\mathcal{F}_a$), respectively, and block-diagonal covariance matrix
$\operatorname{diag} (\Gamb_{\varthetab; K}^{\III},\Gamb_{\varthetab;
K}^{\IV} )$ under both, with
\[
\Gamb_{\varthetab; K}^{\III}:=\frac{\mathcal{J}_k(K)}{k(k+2)}
\mathbf{D}_{k}(\Lamb_{\Vb})\vspace*{-2pt}
\]
and
%
%e6.1 ###
%
\begin{equation}\label{GambK}
\Gamb_{\varthetab; K}^{\IV}:=\frac{\mathcal{J}_k(K)}{4k(k+2)}
\mathbf{G}_{k}^{\betab} \operatorname{diag}\bigl(\nu_{12}^{-1}, \ldots, \nu
_{(k-1)k}^{-1}\bigr)(\mathbf{G}_{k}^{\betab})\pr.\vspace*{-2pt}
\end{equation}
\end{longlist}
\end{Prop}

The proofs of parts (i) and (ii) of this proposition are entirely
similar to those of Lemma 4.1 and Proposition 4.1, respectively, in
\citet{HP06a}, and therefore are omitted.

In case $K=K_{f_1}$ is the score function associated with $f_1\in
\mathcal{F}_a$, and provided that Assumption~\ref{assuA} holds (in order for
the central sequence $\Deltab_{ \varthetab; f_{1}}$ of
Pro\-position~\ref{LAN} to make sense), $\ubDelta_{ \varthetab; K_{f_1},
f_{1}}^{\III}$ and $\ubDelta_{ \varthetab; K_{f_1}, f_{1}}^{\IV}$,
under $\mathrm{P}^{(n)}_{\varthetab; f_{1}}$ clearly coincide with
$\Deltab_{ \varthetab; f_{1}}^{\III}$ and $\Deltab_{ \varthetab;
f_{1}}^{\IV}$. Therefore, ${\utDelta}{}_{ \varthetab; K_{f_1}}^{\III
}$ and ${\utDelta}{}_{ \varthetab; K_{f_1}}^{\IV}$
constitute rank-based, hence distribution-free, versions of
those central sequence components. Exploiting this, we now construct
signed-rank tests for the two problems we are interested in.\vspace*{-2pt}

%s6.2 ###
\subsection{Optimal rank-based tests for eigenvectors}
\label{gsjdlr}
Proposition~\ref{Hajek} provides the theoretical tools for
constructing rank-based tests for ${\mathcal H}_{0}^{\betab}$ and
computing their local powers. Letting again $\varthetab_0:=({\bolds
\theta}\pr, \sigma^2, (\dvecrond\Lamb_\Vb)\pr,(\vecop\betab
_0)\pr)\pr$, with $\betab_0=(\betab^0,\betab_2,\ldots, \betab_k
)$, define the rank-based analog of~(\ref{Qparam})
%
%e6.2 ###
%
\begin{eqnarray}\label{Qrank}\quad
{\utQ}{}_{\varthetab_0 ;K}^{(n)}
:\!&=&
{\utDelta}{}_{\varthetab_0 ; K}^{\IV\prime}
[(\Gamb^{\IV}_{\varthetab_0 ; K})^{-}- \mathbf{P}_{k}^{\betab
_{0}^{}} ( (\mathbf{P}_{k}^{\betab_{0}^{}})\pr\Gamb^{\IV
}_{\varthetab_0 ; K} \mathbf{P}_{k}^{\betab_{0}^{}} )^{-}(\mathbf
{P}{}^{\betab_{0}^{}}_{k})\pr
]
{\utDelta}{}_{\varthetab_0 ; K}^{\IV}
\nonumber\\[-9pt]\\[-9pt]
&=&\frac{nk(k+2)}{ \mathcal{J}_k(f_1)} \sum_{j=2}^{k} \bigl(
\betab_j\pr{\utSb}{}^{(n)}_{\varthetab_0 ;K} \betab^0\bigr)^2,\nonumber\vspace*{-2pt}
\end{eqnarray}
where
$
{\utSb}{}^{(n)}_{\varthetab;K}:=\frac{1}{n} \sum_{i=1}^{n} K
(\frac{R^{(n)}_{i}({\bolds\theta}, {\Vb})}{n+1} ) \Ub
_{i}({\bolds\theta}, {\Vb}) \Ub\pr_{i}({\bolds\theta}, {\Vb}
).
$\vadjust{\goodbreak}

In order\vspace*{-4pt} to turn ${\utQ}{}_{\varthetab_0 ;K}^{(n)}$ into a genuine test
statistic, as in the parametric case, we still have to replace
$\varthetab_0$ with some adequate estimator $\hat{\varthetab}$
satisfying, under as large as possible a class of densities, Assumption
\ref{assuB} for $\mathcal{H}_0^\betab$. In particular, root-$n$ consistency
should hold without any moment assumptions. Denote by $\hat{{\bolds
\theta}}_{\mathrm{HR}}$ the \citet{HR02}
affine-equivariant median, and by $\hats{\Vb}_{\mathrm{Tyler}}$ the
shape estimator of \citet{T87}, normalized so that it has determinant
one: both are root-$n$ consistent under any radial density $g_1$.
Factorize $\hats{\Vb}_{\mathrm{Tyler}}$ into $\hat\betab_{\mathrm
{Tyler}}\hat\Lamb_{\mathrm{Tyler}}\hat\betab_{\mathrm{Tyler}}\pr
$. The estimator we are proposing (among many possible ones) is
$ \hat\varthetab=(\hat{{\bolds\theta}}_{\mathrm{HR}}\pr, {\sigma
}^{2}, (\dvecrond\hat{\Lamb}_{\mathrm{Tyler}})\pr,\break (\vecop
\tilde{\betab}_{0})\pr)\pr$,
where the constrained estimator $\tilde{\betab}_{0} := (\betab
^{0}, \tilde{\betab}_{2}, \ldots, \tilde{\betab}_{k})$ is
constructed from $\hat\betab_{\mathrm{Tyler}}$ via the same
Gram--Schmidt procedure as was applied in Section~\ref{Gaussbetab} to
the eigenvectors $\hat\betab_{\Vb}$ of $\hats{\Vb}:=\Sb^{(n)}/|\Sb
^{(n)}|^{1/k}$; note that $\sigma^2$ does not even appear in ${\utQ}{}
_{\varthetab_0 ;K}^{(n)}$, hence needs not be
estimated.

In view of~(\ref{Qrank}),
%
%e6.3 ###
%
\begin{eqnarray} \label{Qchapeau}
{\utQ}{}_{K}^{(n)}&=&
{\utQ}{}^{(n)}_{\hat\varthetab;K}\nonumber\\
&=& \frac{nk(k+2)}{\mathcal{J}_k(K)} \betab^{0\prime} {\utSb
}{}^{(n)}_{\hat\varthetab;K} (\mathbf{I}_k-\betab^{0}\betab^{0\prime
}) {\utSb}{}_{\hat\varthetab;K}^{(n)}\betab_{}^{0 }
\\
& =&
\frac{nk(k+2)}{\mathcal{J}_k(K)}
\bigl\| [\betab_{}^{0\prime} \otimes(\mathbf{I}_k-\betab
_{}^{0}\betab^{0\prime}) ] \bigl(\vecop\utSb{}^{(n)}_{\hat
\varthetab;K} \bigr)
\bigr\|^2,\nonumber
\end{eqnarray}
where the\vspace*{-2pt} ranks and signs in $\utSb{}_{\hat\varthetab;K}$ are
computed at $\hat\varthetab$, that is, $R^{(n)}_{i}:=R^{(n)}_{i}(\hat
{{\bolds\theta}}_{\mathrm{HR}}, \tilde{\betab}_{0} \hat\Lamb
_{\mathrm{Tyler}}\tilde{\betab}_{0} \pr)$ and ${\Ub}_{i}={\Ub
}_{i}(\hat{{\bolds\theta}}_{\mathrm{HR}}, \tilde{\betab}_{0} \hat
\Lamb_{\mathrm{Tyler}}\tilde{\betab}_{0} \pr)$.

Let us show that
substituting $\hat{\varthetab}$ for $\varthetab_{0}$ in (\ref
{Qrank}) has no asymptotic impact on ${\utQ}{}_{\varthetab_0
;K}^{(n)}$---that is, ${\utQ}{}_{\varthetab_0 ;K}^{(n)}-{\utQ
}{}_{K}^{(n)}=o_\mathrm{P}(1)$ as $\ny$ under\vspace*{-4pt} $\mathrm
{P}^{(n)}_{\varthetab
_0; g_1}$, with $g_{1}\in\mathcal{F}_a$, $\varthetab_0 \in{\mathcal
H}_{0;1}^{\betab\prime}$. The proof, as usual, relies on an
asymptotic linearity property which, in turn, requires ULAN. The ULAN
property of Proposition~\ref{LAN}, which was motivated by optimality
issues in tests involving $\betab$ and $\Lamb_\Vb$, here cannot help
us, as it does not hold under Assumption~\ref{assuApr1}. Another ULAN
property, however, where Assumption~\ref{assuA} is not required, has been
obtained by \citet{HP06a} for another
parametrization---based on $({\bolds\theta},\sigma^2,\Vb)$---of the
same families of distributions, and perfectly fits our needs here.

Defining $\mathbf{J}_{k}:= (\vecop\mathbf{I}_{k})(\vecop\mathbf
{I}_{k})\pr
$ and $\mathbf{J}_{k}^{\perp}:= \mathbf{I}_{k^{2}}- \frac{1}{k} \mathbf
{J}_{k}$, it follows from Proposition A.1 in \citet{HOP06} and Lemma 4.4 in \citet{K87} that, for
any locally asymptotically discrete [Assumption \hyperlink{B3}{(B3)}] and root-$n$
consistent [Assumption \hyperlink{B2}{(B2)}] sequence $( \hat{\bolds\theta}^{(n)},
\hats{\Vb}^{(n)})$ of estimators of location and shape, one
has
%
%e6.4 ###
%
\begin{eqnarray}
\label{HPres}
&&
\mathbf{J}_{k}^{\perp} \sqrt{n} \vecop({\utSb}{}_{\hat\varthetab
;K}- {\utSb}{}_{\varthetab;K})
\nonumber\\
&&\quad{} +
\frac{{\mathcal J}_k(K,g_{1})}{4k(k+2)}
\biggl[\mathbf{I}_{k^{2}}+\mathbf{K}_{k}-\frac{2}{k}\mathbf{J}_{k} \biggr]({\Vb
}^{-1/2})^{\otimes2} n^{1/2}\vecop\bigl(\hats{\Vb}^{(n)}- \Vb\bigr)\\
&&\qquad= o_\mathrm{P}(1)\nonumber
\end{eqnarray}
as $\ny$ under $\mathrm{P}^{(n)}_{\varthetab_0 ;g_{1}}$,
with $\varthetab_0 \in\Thetab$ and $g_{1} \in\mathcal{F}_a$. This
result readily applies to any adequately discretized version of $(\hat
{{\bolds\theta}}_{\mathrm{HR}}, \tilde{\betab}_{0} \hat\Lamb
_{\mathrm{Tyler}}\tilde{\betab}_{0} \pr)$ at $\varthetab_0 \in
{\mathcal H}_{0;1}^{\betab\prime}$. It is well known, however, that
discretization, although necessary for asymptotic statements, is not
required in practice [see pages 125 or 188 of \citet{CY00} for
a discussion on this point]; we therefore do not emphasize
discretization any further in the notation, and henceforth assume that
$\hat\varthetab$, whenever needed, has been adequately discretized.

Using~(\ref{HPres}) and the fact that
$ [\betab^{0\prime} \otimes(\mathbf{I}_k-\betab_{}^{0}\betab
^{0\prime}) ]\mathbf{J}_{k}=\mathbf{0}$,
we obtain,\break under $\mathrm{P}^{(n)}_{\varthetab_0 ;g_{1}}$
with ${{\varthetab}_0} \in{\mathcal H}_{0;1}^{\betab\prime}$ and
$g_{1} \in\mathcal{F}_a$, since $\mathbf{K}_{k}(\vecop\mathbf
{A})=\vecop(\mathbf{A}\pr)$ for any $k\times k$ matrix $\mathbf{A}$
and since $\betab^0$ under ${{\varthetab}_0} \in{\mathcal
H}_{0;1}^{\betab\prime}$ is an eigenvector of ${\Vb}^{-1/2}\tilde
{\betab}_{0} \hat\Lamb_{\mathrm{Tyler}}\tilde{\betab}_{0} \pr
{\Vb}^{-1/2}$,
\begin{eqnarray*}%\label{nullspace}
&&\sqrt{n} [\betab^{0\prime} \otimes(\mathbf{I}_k-\betab
_{}^{0}\betab^{0\prime}) ] \vecop\bigl(\utSb{}^{(n)}_{\hat
\varthetab;K}-\utSb{}^{(n)}_{\varthetab;K}\bigr)
\\
&&\qquad
= -\frac{{\mathcal J}_k(K,g_{1})}{2k(k+2)} n^{1/2} [\betab
^{0\prime} \otimes(\mathbf{I}_k-\betab_{}^{0}\betab^{0\prime}) ]
({\Vb}^{-1/2})^{\otimes2} \vecop(\tilde{\betab}_{0} \hat\Lamb
_{\mathrm{Tyler}}\tilde{\betab}_{0} \pr- \Vb)\\
&&\qquad\quad{} +o_\mathrm{P}(1)
\\
&&\qquad
= -\frac{{\mathcal J}_k(K,g_{1})}{2k(k+2)}n^{1/2}
\vecop\bigl((\mathbf{I}_k-\betab_{}^{0}\betab^{0\prime}) [{\Vb
}^{-1/2}\tilde{\betab}_{0} \hat\Lamb_{\mathrm{Tyler}}\tilde
{\betab}_{0} \pr{\Vb}^{-1/2}-\mathbf{I}_k]\betab^0\bigr)\\
&&\qquad\quad{} +o_\mathrm{P}(1)
\\
&&\qquad
=o_\mathrm{P}(1),
\end{eqnarray*}
as $\ny$. In\vspace*{-6pt} view of~(\ref{Qchapeau}), we therefore conclude that $
{\utQ}{}_{K}^{(n)}-{\utQ}{}_{\varthetab_0 ; K}^{(n)} = o_\mathrm{P}(1)$ as
$\ny$, still under under ${{\varthetab}_0} \in{\mathcal
H}_{0;1}^{\betab\prime}$, as was to be shown.

The following result summarizes the results of this section.
\begin{Prop} \label{ranktestbeta}Let Assumption~\ref{assuS} hold for the score
function $K$. Then:

\begin{longlist}
\item
${\utQ}{}_{K}^{(n)}$ is asymptotically\vspace*{-2pt} chi-square with $(k-1)$
degrees of freedom under $\bigcup_{\varthetab\in{\mathcal
H}_{0;1}^{\betab\prime}} \bigcup_{ g_{1}\in\mathcal{F}_a }
\{ \mathrm{P}^{(n)}_{\varthetab;g_{1}}\}$, and asymptotically noncentral
chi-square, still with $(k-1)$ degrees of
freedom, and noncentrality parameter
\[
\frac{{\mathcal J}^{2}_k(K, g_{1} ) }{4k(k+2)\mathcal{J}_k(K)}
r_{\varthetab; \taub}^{\betab}
\]
under $\mathrm{P}^{(n)}_{\varthetab+n^{-1/2}\taub^{(n)};g_{1}}$, for
$\varthetab\in
{\mathcal H}_{0}^{\betab}$ and $g_{1}\in\mathcal{F}_a$, with $\taub
^{(n)}$ as in Proposition~\ref{pseudogausstestbeta} and $
r_{\varthetab; \taub}^{\betab} $ defined in~(\ref{rbeta});

\item
the sequence\vspace*{-1pt} of tests ${\utphi}{}^{(n)}_{\betab; K}$ rejecting
the null when $Q_{K}^{(n)}$ exceeds the $\alpha$ upper-quantile of the
chi-square distribution with $(k-1)$ degrees of freedom has asymptotic size
$\alpha$ under $\bigcup_{\varthetab\in{\mathcal H}_{0}^{\betab
\prime}}\bigcup_{ g_{1}\in\mathcal{F}_a} \{ \mathrm{P}^{(n)}_{\varthetab
;g_{1}}\}$;

\item
for scores\vspace*{-2pt} $K=K_{f_1}$, with $f_{1}\in\mathcal{F}_a$, ${\utphi
}{}^{(n)}_{\betab; K}$ is locally asymptotically most stringent, at
asymptotic level $\alpha$, for $\bigcup_{\varthetab\in{\mathcal
H}_{0}^{\betab\prime}}\bigcup_{ g_{1}\in\mathcal{F}_a} \{
\mathrm{P}^{(n)}_{\varthetab;g_{1}}\}$ against alternatives of the form
$\bigcup_{\varthetab\notin{\mathcal H}_{0}^{\betab}} \{ \mathrm
{P}^{(n)}_{\varthetab;f_{1}}\}$.\vspace*{-1pt}
\end{longlist}
\end{Prop}

Being measurable with respect to\vspace*{-2pt}
%Because of the use of
signed-ranks, ${\utQ}{}_{K}^{(n)}$ is
asymptotically invariant under continuous monotone
radial transformations, in the sense that it is asymptotically
equivalent (in probability) to a random variable that is strictly
invariant under such transformations. Furthermore, it is easy to show
that it enjoys the same $\mathcal{G}_{\mathrm{rot},\circ}$-invariance
features as the parametric, Gaussian, or pseudo-Gaussian test
statistics.%\looseness=1

%s6.3 ###
\subsection{Optimal rank-based tests for eigenvalues}\label{rankeigvalues}
Finally, still from the results of Proposition~\ref{Hajek}, we
construct signed-rank tests for the null hypothesis ${\mathcal
H}_{0}^{\Lamb}$. A~rank-based counterpart of~(\ref{Tparam1}) and~(\ref{Tparam})
[at $\varthetab_0=({\bolds\theta}\pr, \sigma^2, (\dvecrond\Lamb
_0)\pr,\break
(\vecop\betab)\pr)\pr\in\mathcal{H}^\Lamb_0$] is, writing
$\Vb_0$ for $\betab\Lamb_0\betab\pr$,
%
%e6.5 ###
%
\begin{eqnarray}\label{Trank}\qquad
{\utT}{}_{\varthetab_0; K}^{(n)}
& = & (\operatorname{grad}\pr h(\dvecrond\Lamb_0)
(\Gamb^{\III}_{{\bolds\vartheta}_{0}; K})^{-1}\operatorname{grad}
h(\dvecrond\Lamb_0) )^{-1/2} \nonumber\\
& &{}
\times\operatorname{grad}\pr h(\dvecrond\Lamb_{0}) (\Gamb
^{\III}_{{\bolds\vartheta}_{0}; K})^{-1} {\utDelta}{}^{\III}
_{{\bolds\vartheta}_{0} ; K} \\
& = &
\biggl(\frac{nk(k+2)}{\mathcal{J}_k(K)} \biggr)^{1/2} ( a_{p,q}(\Lamb
_0 ))^{-1/2}
\mathbf{c}_{p,q}\pr\dvec\bigl( \Lamb_0^{1/2}\betab\pr{\utSb
}{}^{(n)}_{\varthetab_{0},K} \betab\Lamb_0^{1/2}\bigr).
\nonumber
\end{eqnarray}
Here again, we have to estimate $\varthetab_0$. Note that, unlike the
quantity\break
$\operatorname{grad}\pr h(\dvecrond\Lamb_{0} ) (\Gamb^{\III
}_{{\bolds\vartheta}_{0}; \phi_1})^{-1} \Deltab^{\III} _{{\bolds
\vartheta}_{0}; \phi_1}$ appearing in the Gaussian or
pseudo-Gaus\-sian~cases,
$\operatorname{grad}\pr h(\dvecrond\Lamb_{0} ) (\Gamb^{\III
}_{{\bolds\vartheta}_{0}; K})^{-1} {\utDelta}{}^{\III} _{{\bolds
\vartheta}_{0}; K}$ does depend\vspace*{-2pt} on $\Lamb_{0}$ [see the comments
below~(\ref{Tgauss})]. Consequently, we have to carefully select an
estimator $\hat\varthetab$ that has no influence on the asymptotic
behavior of ${\utT}{}_{\varthetab_0; K}^{(n)}$ under~${\mathcal
H}_{0;q}^{\Lamb\prime\prime}$.

To this end, consider Tyler's estimator of shape $\hats{\mathbf
{V}}_{\mathrm
{Tyler}}(=:\hat{\betab}_{\mathrm{Tyler}}\hat{\Lamb}_{\mathrm
{Tyler}}\hat{\betab}_{\mathrm{Tyler}}\pr$, with obvious notation)
and define
\[
\dvec(\tilde{\Lamb}_{\mathrm{Tyler}}):= \bigl(\mathbf{I}_k- \mathbf
{c}_{p,q}(\mathbf{c}_{p,q}\pr\mathbf{c}_{p,q})^{-1}\mathbf{c}_{p,q}\pr\bigr)
(\dvec\hat{\Lamb}_{\mathrm{Tyler}}).
\]
Then the estimator of shape
$
\tilde{\Lamb}_{\Vb}:= \tilde{\Lamb}_{\mathrm{Tyler}}/ |\tilde
{\Lamb}_{\mathrm{Tyler}}|^{1/k}
$
is clearly constrained: $\mathbf{c}_{p,q}\pr(\dvec\tilde{\Lamb}_{\Vb
})=0$ and $|\tilde{\Lamb}_{\Vb}|=1$. The resulting preliminary
estimator $\hat{\varthetab}$ is
%
%e6.6 ###
%
\begin{equation} \label{initialresti}
\hat{\varthetab}:=
(\hat{{\bolds\theta}}_{\mathrm{HR}}\pr, {\sigma}^{2},
(\dvecrond\tilde{\Lamb}_{\Vb})\pr, (\vecop\hat{\betab
}_{\mathrm{Tyler}})\pr)\pr,
\end{equation}
where $\hat{{\bolds\theta}}_{\mathrm{HR}}$ still denotes the
\citet{HR02} affine-equivariant median. The test
statistic we propose is then
%
%e6.7 ###
%
\begin{eqnarray}\label{717}
{\utT}{}_{K}^{(n)}
:\!& = &
{\utT}{}_{\hat\varthetab; K}^{(n)}
\nonumber\\
&=& (
\operatorname{grad}\pr h(\dvecrond\tilde{\Lamb}_{\Vb} )
(\Gamb^{\III}_{\hat{\bolds\vartheta}; K})^{-1}
\operatorname{grad} h( \dvecrond\tilde{\Lamb}_{\Vb} )
)^{-1/2}
\nonumber\\
& &{}
\times
\operatorname{grad}\pr h(\dvecrond\tilde{\Lamb}_{\Vb})
(\Gamb^{\III}_{\hat{\bolds\vartheta}; K})^{-1}
{\utDelta}{}^{\III} _{\hat{\bolds\vartheta}; K}
\\
& = &
\biggl(\frac{nk(k+2)}{\mathcal{J}_k(K)} \biggr)^{1/2}
( a_{p,q} (\tilde{\Lamb}_\Vb))^{-1/2}
\mathbf{c}_{p,q} \pr\nonumber\\
&&{}\times\dvec\bigl(\tilde{\Lamb}_{\Vb}^{1/2}\hat{\betab
}_{\mathrm{Tyler}} \pr{\utSb}{}_{\hat{\varthetab};K}^{(n)}\hat
{\betab}_{\mathrm{Tyler}} \tilde{\Lamb}_{\Vb}^{1/2}\bigr), \nonumber
\end{eqnarray}
where ${\utSb}{}_{\hat{\varthetab};K}^{(n)}:=\frac{1}{n} \sum_{i=1}^{n}
K (\frac{\hat R^{(n)}_{i}}{n+1} ) {\hat\Ub}_{i}{\hat\Ub
}_{i}\pr$,
with $ \hat R^{(n)}_{i} := R^{(n)}_{i}(\hat{{\bolds\theta}}_{\mathrm
{HR}},
\hat{\betab}_{\mathrm{Tyler}} \tilde{\Lamb}_{\Vb}
\hat{\betab}_{\mathrm{Tyler}} ^{ \prime}) $ and $ \hat
\Ub_i:=\Ub_i(\hat{{\bolds\theta}}_{\mathrm{HR}} ,
\hat{\betab}_{\mathrm{Tyler}} \tilde{\Lamb}_{\Vb}
\hat{\betab}_{\mathrm{Tyler}} ^{ \prime})$. The following lemma shows
that the substitution of $\hat\varthetab$ for ${\varthetab}$
in~(\ref{initialresti}) has no asymptotic effect on ${\utT
}{}_{\varthetab
; K}^{(n)}$ (see the \hyperref[app]{Appendix} for a proof).
\begin{Lem} \label{alihoprimeprime}
Fix $\varthetab\in{\mathcal H}_{0;q}^{\Lamb\prime\prime}$ and
$g_{1} \in\mathcal{F}_{a}$, and let $\hat{\varthetab}$ be the
estimator in~(\ref{initialresti}). Then ${\utT}{}_{K}^{(n)}-{\utT
}{}_{\varthetab; K}^{(n)}$ is $o_\mathrm{P}(1)$ as $\ny$, under $\mathrm
{P}^{(n)}_{\varthetab; g_{1}}$.
\end{Lem}

The following result summarizes the results of this section.
\begin{Prop} \label{ranktestlambda}
Let Assumption~\ref{assuS} hold for the score function $K$. Then:

\begin{longlist}
\item
${\utT}{}^{(n)}_{K}$ is asymptotically standard normal under
$\bigcup_{\varthetab\in{\mathcal H}_{0;q}^{\Lamb\prime\prime}}
\bigcup_{g_{1} \in{\mathcal F}_a}
\{ \mathrm{P}^{(n)}_{\varthetab;g_{1}}\}$, and
asymptotically normal with mean
\[
\frac{\mathcal{J}_k(K,g_{1})}{\sqrt{4k(k+2) a_{p,q}(\Lamb_\Vb
)\mathcal{J}_k(K)}} r_{\varthetab; \taub}^{\Lamb_{\Vb}}
\]
and variance 1
under $\mathrm{P}^{(n)}_{\varthetab+n^{-1/2}\taub^{(n)};g_{1}}$, with
$\varthetab\in{\mathcal H}_{0}^{\Lamb_{\Vb}}$, $g_{1} \in{\mathcal
F}_a$, $ \taub^{(n)}$ as in Proposition~\ref{pseudogausstestlambda},
and $r_{\varthetab; \taub}^{\Lamb_{\Vb}}$ defined in~(\ref{rLamb});
\item
the sequence of tests ${\utphi}{}^{(n)}_{\Lamb; K}$ rejecting
the null whenever ${\utT}{}^{(n)}_{K}$ is less than the standard
normal $\alpha$-quantile $z_\alpha$ has asymptotic size
$\alpha$ under\break $\bigcup_{\varthetab\in
{\mathcal H}_{0;q}^{\Lamb\prime\prime}}\bigcup_{g_{1}\in\mathcal
{F}_a } \{ \mathrm{P}^{(n)}_{\varthetab;g_{1}}\}$;
\item
for scores\vspace*{-2pt} $K=K_{f_1}$ with $f_{1} \in{\mathcal F}_a$, the
sequence of tests ${\utphi}{}^{(n)}_{\Lamb; K}$ is locally and
asymptotically most powerful,
still at asymptotic level $\alpha$, for $\bigcup_{\varthetab\in
{\mathcal H}_{0;q}^{\Lamb\prime\prime}}\bigcup_{g_{1}\in\mathcal
{F}_a } \{ \mathrm{P}^{(n)}_{\varthetab;g_{1}}\}$ against alternatives of
the form
$\bigcup_{\varthetab\notin{\mathcal H}_{0}^{\Lamb}} \{ \mathrm
{P}^{(n)}_{\varthetab;f_1 }\}$.
\end{longlist}
\end{Prop}
%
%%%%%%%%%%%%%%%%%%%%%%%%%%%%%%%%%%%%%%%%%%%%%%%%%%%%%%%%%%%%%%%%%%%%%%%%%%%%%%%%%%%%%%%%%%%%%%%%%%%%%%%%%%%%%%%%%%%%%%%%%%%%%%%%%%%%%%%%%%%%%%%%

%s7 ###
\section{Asymptotic relative efficiencies}\label{secare}

The asymptotic relative efficiencies (AREs) of the rank-based tests of
Section~\ref{ranktests} %${\utphi}{}^{(n)}_{\betab; K}$
with respect to their Gaussian and pseudo-Gaussian competitors of
Sections~\ref{gausscase} are readily obtained as
ratios of noncentrality parameters under
local alternatives (squared ratios of standardized asymptotic shifts
for the one-sided problems on eigenvalues).
Denoting by $\mathrm{ARE}^{\varthetab,\taub}_{k,g_1}(\phi^{(n)}_1/\phi
^{(n)}_2)$ the ARE, under local\vspace*{-2pt} alternatives of the form $\mathrm
{P}^{(n)}_{\varthetab+n^{-1/2}\taub;g_1}$, of a sequence of tests
$\phi^{(n)}_1$ with respect to the sequence~$\phi^{(n)}_2$, we thus
have the following result.
\begin{Prop} \label{6are}
Let Assumptions~\ref{assuS} and~\ref{assuB} hold for the score function $K$ and (with
the appropriate null hypotheses and densities) for the estimators $\hat
\varthetab$ described in the previous sections. Then, for any $g_1
\in\mathcal{F}_a^{4}$,
\[
\mathrm{ARE}^{\varthetab,\taub}_{k,g_1}\bigl({\utphi}{}^{(n)}_{\betab;
K}/\phi
^{(n)}_{\betab; \mathcal{N}*}\bigr)=\mathrm{ARE}^{\varthetab,\taub
}_{k,g_1}\bigl({\utphi}{}^{(n)}_{\Lamb; K}/\phi
^{(n)}_{\Lamb; \mathcal{N}*}\bigr) :=
\frac{(1+\kappa_{k}(g_1))\mathcal{J}^2_k(K,g_1)}{k(k+2)\mathcal{J}_k(K)}.
\]
\end{Prop}

Table~\ref{AREtable} provides numerical values of these AREs for
various values of the space
%
%t1 ###
%
\begin{table}
\caption{AREs of the van\vspace*{-3pt} der Waerden (vdW), Wilcoxon (W) and Spearman
(SP) rank-based tests~${\utphi}{}^{(n)}_{\betab; K}$ and $
{\utphi}{}^{(n)}_{\Lamb; K}$ with respect to their pseudo-Gaussian
counterparts, under $k$-dimensional Student (with $5$, $8$ and $12$
degrees of freedom), Gaussian, and power-exponential densities (with
parameter $\eta=2,3,5$), for $k=2$, $3$, $4$, $6$, $10$, and $k
\rightarrow\infty$} \label{AREtable}
\begin{tabular*}{\tablewidth}{@{\extracolsep{\fill}}lcccccccc@{}}
\hline
& & \multicolumn{7}{c@{}}{\textbf{Underlying density}}\\[-4pt]
& & \multicolumn{7}{c@{}}{\hrulefill}\\
$\bolds K$& $\bolds k$ & $\bolds{t_{5}}$ & $\bolds{t_{8}}$ &
$\bolds{t_{12}}$ & $\bolds{\mathcal{N}}$ & $\bolds{e_2}$ &
$\bolds{e_3}$ & $\bolds{e_5}$ \\
\hline
vdW
& \phantom{0}2 & 2.204 & 1.215 & 1.078 & 1.000 & 1.129 & 1.308 & 1.637 \\
& \phantom{0}3 & 2.270 & 1.233 & 1.086 & 1.000 & 1.108 & 1.259 & 1.536 \\
& \phantom{0}4 & 2.326 & 1.249 & 1.093 & 1.000 & 1.093 & 1.223 & 1.462 \\
& \phantom{0}6 & 2.413 & 1.275 & 1.106 & 1.000 & 1.072 & 1.174 & 1.363 \\
& 10 & 2.531 & 1.312 & 1.126 & 1.000 & 1.050 & 1.121 & 1.254 \\
& $\infty$ & 3.000 & 1.500 & 1.250 & 1.000 & 1.000 & 1.000 & 1.000 \\
[3pt]
W
& \phantom{0}2 & 2.258 & 1.174 & 1.001 & 0.844 & 0.789 & 0.804 & 0.842 \\
& \phantom{0}3 & 2.386 & 1.246 & 1.068 & 0.913 & 0.897 & 0.933 & 1.001 \\
& \phantom{0}4 & 2.432 & 1.273 & 1.094 & 0.945 & 0.955 & 1.006 & 1.095 \\
& \phantom{0}6 & 2.451 & 1.283 & 1.105 & 0.969 & 1.008 & 1.075 & 1.188 \\
& 10 & 2.426 & 1.264 & 1.088 & 0.970 & 1.032 & 1.106 & 1.233 \\
& $\infty$ & 2.250 & 1.125 & 0.938 & 0.750 & 0.750 & 0.750 & 0.750 \\
[3pt]
SP
& \phantom{0}2 & 2.301 & 1.230 & 1.067 & 0.934 & 0.965 & 1.042 & 1.168 \\
& \phantom{0}3 & 2.277 & 1.225 & 1.070 & 0.957 & 1.033 & 1.141 & 1.317 \\
& \phantom{0}4 & 2.225 & 1.200 & 1.051 & 0.956 & 1.057 & 1.179 & 1.383 \\
& \phantom{0}6 & 2.128 & 1.146 & 1.007 & 0.936 & 1.057 & 1.189 & 1.414 \\
& 10 & 2.001 & 1.068 & 0.936 & 0.891 & 1.017 & 1.144 & 1.365 \\
& $\infty$ & 1.667 & 0.833 & 0.694 & 0.556 & 0.556 & 0.556 & 0.556 \\
\hline
\end{tabular*}
\end{table}
dimension $k$ and selected radial densities $g_1$ (Student, Gaussian
and power-expo\-nential), and for the van der
Waerden tests ${\utphi}{}^{(n)}_{ \betab;{\mathrm{vdW}}}$ and ${\utphi
}{}^{(n)}_{ \Lamb;{\mathrm{vdW}}}$, the\vspace*{-2pt} Wilcoxon tests
${\utphi}{}^{(n)}_{ \betab; K_1}$ and ${\utphi}{}^{(n)}_{ \Lamb
; K_1}$, and the Spearman tests ${\utphi}{}^{(n)}_{ \betab; K_2}$
and ${\utphi}{}^{(n)}_{ \Lamb; K_2}$ (the score functions $K_a$,
$a>0$ were defined in Section~\ref{scores}). These values
coincide with the ``AREs for shape'' obtained in \citet{HP06a}, which implies
[\citet{P06}] that the AREs of
van der Waerden tests with respect to their pseudo-Gaussian
counterparts are uniformly larger than or equal to one (an extension
of the classical Chernoff--Savage property):
\[
\inf_{g_1} \operatorname{ARE}^{\varthetab,\taub}_{k,g_1}\bigl({\utphi
}{}^{(n)}_{\betab; \mathrm{vdW}}/\phi^{(n)}_{\betab; \mathcal{N}*}\bigr)=\inf_{g_1}
\operatorname{ARE}^{\varthetab,\taub}_{k,g_1}\bigl({\utphi}{}^{(n)}_{\Lamb;
\mathrm{vdW}}/\phi^{(n)}_{\Lamb; \mathcal{N}*}\bigr)=1.
\]

%s8 ###
\section{Simulations} \label{secsimu}

In this section, we investigate via simulations the finite-sample
performances of the following tests:

(i) the Anderson test $\phi_{ \betab;\mathrm{Anderson}}^{(n)}$,
the optimal Gaussian test $\phi_{ \betab; {\mathcal N}}^{(n)}$,
the pseudo-Gaussian test $\phi_{ \betab; {\mathcal N}{*}}^{(n)}$,
the robust test $\phi^{(n)}_{\betab;\mathrm{Tyler}}$ based on
$Q_{\mathrm{Tyler}}^{(n)}$, and various rank-based tests
$\phi_{ \betab; K}^{(n)}$ (with van der Waerden, Wilcoxon, Spearman
and sign scores, but also with scores achieving optimality at $t_1$,
$t_3$ and $t_5$ densities), all for the null hypothesis ${\mathcal
H}_{0}^{\betab}$ on eigenvectors;

(ii) the optimal Anderson test $\phi_{ \Lamb;\mathrm
{Anderson}}^{(n)}=\phi_{ \Lamb; {\mathcal N}}^{(n)}$,
the pseudo-Gaussian test $\phi_{ \Lamb; {\mathcal N}{*}}^{(n)}
%$,
%the robust test $
=\phi^{(n)}_{ \Lamb; \mathrm{Davis}}$ based on $T_{\mathrm
{Davis}}^{(n)}$, and various rank-based tests
$\phi_{ \Lamb; K}^{(n)}$ (still with van der Waerden, Wilcoxon,
Spearman, sign, $t_1$, $t_3$ and $t_5$ scores), for the null
hypothesis ${\mathcal H}_{0}^{ \Lamb} $ on eigenvalues.

Simulations were conducted as follows. We generated
$N = 2500$
mutually independent samples of i.i.d. trivariate ($k=3$) random vectors
$
{\bolds\varepsilon}_{\ell;j}$,
$\ell=1, 2, 3, 4, j=1,\ldots, n=100,
$
with spherical Gaussian (${\bolds\varepsilon}_{1;j}$), $t_{5}$
(${\bolds\varepsilon}_{2;j}$), $t_3$ (${\bolds\varepsilon}_{3;j}$)
and $t_{1}$ (${\bolds\varepsilon}_{4;j}$) densities, respectively. Letting
\[
\Lamb:= \pmatrix{ 10 & 0 & 0
\cr0 & 4 & 0 \cr0 & 0 & 1},\qquad
\mathbf{B}_{\xi}:= \pmatrix{
\cos(\pi\xi/12) & - \sin(\pi\xi/12) & 0
\cr\sin(\pi\xi/12) & \cos(\pi\xi/12) & 0 \cr0& 0& 1}
\]
and
\[
\mathbf{L}_{\xi}:= \pmatrix{ 3 \xi& 0 & 0
\cr0 & 0 & 0 \cr0& 0& 0},
\]
each ${\bolds\varepsilon}_{\ell;j}$ was successively transformed into
%
%e8.1 ###
%
\begin{equation}\label{sample1}\quad
\Xb_{\ell;j ;\xi}= \mathbf{B}_{\xi} \Lamb^{1/2}{\bolds\varepsilon}_{\ell
;j},\qquad \ell=1,
2, 3, 4, j=1,\ldots,n, \xi=0,\ldots,3 ,
\end{equation}
and
%
%e8.2 ###
%
\begin{equation}\label{sample2}\quad\qquad
\Yb_{\ell;j;\xi}= (\Lamb+\mathbf{L}_{\xi})^{1/2}{\bolds\varepsilon
}_{\ell;j},\qquad
\ell=1, 2, 3, 4, j=1,\ldots,n, \xi=0,\ldots, 3 .
\end{equation}
The value $ \xi=0$ corresponds to the null hypothesis ${\mathcal
H}_{0}^{\betab}\dvtx\betab_{1}=(1, 0, 0)\pr$ for the $\Xb_{\ell;j;\xi
}$'s and the null hypothesis ${\mathcal H}_{0}^{\Lamb}\dvtx{\sum
_{j=q+1}^{k} \lambda_{j; \Vb}}/{\sum_{j=1}^{k} \lambda_{j; \Vb
}}=1/3$ (with $q=1$ and $k=3$) for the $\Yb_{\ell;j;\xi}$'s; $\xi
=1,2,3$ characterizes increasingly distant alternatives. We then
performed the tests listed under (i) and (ii) above in $N=2500$
independent replications of such samples.
Rejection frequencies are reported in Table~\ref{simuresu2}
for ${\mathcal H}_{0}^{\betab}$ and in Table~\ref{simuresu3} for
${\mathcal H}_{0}^{\Lamb}$.

%t2 ###
%
\begin{table}
\caption{Rejection frequencies (out of $N=2500$
replications), under the null ${\mathcal H}_{0}^{\betab}$ and
increasingly distant alternatives (see Section \protect\ref{secsimu} for
details), of the Anderson test $\phi^{(n)}_{\betab; \mathrm
{Anderson}}$, the
Tyler test $\phi^{(n)}_{\betab; \mathrm{Tyler}}$, the parametric Gaussian
test $\phi^{(n)}_{\betab; \mathcal{N}}$, its pseudo-Gaussian
version $\phi^{(n)}_{\betab; \mathcal{N}*}$, and the signed-rank tests
with van der Waerden, $t_\nu$ ($\nu=1$, $3$, $5$), sign, Wilcoxon,
and Spearman scores, ${\utphi}{}^{(n)}_{\betab; \mathrm{vdW}}$,
${\utphi}{}^{(n)}_{\betab; t_{1, \nu}}$, ${\utphi}{}^{(n)}_{\betab;
\mathrm{S}}$, ${\utphi}{}^{(n)}_{\betab; \mathrm{W}}$, and
${\utphi}{}^{(n)}_{\betab; \mathrm{SP}}$, respectively. Sample size is
$n=100$. All tests were based on asymptotic 5\% critical values}
\label{simuresu2}
{\fontsize{8.9}{13}\selectfont{
\begin{tabular*}{\tablewidth}{@{\extracolsep{4in minus 4in}}lcccccccc@{}}
\hline
& \multicolumn{8}{c@{}}{$\bolds\xi$}\\[-4pt]
& \multicolumn{8}{c@{}}{\hrulefill}\\
\textbf{Test} & \textbf{0} & \textbf{1} & \textbf{2}
& \textbf{3} & \textbf{0} & \textbf{1} & \textbf{2} & \textbf{3} \\
\hline
%%%%%%%%%%%%%%%%%%%%%%%%%%%%%%%%%%%%%%%%%%%%%%
& \multicolumn{4}{c}{$\mathcal{N}$} &
\multicolumn{4}{c@{}}{$t_{5}$}\\[-4pt]
& \multicolumn{4}{c}{\hrulefill} & \multicolumn{4}{c@{}}{\hrulefill}\\
$\phi^{(n)}_{\betab; \mathrm{Anderson}}$ & 0.0572 & 0.3964 & 0.8804 &
0.9852 & 0.2408 & 0.4940 & 0.8388 & 0.9604 \\[2pt]
$\phi^{(n)}_{\betab; \mathcal{N}}$ & 0.0528 & 0.3724 & 0.8568 & 0.9752
& 0.2284 & 0.4716 & 0.8168 & 0.9380 \\[2pt]
$\phi^{(n)}_{\betab; \mathrm{Tyler}}$ & 0.0572 & 0.3908 & 0.8740 &
0.9856 & 0.0612 & 0.2520 & 0.6748 & 0.8876 \\[2pt]
$\phi^{(n)}_{\betab; \mathcal{N}*}$ & 0.0524 & 0.3648 & 0.8512 &
0.9740 & 0.0544 & 0.2188 & 0.6056 & 0.8156 \\[2pt]
${\utphi}{}^{(n)}_{\betab; \mathrm{vdW}}$ & 0.0368 & 0.2960 & 0.8032 &
0.9608 & 0.0420 & 0.2328 & 0.6908 & 0.9056 \\
${\utphi}{}^{(n)}_{\betab; t_{1,5}}$ & 0.0452 & 0.3204
& 0.8096 & 0.9596 & 0.0476 & 0.2728 & 0.7440 & 0.9284 \\
${\utphi}{}^{(n)}_{\betab; t_{1,3}}$ & 0.0476 & 0.3104 & 0.7988 &
0.9532 & 0.0496 & 0.2760 & 0.7476 & 0.9280\\
${\utphi}{}^{(n)}_{\betab; t_{1,1}}$ & 0.0488 & 0.2764 & 0.7460 &
0.9220 & 0.0552 & 0.2652 & 0.7184 & 0.9024 \\
${\utphi}{}^{(n)}_{\betab; \mathrm{S}}$ & 0.0448 & 0.2268 & 0.6204 & 0.8392
& 0.0496 & 0.2164 & 0.6236 & 0.8324 \\
${\utphi}{}^{(n)}_{\betab; \mathrm{W}}$ & 0.0456 & 0.3144 & 0.8012 & 0.9556
& 0.0484 & 0.2808 & 0.7464 & 0.9320 \\
${\utphi}{}^{(n)}_{\betab; \mathrm{SP}}$ & 0.0444 & 0.3096 & 0.8160 & 0.9576
& 0.0464 & 0.2548 & 0.7068 & 0.9152 \\
[2pt]
& \multicolumn{4}{c}{$t_{3}$} & \multicolumn{4}{c@{}}{$t_{1}$}\\[-4pt]
& \multicolumn{4}{c}{\hrulefill} & \multicolumn{4}{c@{}}{\hrulefill}\\
%%%%%%%%%%%%%%%%%%%%%%%%%%%%%%%%%%%%%%%%%%%%%%
%%%%%%%%%%%%%%%%%%%%%%%%%%%%%%%%%%%%%%%%%%%%%%
$\phi^{(n)}_{\betab; \mathrm{Anderson}}$ & 0.4772 & 0.6300 & 0.8532 &
0.9452 & 0.9540 & 0.9580 & 0.9700 & 0.9740 \\[2pt]
$\phi^{(n)}_{\betab; \mathcal{N}}$ & 0.4628 & 0.6040 & 0.8304 & 0.9168
& 0.9320 & 0.9384 & 0.9472 & 0.9480 \\[2pt]
$\phi^{(n)}_{\betab; \mathrm{Tyler}}$ & 0.0892 & 0.2248 & 0.5364 &
0.7508 & 0.5704 & 0.5980 & 0.6584 & 0.7444 \\[2pt]
$\phi^{(n)}_{\betab; \mathcal{N}*}$ & 0.0616 & 0.1788 & 0.4392 &
0.6092 & 0.4516 & 0.4740 & 0.5160 & 0.5624 \\[2pt]
${\utphi}{}^{(n)}_{\betab; \mathrm{vdW}}$ & 0.0444 & 0.2172 & 0.6464 &
0.8676 & 0.0472 & 0.1656 & 0.5104 & 0.7720 \\
${\utphi}{}^{(n)}_{\betab; t_{1,5}}$ & 0.0488 & 0.2628 &
0.7120 & 0.9076 & 0.0560 & 0.2100 & 0.6068 & 0.8508 \\
${\utphi}{}^{(n)}_{\betab; t_{1,3}}$ & 0.0500 & 0.2728 & 0.7156 &
0.9116 & 0.0576 & 0.2156 & 0.6292 & 0.8672 \\
${\utphi}{}^{(n)}_{\betab; t_{1,1}}$ & 0.0476 & 0.2688 & 0.7100 &
0.9084 & 0.0548 & 0.2256 & 0.6600 & 0.8856 \\
${\utphi}{}^{(n)}_{\betab; \mathrm{S}}$ & 0.0492 & 0.2202 & 0.6188 & 0.8352
& 0.0512 & 0.2116 & 0.6172 & 0.8448 \\
${\utphi}{}^{(n)}_{\betab; \mathrm{W}}$ & 0.0520 & 0.2708 & 0.7136 & 0.9120
& 0.0552 & 0.2148 & 0.6148 & 0.8604 \\
${\utphi}{}^{(n)}_{\betab; \mathrm{SP}}$ & 0.0544 & 0.2436 & 0.6648 & 0.8776
& 0.0580 & 0.1824 & 0.5200 & 0.7740 \\
\hline
\end{tabular*}}}
\end{table}

%t3 ###
%
\begin{table}
\caption{Rejection frequencies (out of $N=2500$ replications), under
the null ${\mathcal H}_{0}^{\Lamb}$ and increasingly distant
alternatives (see Section \protect\ref{secsimu}), %for details
of the optimal Gaussian test $\phi^{(n)}_{\Lamb; \mathcal{N}}
= \phi^{(n)}_{\Lamb; \mathrm{Anderson}}$,
its pseudo-Gaussian version $\phi^{(n)}_{\Lamb; \mathcal{N}*} =
\phi^{(n)}_{ \Lamb; \mathrm{Davis}}$, and the signed-rank
tests with van der Waerden, $t_\nu$ ($\nu=
1$, $3$, $5$), sign, Wilcoxon, and Spearman sco\-res ${\utphi
}{}^{(n)}_{\Lamb; \mathrm{vdW}}$, ${\utphi}{}^{(n)}_{\Lamb; t_{1, \nu
}}$, ${\utphi}{}^{(n)}_{\Lamb; \mathrm{S}}$, ${\utphi}{}^{(n)}_{\Lamb;
\mathrm{W}}$, ${\utphi}{}^{(n)}_{\Lamb; \mathrm{SP}}$. Sample size is
$n =
100$. All tests were based on asymptotic 5\%
critical values and (in parentheses) simulated ones}
\label{simuresu3}
\begin{tabular*}{\tablewidth}{@{\extracolsep{\fill}}lllll@{}}
\hline
& \multicolumn{4}{c@{}}{$\bolds\xi$}\\[-4pt]
& \multicolumn{4}{c@{}}{\hrulefill}\\
\textbf{Test} & \multicolumn{1}{c}{\textbf{0}}
& \multicolumn{1}{c}{\textbf{1}}
& \multicolumn{1}{c}{\textbf{2}}
& \multicolumn{1}{c@{}}{\textbf{3}} \\
\hline
& \multicolumn{4}{c@{}}{$\mathcal{N}$}\\[-4pt]
& \multicolumn{4}{c@{}}{\hrulefill}\\
%%%%%%%%%%%%%%%%%%%%%%%%%%%%%%%%%%%%%%%%%%%%%%
$\phi^{(n)}_{\Lamb; \mathcal{N}}=\phi^{(n)}_{\Lamb; \mathrm
{Anderson}}$ & 0.0460 & 0.4076 & 0.8308 & 0.9604 \\[2pt]
$\phi^{(n)}_{\Lamb; \mathcal{N}*}=\phi^{(n)}_{ \Lamb; \mathrm{
Davis}}$ & 0.0432 & 0.3976 & 0.8220 & 0.9572 \\[2pt]
${\utphi}{}^{(n)}_{\Lamb; \mathrm{vdW}}$ & 0.0608 (0.0480) &
0.4604 (0.4116) &
0.8576 (0.8280) & 0.9668 (0.9596) \\
${\utphi}{}^{(n)}_{\Lamb; t_{1,5}}$ & 0.0728 (0.0480) & 0.4804 (0.3972) &
0.8572 (0.8116) & 0.9644 (0.9504) \\
${\utphi}{}^{(n)}_{\Lamb; t_{1,3}}$ &  0.0748 (0.0496) &
0.4804 (0.3884) & 0.8524 (0.7964) & 0.9612 (0.9432) \\
${\utphi}{}^{(n)}_{\Lamb; t_{1,1}}$ & 0.0780 (0.0504) & 0.4532 (0.3572) &
0.8160 (0.7320) & 0.9448 (0.9112) \\
${\utphi}{}^{(n)}_{\Lamb; \mathrm{S}}$ & 0.0864 (0.0508) &
0.3980 (0.3088) &
0.7384 (0.6408) & 0.9028 (0.8552) \\
${\utphi}{}^{(n)}_{\Lamb; \mathrm{W}}$ & 0.0744 (0.0480) &
0.4816 (0.3908) &
0.8544 (0.8012) & 0.9640 (0.9464) \\
${\utphi}{}^{(n)}_{\Lamb; \mathrm{SP}}$ & 0.0636 (0.0460) &
0.4664 (0.4096) &
0.8564 (0.8200) & 0.9668 (0.9584) \\
& \multicolumn{4}{c@{}}{$t_{5}$}\\[-4pt]
& \multicolumn{4}{c@{}}{\hrulefill}\\
%%%%%%%%%%%%%%%%%%%%%%%%%%%%%%%%%%%%%%%%%%%%%%
$\phi^{(n)}_{\Lamb; \mathcal{N}}=\phi^{(n)}_{\Lamb; {\mathrm
{Anderson}}}$ & 0.1432 & 0.4624 & 0.7604 & 0.9180 \\[2pt]
$\phi^{(n)}_{\Lamb; \mathcal{N}*}=\phi^{(n)}_{ \Lamb; \mathrm{
Davis}}$ & 0.0504 & 0.2768 & 0.5732 & 0.7988 \\[2pt]
${\utphi}{}^{(n)}_{\Lamb; \mathrm{vdW}}$ & 0.0692 (0.0548) &
0.4256 (0.3772) &
0.7720 (0.7404) & 0.9444 (0.9324) \\
${\utphi}{}^{(n)}_{\Lamb; t_{1,5}}$ & 0.0736 (0.0492) & 0.4544 (0.3772) &
0.7980 (0.7372) & 0.9524 (0.9332) \\
${\utphi}{}^{(n)}_{\Lamb; t_{1,3}}$ & 0.0732 (0.0452) &
0.4576 (0.3748) & 0.7968 (0.7320) & 0.9524 (0.9288) \\
${\utphi}{}^{(n)}_{\Lamb; t_{1,1}}$ & 0.0776 (0.0416) & 0.4448 (0.3484) &
0.7832 (0.6952) & 0.9436 (0.9116) \\
${\utphi}{}^{(n)}_{\Lamb; \mathrm{S}}$ & 0.0768 (0.0436) &
0.4060 (0.3172) &
0.7180 (0.6360) & 0.9100 (0.8592) \\
${\utphi}{}^{(n)}_{\Lamb; \mathrm{W}}$ & 0.0732 (0.0456) &
0.4512 (0.3756) &
0.7972 (0.7364) & 0.9524 (0.9308) \\
${\utphi}{}^{(n)}_{\Lamb; \mathrm{SP}}$ & 0.0764 (0.0544) &
0.4360 (0.3736) &
0.7776 (0.7304) & 0.9480 (0.9300) \\
\hline
\end{tabular*}
\end{table}

\setcounter{table}{2}
\begin{table}
\caption{(Continued.)}
\begin{tabular*}{\tablewidth}{@{\extracolsep{\fill}}lllll@{}}
\hline
& \multicolumn{4}{c@{}}{$\bolds\xi$}\\[-4pt]
& \multicolumn{4}{c@{}}{\hrulefill}\\
\textbf{Test} & \multicolumn{1}{c}{\textbf{0}}
& \multicolumn{1}{c}{\textbf{1}} & \multicolumn{1}{c}{\textbf{2}}
& \multicolumn{1}{c@{}}{\textbf{3}} \\
\hline
& \multicolumn{4}{c@{}}{$t_{3}$}\\[-4pt]
& \multicolumn{4}{c@{}}{\hrulefill}\\
%%%%%%%%%%%%%%%%%%%%%%%%%%%%%%%%%%%%%%%%%%%%%%
$\phi^{(n)}_{\Lamb; \mathcal{N}}=\phi^{(n)}_{\Lamb; \mathrm
{Anderson}}$ & 0.2572 & 0.5308 & 0.7200 & 0.8596 \\[2pt]
$\phi^{(n)}_{\Lamb; \mathcal{N}*}=\phi^{(n)}_{ \Lamb; \mathrm{
Davis}}$ & 0.0368 & 0.1788 & 0.3704 & 0.5436 \\[2pt]
${\utphi}{}^{(n)}_{\Lamb; \mathrm{vdW}}$ & 0.0708 (0.0560) &
0.4088 (0.3260) &
0.7540 (0.7040) & 0.9304 (0.9140) \\
${\utphi}{}^{(n)}_{\Lamb; t_{1,5}}$ & 0.0812 (0.0544) & 0.4472 (0.3524) &
0.7936 (0.7240) & 0.9416 (0.9208) \\
${\utphi}{}^{(n)}_{\Lamb; t_{1,3}}$ & 0.0832 (0.0560) &
0.4556 (0.3568) & 0.7944 (0.7256) & 0.9452 (0.9192) \\
${\utphi}{}^{(n)}_{\Lamb; t_{1,1}}$ & 0.0924 (0.0548) & 0.4464 (0.3400) &
0.7812 (0.7024) & 0.9364 (0.8996) \\
${\utphi}{}^{(n)}_{\Lamb; \mathrm{S}}$ & 0.0936 (0.0604) &
0.4104 (0.2928) &
0.7320 (0.6404) & 0.9012 (0.8528) \\
${\utphi}{}^{(n)}_{\Lamb; \mathrm{W}}$ & 0.0832 (0.0572) &
0.4488 (0.3580) &
0.7956 (0.7272) & 0.9448 (0.9180) \\
${\utphi}{}^{(n)}_{\Lamb; \mathrm{SP}}$ & 0.0796 (0.0576) &
0.4212 (0.3412) &
0.7572 (0.7020) & 0.9276 (0.9044) \\
& \multicolumn{4}{c@{}}{$t_{1}$}\\[-4pt]
& \multicolumn{4}{c@{}}{\hrulefill}\\
%%%%%%%%%%%%%%%%%%%%%%%%%%%%%%%%%%%%%%%%%%%%%%
$\phi^{(n)}_{\Lamb; \mathcal{N}}=\phi^{(n)}_{\Lamb; \mathrm
{Anderson}}$ & 0.7488 & 0.8000 & 0.8288 & 0.8528 \\[2pt]
$\phi^{(n)}_{\Lamb; \mathcal{N}*}=\phi^{(n)}_{ \Lamb; \mathrm{
Davis}}$ & 0.0072 & 0.0080 & 0.0172 & 0.0296 \\[2pt]
${\utphi}{}^{(n)}_{\Lamb; \mathrm{vdW}}$ & 0.0724 (0.0596) &
0.3500 (0.3032) &
0.6604 (0.6176) & 0.8600 (0.8332) \\
${\utphi}{}^{(n)}_{\Lamb; t_{1,5}}$ & 0.0824 (0.0512) & 0.3836 (0.3120) &
0.7312 (0.6492) & 0.9036 (0.8664) \\
${\utphi}{}^{(n)}_{\Lamb; t_{1,3}}$ & 0.0828 (0.0532)&
0.3936 (0.3108) & 0.7488 (0.6644) & 0.9168 (0.8776) \\
${\utphi}{}^{(n)}_{\Lamb; t_{1,1}}$ & 0.0864 (0.0532) & 0.4088 (0.3104) &
0.7612 (0.6720) & 0.9264 (0.8832) \\
${\utphi}{}^{(n)}_{\Lamb; \mathrm{S}}$ & 0.0920 (0.0556) &
0.3896 (0.3028) &
0.7336 (0.6488) & 0.9092 (0.8564) \\
${\utphi}{}^{(n)}_{\Lamb; \mathrm{W}}$ & 0.0824 (0.0524) &
0.3872 (0.3072) &
0.7376 (0.6552) & 0.9108 (0.8728) \\
${\utphi}{}^{(n)}_{\Lamb; \mathrm{SP}}$ & 0.0752 (0.0588) &
0.3536 (0.2992) &
0.6648 (0.6064) & 0.8604 (0.8220) \\
\hline
\end{tabular*}\vspace*{-3pt}
\end{table}

Inspection of Table~\ref{simuresu2} confirms our theoretical results.
Anderson's $\phi_{ \betab;\mathrm{Anderson}}^{(n)}$ meets the level
constraint at\vspace*{2pt} Gaussian densities only; $\phi^{(n)}_{\betab;\mathrm
{Tyler}}$ (equivalently, $\phi^{(n)}_{\betab; \mathcal{N}*}$)
further survives the $t_5$ but not the $t_3$ or $t_1$ densities which
have infinite fourth-order moments. In contrast, the rank-based tests
for eigenvectors throughout satisfy the nominal asymptotic level
condition (a 95\% confidence interval here has half-width 0.0085).
Despite the relatively small sample size $n=100$, empirical power and
ARE rankings almost perfectly agree.

The results for eigenvalues, shown in Table~\ref{simuresu3}, are
slightly less auspicious. While the Gaussian and pseudo-Gaussian tests
remain hopelessly sensitive to the violations of Gaussian and
fourth-order moments, respectively, the rank tests, when based on
asymptotic critical values, all significantly overreject, indicating
that asymptotic conditions are not met for $n=100$. We therefore
propose an alternative construction for critical values. Lemma \ref
{alihoprimeprime} indeed implies that the asymptotic distribution of
the test statistic ${\utT}{}^{(n)}_{K}$ (based\vspace*{-3pt} on the ranks and signs
of estimated residuals) is the same, under $\mathrm{P}^{(n)}_{\varthetab
_0; g_{1}}$, $\varthetab_0 \in{\mathcal H}_{0;q}^{\Lamb\prime\prime
}$, as that of ${\utT}{}_{\varthetab_0, K}^{(n)}$ (based on the ranks
and signs of \textit{exact} residuals, which are distribution-free). The
latter distribution can be simulated, and its simulated quantiles
provide valid approximations of the exact ones. The following critical
values were obtained from $M=100$,000 replications:\vadjust{\goodbreak} $-1.7782$ for van
der Waerden, $-1.8799$ for $t_5$-scores, $-1.8976$ for $t_3$-scores,
$-1.9439$ for $t_1$-scores, $-1.9320$ for sign scores, $-1.8960$ for
Wilcoxon and $-1.8229$ for Spearman. Note that they all are smaller
than $-1.645$, which is consistent with the overrejection phenomenon.
The corresponding rejection frequencies are reported in parentheses in
Table~\ref{simuresu3}. They all are quite close to the nominal
probability level $\alpha= 5\%$, while empirical powers are in line
with theoretical ARE values.

\begin{appendix}\label{app}
\section*{Appendix}

We start with the proof of Proposition~\ref{LAN}. To this end, note
that although generally stated as a property of a parametric sequence
of families of the form $\mathcal{P}^{(n)}= \{\mathrm{P}^{(n)}_{\bolds
\omega} \vert{\bolds\omega}\in{\bolds\Omega}\}$ ($n\in\N$),
LAN (ULAN)\vadjust{\goodbreak} actually is a property of the parametrization ${\bolds
\omega}\mapsto\mathrm{P}^{(n)}_{\bolds\omega}$, ${\bolds\omega}\in
{\bolds\Omega}$ of $\mathcal{P}^{(n)}$ (i.e., of a bijective map
from $\Omegab$ to $\mathcal{P}^{(n)}$). When parametrized with
${\bolds\omega}:=({\bolds\theta}\pr, (\operatorname{vech} \Sigb
)\pr)\pr$, ${\bolds\omega}\in\Omegab:=\R^k\times\operatorname
{vech}(\mathcal{S}_k)$, where $\mathcal{S}_k$ stands for the class of
positive definite symmetric real $k\times k$ matrices, the elliptical
families we are dealing with
%in this paper
here have been shown to be ULAN in \citet{HP06a}, with
central sequence
%
%e8.3 ###
%
\setcounter{equation}{0}
\begin{equation} \label{HPcentral}
\Deltab_{{\bolds\omega}}^{(n)}:=
\pmatrix{\displaystyle n^{-1/2}
\sum_{i=1}^{n} \varphi_{f_1}(d_{i})
{\Sigb}^{-1/2}\mathbf{U}_{i} \vspace*{2pt}\cr
\displaystyle\frac{ 1}{2\sqrt{n}}\mathbf{P}_{k}
( {\Sigb}^{\otimes2} )^{ -1/2} \sum_{i=1}^{n}
\vecop\bigl( \varphi_{f_1} (d_{i}) d_{i} \mathbf{U}_{i} \mathbf{U}_{i} \pr
-\mathbf{I}_{k}\bigr)},
\end{equation}
with $d_{i} =d_{i}({\bolds\theta},\Sigb)$ and $\mathbf{U}_{i}=\mathbf
{U}_{i}({\bolds\theta},\Sigb)$, where $\mathbf{P}_{k}\pr$ denotes the
\textit{duplication matrix} [such that $\mathbf{P}_{k}\pr\vechop(\mathbf
{A})=\vecop(\mathbf{A})$ for any $k \times k$ symmetric matrix $\mathbf{A}$].

The families we are considering in this proposition are slightly
different, because the $\varthetab$-parametrization requires $k$
identifiable eigenvectors. However, denoting by $\Omegab^B:=\R
^k\times\operatorname{vech}(\mathcal{S}_k^B)$, where $\mathcal
{S}^B_k$ is the set of all matrices in $\mathcal{S}_k$ compatible with
Assumption~\ref{assuA}, the
mapping
\begin{eqnarray*}
&&\dbar\dvtx{\bolds\omega} = ({\bolds\theta}\pr, (\operatorname{vech}
\Sigb)\pr)\pr\in\Omegab^B \\
&&\quad\mapsto\quad\dbar({\bolds\omega}) := \bigl({\bolds\theta}\pr,
(\operatorname{det} \Sigb)^{1/k} , (\dvecrond\Lamb_\Sigb)\pr/
(\operatorname{det} \Sigb)^{1/k}, (\vecop\betab)\pr\bigr)\pr\in
\Thetab
\end{eqnarray*}
from the open subset $\Omegab^B$ of $\R^{k+k(k+1)/2}$ to $\Thetab$
is a differentiable mapping such that, with a small abuse of notation,
$\mathrm{P}^{(n)}_{{\bolds\omega}; f_1} = \mathrm{P}^{(n)}_{\varthetab
=\hspace*{1pt}\dbarr({\bolds\omega}); f_1} $ and $\mathrm{P}^{(n)}_{\varthetab; f_1}
= \mathrm{P}^{(n)}_{{\bolds\omega}=\hspace*{1pt}\dbarr^{-1} (\varthetab); f_1} $
(with $\hspace*{2pt}\dbar^{-1}$ defined on $\Thetab$ only). The proof of
Proposition~\ref{LAN} consists in showing how ULAN in the ${\bolds
\omega}$-parametrization implies ULAN in the $\varthetab
$-parametrization, and how the central sequences and information
matrices are related to each other. Let us start with a lemma.
\begin{Lem}\label{LElemme}
Let the parametrization ${\bolds\omega}\mapsto\mathrm{P}^{(n)}_{\bolds
\omega}$, ${\bolds\omega}\in{\bolds\Omega}$, where $\Omegab$ is
an open subset of $\R^{k_1}$ be ULAN for $\mathcal{P}^{(n)}= \{\mathrm
{P}^{(n)}_{\bolds\omega} \vert{\bolds\omega}\in{\bolds\Omega
}\}$, with central sequence $\Deltab^{(n)}_{\bolds\omega}$ and
information matrix $\Gamb_{\bolds\omega}$. Let $\dbar\dvtx{\bolds
\omega}\mapsto\varthetab:= \dbar({\bolds\omega})$ be a
continuously differentiable mapping from $\R^{k_1}$ to $\R^{k_2}$
($k_2\geq k_1$) with full column rank Jacobian matrix $D\dbar({\bolds
\omega})$ at every ${\bolds\omega}$. Write $\Thetab:= \dbar
(\Omegab)$, and assume that $\varthetab\mapsto\mathrm
{P}^{(n)}_\varthetab$, $\varthetab\in{\Thetab}$ provides another
parametrization of $\mathcal{P}^{(n)}$.
Then, $\varthetab\mapsto\mathrm{P}^{(n)}_\varthetab$, $\varthetab\in
{\Thetab}$ is also ULAN,
with [at $\varthetab=\dbar({\bolds\omega})$] central sequence
$\Deltab^{(n)}_\varthetab= (D^-\dbar({\bolds\omega}) )\pr\Deltab
^{(n)}_{\bolds\omega}$
and information matrix $\Gamb_\varthetab= (D^-\dbar({\bolds\omega
}))\pr\Gamb_{\bolds\omega}D^-\dbar({\bolds\omega})$,
where $D^-\dbar({\bolds\omega}):=((D\dbar({\bolds\omega}) )\pr
D\dbar({\bolds\omega}))^{-1}(D\dbar({\bolds\omega}))\pr$ is the
Moore--Penrose inverse of $D\dbar({\bolds\omega})$.
\end{Lem}
\begin{pf}
Throughout, let $\varthetab$ and ${\bolds\omega}$ be such that
$\varthetab= \dbar({\bolds\omega})$. %ULAN for the $
Consider $\varthetab\in\Thetab$ and an arbitrary sequence
$\varthetab^{(n)}=\varthetab+ O(n^{-1/2})\in\Thetab$. The
characterization of ULAN for the $\varthetab$-parametrization involves
bounded sequence $\taub_{**}^{(n)}\in\R^{k_2}$ such that the
perturbation $\varthetab^{(n)}+ n^{-1/2}\taub_{**}^{(n)}$ still
belongs to $\Thetab$.
%$$\log(d\mathrm{P}\n_{\varthetab\n+ n^{-1/2}\taub_{**}\n})d
%$$
In order for $\varthetab^{(n)}+ n^{-1/2}\taub_{**}^{(n)}$ to belong
to $\Thetab$, it is necessary that $\taub_{**}^{(n)}$ be of the form
$\taub_*^{(n)}+ o(1)$, with $\taub_*^{(n)}$ in the tangent space to
$\Thetab$ at $\varthetab^{(n)}$, hence of the form $\taub^{(n)}+
o(1)$ with $\taub^{(n)}$ in the tangent space to $\Thetab$ at
$\varthetab$, that is, $\taub^{(n)}=D\dbar({\bolds\omega})\wb
^{(n)}$ for some bounded sequence $\wb^{(n)}\in\R^{k_1}$. It follows
from differentiability that, letting ${\bolds\omega}^{(n)}=\dbar
^{-1}(\varthetab^{(n)})$,
%
%e8.4 ###
%
\begin{eqnarray} \label{c'estlui}
\varthetab^{(n)} + n^{-1/2}\taub_{**}^{(n)} & = &
\varthetab^{(n)}+ n^{-1/2}D\dbar({\bolds\omega})\wb^{(n)}+ o(n^{-1/2})
\nonumber\\
& = &
\dbar\bigl({\bolds\omega}^{(n)}\bigr) + n^{-1/2}D\dbar({\bolds\omega})\wb
^{(n)} + o(n^{-1/2}) \nonumber\\[-8pt]\\[-8pt]
& = & \dbar\bigl({\bolds\omega}^{(n)}\bigr) + n^{-1/2}D\dbar\bigl({\bolds
\omega}^{(n)}\bigr)\wb^{(n)}+ o(n^{-1/2}) \nonumber\\
& = &
\dbar\bigl({\bolds\omega}^{(n)}+ n^{-1/2} \wb^{(n)}+ o(n^{-1/2})\bigr).
\nonumber
\end{eqnarray}
Hence, turning to local log-likelihood ratios, in view of ULAN for the
${\bolds\omega}$-parametrization,
\begin{eqnarray}\label{onyest}
&&\log\bigl(d\mathrm{P}{}_{\varthetab^{(n)}+ n^{-1/2}\taub
_{**}^{(n)}}^{(n)}/d\mathrm{P}^{(n)}_{\varthetab^{(n)}} \bigr)
\nonumber\\
&&\qquad=
\log\bigl(d\mathrm{P}^{(n)}_{ {\bolds\omega}^{(n)}+ n^{-1/2} \wb
^{(n)}+ o(n^{-1/2})}/d\mathrm{P}^{(n)}_{{\bolds\omega}^{(n)}} \bigr)
\\
&&\qquad= \wb^{(n)\prime}\Deltab^{(n)}_{{\bolds\omega}^{(n)}} -\tfrac
{1}{2}
\wb^{(n)\prime}\Gamb_{\bolds\omega}\wb^{(n)}+o_\mathrm{P}(1)\nonumber
%&=&
\end{eqnarray}
under
$\mathrm{P}^{(n)}_{ {\bolds\omega}^{(n)}} = \mathrm{P}^{(n)}_{
\varthetab
^{(n)}}$-probability,
as $\ny$.
Now, the LAQ part of ULAN for the $\varthetab$-parametrization
requires, for some random vector $\Deltab^{(n)}_{\varthetab^{(n)}}$
and constant matrix $\Gamb_\varthetab$,
%
%e8.5 ###
%
\begin{equation}\label{LAQvartheta}\qquad
\log\bigl(d\mathrm{P}^{(n)}_{\varthetab^{(n)}+ n^{-1/2}\taub
_{**}^{(n)}}/d\mathrm{P}^{(n)}_{\varthetab^{(n)}} \bigr) = \taub
_{**}^{(n)\prime}\Deltab^{(n)}_{\varthetab^{(n)}}-\tfrac{1}{2}
\taub_{**}^{(n)\prime}\Gamb_\varthetab\taub_{**}^{(n) } +o_\mathrm{P}(1)
\end{equation}
under the same
$\mathrm{P}^{(n)}_{ {\bolds\omega}^{(n)}} = \mathrm{P}^{(n)}_{
\varthetab
^{(n)}}$ probability distributions with, in view of~(\ref{c'estlui}),
$\taub_{**}^{(n)}= D\dbar({\bolds\omega})\wb^{(n)} + o(1)$.
Identifying~(\ref{onyest}) and~(\ref{LAQvartheta}), we obtain that
LAQ is satisfied for the $\varthetab$-parametrization, with
any $\Deltab^{(n)}_\varthetab$ satisfying
%
%e8.6 ###
%
\begin{equation}\label{LAQvartheta2}
(D\dbar({\bolds\omega}))\pr\Deltab^{(n)}_\varthetab= \Deltab
^{(n)}_{\bolds\omega}.
\end{equation}

Now, let $\mathbf{t}_i$ be the $i$th column of $D\dbar({\bolds\omega
})$, $i=1,\ldots,k_1$, and choose $\mathbf{t}_{k_1+1},\ldots,\break
\mathbf{t}_{k_2}\in\mathbb{R}^{k_2}$ in such a way that they span the
orthogonal complement of $\mathcal{M}(D\dbar({\bolds\omega}))$.
Then $\{\mathbf{t}_{i}$, $i=1,\ldots,k_2\}$ is a basis of $ \mathbb
{R}^{k_2}$, so that there exists a unique $k_2$-tuple $(\delta
_{\varthetab;1}^{(n)},\ldots, \delta^{(n)}_{\varthetab; k_2})\pr$
such that $\Deltab^{(n)}_\varthetab=\sum_{i=1}^{k_2} \delta
^{(n)}_{\varthetab;i} \mathbf{t}_i$. With this notation,~(\ref{LAQvartheta2})
yields
\begin{eqnarray*}
\Deltab^{(n)}_{\bolds\omega}
&=&(D\dbar({\bolds\omega}))\pr\Deltab^{(n)}_\varthetab\\
&=&\sum_{i=1}^{k_2} \delta^{(n)}_{\varthetab;i} (D\dbar({\bolds
\omega}))\pr\mathbf{t}_i
= \sum_{i=1}^{k_1} \delta^{(n)}_{\varthetab;i} (D\dbar({\bolds
\omega}))\pr\mathbf{t}_i\\
&=& (D\dbar({\bolds\omega}))\pr D\dbar({\bolds\omega}) \underline
\Deltab^{(n)}_{\varthetab},
\end{eqnarray*}
where we let $\underline\Deltab^{(n)}_{\varthetab}:=(\delta
^{(n)}_{\varthetab;1} ,\ldots,\delta^{(n)}_{\varthetab;k_1})\pr$.
Since $D\dbar({\bolds\omega})$ has full column rank, this entails (i)
$\Deltab^{(n)}_\varthetab=D\dbar({\bolds\omega}) \underline
\Deltab^{(n)}_{\varthetab}$ and (ii) %\vadjust
$\underline\Deltab^{(n)}_{\varthetab}=((D\dbar({\bolds\omega
}))\pr D\dbar({\bolds\omega}))^{-1} \Deltab^{(n)}_{\bolds\omega
}$, hence $\Deltab^{(n)}_\varthetab=(D^-\dbar({\bolds\omega}))\pr
\Deltab^{(n)}_{\bolds\omega}$. As a linear transformation of
$\Deltab^{(n)}_{\bolds\omega}$, $\Deltab^{(n)}_\varthetab$ clearly
also satisfies the asymptotic normality part of ULAN, with the desired
$\Gamb_\varthetab$.
\end{pf}

The following slight extension of Lemma~\ref{LElemme} plays a role in
the proof of Proposition~\ref{LAN} below. Consider a parametrization
${\bolds\omega}=({\bolds\omega}_a\pr,{\bolds\omega}_b\pr)\pr
\mapsto\mathrm{P}^{(n)}_{\bolds\omega}$, ${\bolds\omega}\in{\bolds
\Omega}\times\mathcal{V}$, where $\Omegab$ is an open subset of $\R
^{k_1}$ and $\mathcal{V}\subset\R^m$ is a $\ell$-dimensional
manifold in $\R^m$, and assume that it is ULAN for $\mathcal
{P}^{(n)}= \{\mathrm{P}^{(n)}_{\bolds\omega} \vert{\bolds\omega
}\in{\bolds\Omega}\times\mathcal{V} \}$, with central sequence
$\Deltab^{(n)}_{\bolds\omega}$ and information matrix $\Gamb
_{\bolds\omega}$. Let $\dbar_a$ be a continuously differentiable
mapping from $\R^{k_1}$ to $\R^{k_2}$ ($k_2\geq k_1$) with full
column rank Jacobian matrix $D\dbar_a ({\bolds\omega}_a)$ at
every~${\bolds\omega}_a$, and assume that $\varthetab:= \dbar({\bolds
\omega})\mapsto\mathrm{P}^{(n)}_\varthetab$, $\varthetab\in{\Thetab}
\times\mathcal{V}$ [with $ {\Thetab} := \dbar_a(\Omegab)$], where
\[
\dbar\dvtx
\Omegab\times\mathcal{V} \to\Thetab\times\mathcal{V}\qquad
{\bolds\omega}=
( {\bolds\omega}_a,
{\bolds\omega}_b
)\pr
\quad\mapsto\quad\dbar({\bolds\omega})=
(
\dbar_a({\bolds\omega}_a) ,
{\bolds\omega}_b
)\pr
\]
provides another parametrization of $\mathcal{P}^{(n)}$. Then the
proof of Lemma~\ref{LElemme} straightforwardly extends to show that
$\varthetab\mapsto\mathrm{P}^{(n)}_\varthetab$, $\varthetab\in
{\Thetab}\times\mathcal{V}$ is also ULAN, still with [at $\varthetab
=\dbar({\bolds\omega})$] central sequence $\Deltab^{(n)}_\varthetab
= (D^-\dbar({\bolds\omega}) )\pr\Deltab^{(n)}_{\bolds\omega}$
and information matrix $\Gamb_\varthetab= (D^-\dbar({\bolds\omega
}))\pr\Gamb_{\bolds\omega}D^-\dbar({\bolds\omega})$.
\begin{pf*}{Proof of Proposition~\ref{LAN}}
% (with an explicit computation of $\Deltab_\varthetab\n$ and $\Gamb_
Consider the differentiable mappings $\dbar_1 \dvtx\break{\bolds\omega
}:=({\bolds\theta}\pr, (\vechop\Sigb)\pr)\pr\mapsto\dbar_1
(\omega)=({\bolds\theta}', (\dvec\Lamb_{\Sigb})\pr, (\vecop
\betab)\pr)'$ and $\dbar_2 \dvtx\dbar_1 ({\bolds\omega})=({\bolds
\theta}'$,\break $(\dvec\Lamb_{\Sigb})\pr, (\vecop\betab)\pr)' \mapsto
\dbar_2 (\dbar_1 ({\bolds\omega}))=({\bolds\theta}', \sigma
^{2},(\dvecrond\Lamb_{\Vb})\pr, (\vecop\betab)\pr)' \in
\Thetab$, the latter being invertible.
Applying Lemma~\ref{LElemme} twice (the second time in its ``extended
form,'' since the $\betab$-part of
the parameter is invariant under $\dbar_2$) then yields
\begin{eqnarray*}
\Deltab_{\varthetab}^{(n)}
&=&
(D\dbar_2 (\dbar_1({\bolds\omega})))^{\prime-1} D \dbar_1
({\bolds\omega})((D \dbar_1 ({\bolds\omega}))\pr D \dbar_1
({\bolds\omega}))^{-1} \Deltab^{(n)}_{{\bolds\omega}} \\
&=&
(D \dbar_2^{-1} (\dbar({\bolds\omega})))\pr D \dbar_1 ({\bolds
\omega})((D \dbar_1 ({\bolds\omega}))\pr D \dbar_1 ({\bolds\omega
}))^{-1} \Deltab^{(n)}_{{\bolds\omega}} .
\end{eqnarray*}
In view of the definition of $\mathbf{M}_{k}^{\Lamb_{\Vb}}$
(Section~\ref{curvLAN}), the Jacobian matrix, computed at~$\varthetab
$, of the inverse mapping $\dbar_2 ^{-1}$ is
\[
D \dbar_2 ^{-1}( \varthetab)=\pmatrix{
\mathbf{I}_{k} & \mathbf{0} & \mathbf{0} & \mathbf{0} \cr\mathbf{0} &
\dvec(\Lamb_{\Vb}) & \sigma^{2} (\mathbf{M}_{k}^{\Lamb_{\Vb}})
\pr& \mathbf{0} \cr\mathbf{0} & \mathbf{0} & \mathbf{0} & \mathbf{I}_{k^{2}}}.
\]
An explicit expression for $D \dbar_1 ({\bolds\omega})$ was obtained
by Kollo and Neudecker [(\citeyear{KN93}), page 288]:
%
%e8.7 ###
%
\begin{eqnarray}
\label{ttttt}
D \dbar_1 ({\bolds\omega}) = \pmatrix{
\mathbf{I}_{k} & \mathbf{0} \cr\mathbf{0} & \Xib_{\betab, \Lamb_{\Sigb}}
\mathbf{P}_{k}\pr},\nonumber\\[-8pt]\\[-8pt]
\eqntext{\mbox{with }
\Xib_{\betab, \Lamb_{\Sigb}}:=
\pmatrix{\displaystyle
\mathbf{H}_{k}(\betab\pr)^{\otimes2} \cr
\displaystyle\betab_{1}\pr\otimes[\betab(\lambda_{1;\Sigb}\mathbf{I}_{k} - \Lamb
_{\Sigb})^{-}\betab\pr] \cr
\vdots\cr
\displaystyle\betab_{k}\pr\otimes[\betab(\lambda_{k;\Sigb}\mathbf{I}_{k} - \Lamb
_{\Sigb})^{-}\betab\pr]}.}
\end{eqnarray}
The result then follows from a direct, though painful, computation,
using the fact that
\[
(\mathbf{P}_{k}\Xib_{\betab, \Lamb_{\Sigb}}\pr\Xib_{\betab, \Lamb
_{\Sigb}}\mathbf{P}_{k}\pr)^{-1}= (\mathbf{P}_{k}\pr)^{-} (\betab
\otimes\betab) \operatorname{diag}(l_{11 ; \Sigb}, l_{12; \Sigb}, \ldots
, l_{kk; \Sigb}) (\betab\pr\otimes\betab\pr)\mathbf{P}_{k}^{-},
\]
with $l_{ij ; \Sigb}= 1$ if $i=j$ and $l_{ij; \Sigb}= (\lambda_{i;
\Sigb}- \lambda_{j; \Sigb})^{-2}$ if $i \neq j$; $(\mathbf{P}_{k}\pr
)^{-}$ here stands for the Moore--Penrose inverse of $\mathbf{P}_{k}$ [note
that $(\mathbf{P}_{k}\pr)^{-}$ is such that\break $\mathbf{P}_{k}\pr(\mathbf
{P}_{k}\pr)^{-} \vecop(\mathbf{A})= \vecop(\mathbf{A})$ for any symmetric
matrix $\mathbf{A}$].
\end{pf*}
%
%%%%%%%%%%%%%%%%%%%%%%%%%%%%%%%%%%%%%%%%%%%%%%%%%%%%%%%%%%%%%%%%%%%%%%%%%%%%%%%%%%%%%%%%%%%%%%%%%%%%%%%%%%%%%%%%%%%%%%%%%%%%%%%%
%%%%%%%%%%%%%%%%%%%%%%%%%%%%%%%%%%%%%%%%%%%%%%%%%%%%%%%%%%%%%%%%%%%%%%%%%%%%%%%%%%%%%%%%%%%%%%%%%%%%%%%%%%%%%%%%%%%%%%%%%%%%%%%%%
%%%%%%%%%%%%%%%%%%%%%%%%%%%%%%%%%%%%%%%%%%%%%%%%%%%%%%%%%%%%%%%%%%%%%%%%%%%%%%%%%%%%%%%%%%%%%%%%%%%%%%%%%%%%%%%%%%%%%%%%%%%%%%%%
\begin{pf*}{Proof of Proposition~\ref{Optitest}}
Proceeding as in the proof of Lemma~\ref{LElemme}, let $\mathbf{v}_i$ be
the $i$th column of $D\hbar({\bolds\xi}_0)$, $i=1,\ldots,m$, and
choose $\mathbf{v}_{m+1},\ldots,\mathbf{v}_p\in\mathbb{R}^p$ spanning the
orthogonal complement of $\mathcal{M}(D\hbar({\bolds\xi}_0))$. Then
there exists a unique $p$-tuple $(\delta_{\varthetab_0;1},\ldots,
\delta_{\varthetab_0; p})\pr$ such that $\Deltab_{{\varthetab
_0}}=\sum_{i=1}^p \delta_{\varthetab_0;i} \mathbf{v}_i$ (since $\mathbf
{v}_i, i=1,\ldots,p$ spans $\mathbb{R}^p$) and
%
%e8.8 ###
%
\begin{eqnarray} \label{ma}
\Deltab_{{\bolds\xi}_0}
&=& D\hbar\pr({\bolds\xi}_0)\Deltab_{\varthetab_0}
= \sum_{i=1}^p \delta_{\varthetab_0;i} D\hbar\pr({\bolds\xi}_0)
\mathbf{v}_i
= \sum_{i=1}^m \delta_{\varthetab_0;i} D\hbar\pr({\bolds\xi}_0)
\mathbf{v}_i\nonumber\\[-8pt]\\[-8pt]
&=& \mathbf{C}_{\hbar}({\bolds\xi}_0)
\Deltab^m_{\varthetab_0},\nonumber
\end{eqnarray}
where $\mathbf{C}_{\hbar}({\bolds\xi}_0):=D\hbar\pr({\bolds\xi}_0)
D\hbar({\bolds\xi}_0)$ %(of full rank $m$)
and $\Deltab^m_{\varthetab_0}:=(\delta_{\varthetab_0;1} ,\ldots
,\delta_{\varthetab_0;m})\pr$. Hence, we also have
$\Gamb_{{\bolds\xi}_0}=\mathbf{C}_{\hbar}({\bolds\xi}_0) \Gamb
^m_{\varthetab_0}\mathbf{C}_{\hbar}({\bolds\xi}_0)$, where $ \Gamb
^m_{\varthetab_0}$ is the asymptotic covariance matrix of $ \Deltab
^m_{\varthetab_0}$ under $\mathrm{P}^{(n)}_{\varthetab_0}$. Using
the fact that $\mathbf{C}_{\hbar}({\bolds\xi}_0)$ is invertible, this yields
\begin{eqnarray*}
Q_{{\bolds\xi}_0}:\!&=&
(\Deltab^m_{\varthetab_0})\pr\mathbf{C}_{\hbar}({\bolds\xi}_0)
(\mathbf{C}_{\hbar}({\bolds\xi}_0)\Gamb^m_{\varthetab_0}\mathbf
{C}_{\hbar}({\bolds\xi}_0))^{-1}
\mathbf{C}_{\hbar}({\bolds\xi}_0) \Deltab^m_{\varthetab_0}
\\
& &{}
-
(\Deltab^m_{\varthetab_0})\pr\mathbf{C}_{\hbar}({\bolds\xi}_0)
D\lbar({\bolds\alpha}_0)
( D\lbar\pr({\bolds\alpha}_0) \mathbf{C}_{\hbar}({\bolds\xi}_0)
\Gamb^m_{\varthetab_0} \mathbf{C}_{\hbar}({\bolds\xi}_0) D\lbar
({\bolds\alpha}_0))^{-1} \\
&&\hspace*{9.8pt}{}\times D\lbar\pr({\bolds\alpha}_0)
\mathbf{C}_{\hbar}({\bolds\xi}_0) \Deltab^m_{\varthetab_0}\\
&=& (\Deltab^m_{\varthetab_0})\pr
(\Gamb^m_{\varthetab_0})^{-1}
\Deltab^m_{\varthetab_0}\\
&&{}
-
(\Deltab^m_{\varthetab_0})\pr(\Gamb^m_{\varthetab_0})^{-1/2}
{\bolds\Pi}((\Gamb^m_{\varthetab_0})^{1/2} \mathbf{C}_{\hbar}({\bolds
\xi}_0) D\lbar({\bolds\alpha}_0))
(\Gamb^m_{\varthetab_0})^{-1/2} \Deltab^m_{\varthetab_0}\\
&=&\!: Q_{{\bolds\xi}_0,1}-Q_{{\bolds\xi}_0,2},
\end{eqnarray*}
where ${\bolds\Pi}(\mathbf{P}):=\mathbf{P}(\mathbf{P}\pr\mathbf
{P})^{-1}\mathbf{P}\pr$ denotes the projection matrix on $\mathcal
{M}(\mathbf{P})$.\vadjust{\goodbreak}

Let $\bbar\dvtx A\subset\mathbb{R}^\ell\to\mathbb{R}^p$ be a local
(at $\varthetab_0$) chart for the manifold $C\cap\Thetab$,
and assume, without loss of generality, that $\etab_0=\bbar
^{-1}(\varthetab_0)$.
Since $D\hbar({\bolds\xi}_0)$ has maximal rank, it follows
from~(\ref{ma}) that $\Deltab_{\varthetab_0}= D\hbar({\bolds\xi
}_0) \Deltab^m_{\varthetab_0}$. Hence, the statistic
%
%e8.9 ###
%
\begin{equation}\label{Qxibar}\quad
\bar{Q}_{\varthetab_0}
:=
\Deltab_{\varthetab_0}\pr
\bigl( \Gamb_{\varthetab_0}^{-} - D\bbar(\etab_0) ( D\bbar\pr(\etab
_0)\Gamb_{\varthetab_0}D\bbar(\etab_0))^{-} D\bbar\pr(\etab_0) \bigr)
\Deltab_{\varthetab_0}
\end{equation}
[the squared Euclidean norm of the orthogonal projection, onto the
linear space orthogonal to $ \Gamb_{\varthetab_0}^{1/2}D\bbar(\etab
_0)$, of the standardized central sequence $(\Gamb_{\varthetab
_0}^{1/2})^-\Deltab_{\varthetab_0}$] can be written as
\begin{eqnarray*}
\bar{Q}_{\varthetab_0}
&=&
(\Deltab^m_{\varthetab_0})\pr D\hbar\pr({\bolds\xi}_0)
( D\hbar({\bolds\xi}_0)\Gamb_{\varthetab_0}^m D\hbar\pr({\bolds
\xi}_0))^{-}
D\hbar({\bolds\xi}_0)\Deltab_{\varthetab_0}^m \\
& &{} -
(\Deltab_{\varthetab_0}^m)\pr D\hbar\pr({\bolds\xi}_0)
D\bbar(\etab_0) ( D\bbar\pr(\etab_0) D\hbar({\bolds\xi
}_0)\Gamb_{\varthetab_0}^m D\hbar\pr({\bolds\xi}_0) D\bbar(\etab
_0))^{-}\\
&&\hspace*{9.8pt}{}\times D\bbar\pr(\etab_0)
D\hbar({\bolds\xi}_0)\Deltab_{\varthetab_0}^m\\
&=&
(\Deltab^m_{\varthetab_0})\pr D\hbar\pr({\bolds\xi}_0)
( D\hbar({\bolds\xi}_0)\Gamb_{\varthetab_0}^m D\hbar\pr({\bolds
\xi}_0))^{-}
D\hbar({\bolds\xi}_0)\Deltab_{\varthetab_0}^m \\
& &{} -
(\Deltab_{\varthetab_0}^m)\pr(\Gamb^m_{\varthetab_0})^{-1/2}
{\bolds\Pi}( (\Gamb_{\varthetab_0}^m)^{1/2} D\hbar\pr({\bolds\xi
}_0) D\bbar(\etab_0))
(\Gamb^m_{\varthetab_0})^{-1/2} \Deltab_{\varthetab_0}^m\\
&=&\!: \bar{Q}_{\varthetab_0,1}-\bar{Q}_{\varthetab_0,2}.
\end{eqnarray*}

Since $D\hbar({\bolds\xi}_0)$ has full rank, the standard properties
of Moore--Penrose inverses entail $Q_{{\bolds\xi}_0,1}=\bar
{Q}_{\varthetab_0,1}$. As for $Q_{{\bolds\xi}_0,2}$ and $\bar
{Q}_{\varthetab_0,2}$, they are equal if
\[
\mathcal{M}((\Gamb^m_{\varthetab_0})^{1/2} \mathbf{C}_{\hbar}({\bolds
\xi}_0) D\lbar({\bolds\alpha}_0))
=
\mathcal{M}((\Gamb_{\varthetab_0}^m)^{1/2} D\hbar\pr({\bolds\xi
}_0) D\bbar(\etab_0)).
\]
Since $\Gamb_{\varthetab_0}^m$ and $\mathbf{C}_{\hbar}({\bolds\xi
}_0)$ are invertible, the latter equality holds if
$
\mathcal{M}( D\lbar({\bolds\alpha}_0))
=
\mathcal{M}( (\mathbf{C}_{\hbar}({\bolds\xi}_0))^{-1} D\hbar\pr
({\bolds\xi}_0) D\bbar(\etab_0)),
$
or, since $D\hbar({\bolds\xi}_0)$ has full rank,
if
\begin{eqnarray*}
&&\mathcal{M}( D\hbar({\bolds\xi}_0) D\lbar({\bolds\alpha}_0))
=
\mathcal{M}( D\hbar({\bolds\xi}_0) (\mathbf{C}_{\hbar}({\bolds\xi
}_0))^{-1} D\hbar\pr({\bolds\xi}_0) D\bbar(\etab_0)) \\
&&\hspace*{87.2pt}\bigl(= \mathcal{M}( {\bolds\Pi}(D\hbar({\bolds\xi}_0)) D\bbar(\etab
_0)) \bigr),
\end{eqnarray*}
which trivially holds true. Hence, $Q_{{\bolds\xi}_0,2}=\bar
{Q}_{\varthetab_0,2}$, so that $Q_{{\bolds\xi}_0}=\bar
{Q}_{\varthetab_0}$.

Eventually, the linear spaces orthogonal to $\Gamb^{1/2}_{\varthetab
_0}D\bbar(\etab_0 )$ and to $\Gamb^{1/2}_{\varthetab_0
}D\btilde(\etab_0 )$ do coincide, so that the
statistic $Q_{\varthetab_0}$, which is obtained by
substituting $\btilde$ for $\bbar$ in~(\ref{Qxibar}), is equal
to $Q_{\varthetab_0}$(\mbox{$=$}$Q_{{\bolds\xi}_0})$.
This establishes the result.
\end{pf*}
%
%%%%%%%%%%%%%%%%%%%%%%%%%%%%%%%%%%%%%%%%%%%%%%%%%%%%%%%%%%%%%%%%%%%%%%%%%%%%%%%%%%%%%%%%%%%%%%%%%%%%%%%%%%%%%%%%%%%%%%%%%%%%%%%%%
%%%%%%%%%%%%%%%%%%%%%%%%%%%%%%%%%%%%%%%%%%%%%%%%%%%%%%%%%%%%%%%%%%%%%%%%%%%%%%%%%%%%%%%%%%%%%%%%%%%%%%%%%%%%%%%%%%%%%%%%%%%%%%%

We now turn to the proofs of Lemmas~\ref{infoinverse} and \ref
{parametricasymplin}.
\begin{pf*}{Proof of Lemma~\ref{infoinverse}}
The proof consists in checking
that postmultiplying $\mathbf{D}_{k}(\Lamb_{\Vb})$ with
$\mathbf{N}_{k}\mathbf{H}_{k} \mathbf{P}_{k}^{\Lamb_\Vb} (\mathbf{I}_{k^{2}}+
\mathbf{K}_{k}) \Lamb_{\Vb}^{\otimes2} (\mathbf{P}_{k}^{\Lamb_\Vb
} )\pr\mathbf{H}_{k}\pr\mathbf{N}_{k}\pr$ yields the
$(k-1)$-dimensional identity matrix ($\mathbf{P}_{k}^{\Lamb_\Vb}$ and
$\mathbf{N}_{k}$ are defined in the statement of the lemma). That is, we
show that
%
%e8.10 ###
%
\begin{eqnarray} \label{idone}
&&
\tfrac{1}{4} {\Mb}_{k}^{\Lamb_{\Vb}} \mathbf{H}_{k}
(\mathbf{I}_{k^2} + \mathbf{K}_k)
(\Lamb_{\Vb}^{-1})^{\otimes2} \mathbf{H}_{k}\pr
({\Mb}_{k}^{\Lamb_{\Vb}})\pr
\mathbf{N}_{k}\mathbf{H}_{k}
\nonumber\\
& &\quad{}
\times\mathbf{P}_{k}^{\Lamb_\Vb} (\mathbf{I}_{k^{2}}+ \mathbf{K}_{k})
\Lamb_{\Vb}^{\otimes2} (\mathbf{P}_{k}^{\Lamb_\Vb} )\pr\mathbf{H}_{k}\pr
\mathbf{N}_{k}\pr\\
&&\qquad= \mathbf{I}_{k-1}.\nonumber
\end{eqnarray}

First of all, note that the definition of ${\Mb}_{k}^{\Vb}$
(see Section~\ref{curvLAN}) entails that,
for any $k \times k$ real matrix $\mathbf{l}$ such that $\operatorname
{tr}({\Lamb}_{\Vb}^{-1}\mathbf{l})=0$,
$
(\mathbf{M}_{k}^{\Lamb_{\Vb}} )\pr\mathbf{N}_{k} \mathbf{H}_{k}
(\vecop\mathbf{l})
=
(\mathbf{M}_{k}^{\Lamb_{\Vb}} )\pr(\dvecrond\mathbf{l})
=
\dvec(\mathbf{l})
=
\mathbf{H}_{k} (\vecop\mathbf{l}).
$
Hence, since (letting $\mathbf{E}_{ij}:=\mathbf{e}_{i}\mathbf{e}_{j}\pr
+ \mathbf{e}_{j}\mathbf{e}_{i}\pr$)
\begin{eqnarray*}
\mathbf{P}_{k}^{\Lamb_\Vb} ( \mathbf{I}_{k^{2}}+ \mathbf{K}_{k})
&=& \mathbf{I}_{k^{2}}+ \mathbf{K}_{k}- \frac{2}{k} \Lamb_{\Vb}^{\otimes
2} \vecop(\Lamb_{\Vb}^{-1})(\vecop(\Lamb_{\Vb}^{-1}))\pr
\\[-2pt]
&=& \sum_{i,j=1}^k \vecop\biggl( \frac{1}{2} \mathbf{E}_{ij}- \frac
{1}{k}{(\Lamb_{\Vb}^{-1})_{ij}} \Lamb_{\Vb} \biggr) (\vecop\mathbf{E}_{ij})\pr
\\[-2pt]
&=&\!: \sum_{i,j=1}^k ( \vecop\mathbf{F}_{ij}^{\Lamb_{\Vb}} )
(\vecop\mathbf{E}_{ij})\pr,
\end{eqnarray*}
with
$
\operatorname{tr}( \Lamb_{\Vb}^{-1} \mathbf{F}_{ij}^{\Lamb_{\Vb
}})=0$, for all $i,j=1,\ldots,k$,
we obtain that
$
({\Mb}_{k}^{\Lamb_{\Vb}})\pr
\mathbf{N}_{k}\mathbf{H}_{k} \times\mathbf{P}_{k}^{\Lamb_\Vb}
(\mathbf{I}_{k^{2}}+ \mathbf{K}_{k})
=
\mathbf{H}_{k} \mathbf{P}_{k}^{\Lamb_\Vb} (\mathbf{I}_{k^{2}}+ \mathbf{K}_{k})
$.
Now, using the fact that\break
$\mathbf{H}_{k}\pr\mathbf{H}_{k} (\Lamb_{\Vb}^{-1})^{\otimes2}(\mathbf
{I}_{k^2} + \mathbf{K}_k)
\mathbf{H}_{k}\pr
=(\Lamb_{\Vb}^{-1})^{\otimes2}(\mathbf{I}_{k^2} + \mathbf{K}_k) \mathbf
{H}_{k}\pr$,
the left-hand side\break of~(\ref{idone}) reduces to
%
%e8.11 ###
%
\begin{equation}\label{idtwo}
\tfrac{1}{4} {\Mb}_{k}^{\Lamb_{\Vb}} \mathbf{H}_{k}
(\mathbf{I}_{k^2} + \mathbf{K}_k)
(\Lamb_{\Vb}^{-1})^{\otimes2}
\mathbf{P}_{k}^{\Lamb_\Vb}
(\mathbf{I}_{k^{2}}+ \mathbf{K}_{k}) \Lamb_{\Vb}^{\otimes2}
(\mathbf{P}_{k}^{\Lamb_\Vb} )\pr
\mathbf{H}_{k}\pr\mathbf{N}_{k}\pr.
\end{equation}
After straightforward computation, using essentially the well-known
property of the Kronecker product $\vecop(\mathbf{A}\mathbf{B}\mathbf
{C})=(\mathbf{C}\pr\otimes\mathbf{A}) \vecop(\mathbf{B})$ and the fact that
$\mathbf{M}_{k}^{\Lamb_{\Vb}}\times\mathbf{H}_{k} (\vecop{\Lamb}_{\Vb
}^{-1})=\mathbf{0}$ and $\mathbf{H}_{k} \mathbf{K}_{k}= \mathbf{H}_{k}$,
(\ref{idtwo}) reduces to\break
${\Mb}_{k}^{\Lamb_{\Vb}} \mathbf{H}_{k}
\mathbf{H}_{k}\pr\mathbf{N}_{k}\pr$.
The result follows, since $\mathbf{H}_{k}\mathbf{H}_{k}\pr= \mathbf
{I}_{k}$ and
${\Mb}_{k}^{\Lamb_{\Vb}}\mathbf{N}_{k}\pr=\mathbf{I}_{k-1}$.
\end{pf*}
%
%%%%%%%%%%%%%%%%%%%%%%%%%%%%%%%%%%%%%%%%%%%%%%%%%%%%%%%%%%%%%%%%%%%%%%%%%%%%%%%%%%%%%%%%%%%%%%%%%%%%%%%%%%%%%%%%%%%%%%%%%%%%%%%%
%%%%%%%%%%%%%%%%%%%%%%%%%%%%%%%%%%%%%%%%%%%%%%%%%%%%%%%%%%%%%%%%%%%%%%%%%%%%%%%%%%%%%%%%%%%%%%%%%%%%%%%%%%%%%%%%%%%%%%%%%%%%%%%%%
%%%%%%%%%%%%%%%%%%%%%%%%%%%%%%%%%%%%%%%%%%%%%%%%%%%%%%%%%%%%%%%%%%%%%%%%%%%%%%%%%%%%%%%%%%%%%%%%%%%%%%%%%%%%%%%%%%%%%%%%%%%%%%%%
\begin{pf*}{Proof of Lemma~\ref{parametricasymplin}}
All stochastic convergences in this proof are as $\ny$ under $\mathrm
{P}^{(n)}_{\varthetab; g_{1}}$, for some fixed $\varthetab\in\Thetab
$ and $g_{1} \in{\mathcal F}_{1}^{4}$. It follows from
%
%e8.12 ###
%
\begin{equation} \label{lll1}
\Mb_{k}^{\Lamb_{\Vb}}\mathbf{H}_{k}(\betab\pr)^{\otimes2}(\Vb^{-1})
^{\otimes2}\operatorname{vec} \Vb= \Mb_{k}^{\Lamb_{\Vb}}\mathbf{H}_{k}
(\operatorname{vec} \Lamb_{\Vb} ^{-1})=\mathbf{0}
\end{equation}
and
%
%e8.13 ###
%
\begin{equation} \label{lll2}
\mathbf{L}_{k}^{\betab,\Lamb_{\Vb}} (\Vb^{-1} )^{\otimes2}
\operatorname{vec} {\Vb}= \mathbf{L}_{k}^{\betab,\Lamb_{\Vb}} (\vecop
{\Vb
}^{-1})=\mathbf{0},
\end{equation}
that
\begin{eqnarray*}%\label{centrgaussshape}
\Deltab^{\III} _{\varthetab;\phi_1 }
&=& \frac{a_{k}}{2\sqrt{n}} \Mb_{k}^{\Lamb_{\Vb}}\mathbf{H}_{k}(\Lamb
_{\Vb}^{-1/2} \betab\pr)^{\otimes2}\\[-2pt]
&&{}\times\sum_{i=1}^{n}
\frac{d_{i}^2({\bolds\theta}, \Vb)}{\sigma^2}
\operatorname{vec}(\Ub_{i}({\bolds\theta}, \Vb)\Ub_{i}\pr({\bolds\theta
}, \Vb))\\[-2pt]
&=& \frac{a_{k}}{2\sqrt{n}\sigma^2} \Mb_{k}^{\Lamb_{\Vb}}\mathbf
{H}_{k}(\betab\pr)^{\otimes2}({\Vb}^{-1})^{\otimes2} \\[-2pt]
& &{}
\times\sum_{i=1}^{n} \operatorname{vec}\bigl((\Xb_{i}-{\bolds\theta})(\Xb
_{i}-{\bolds\theta})\pr- \bigl(D_{k}(g_{1})/k\bigr) \Sigb\bigr)
\end{eqnarray*}
and
\begin{eqnarray*}%\label{centrgaussshape}
\Deltab^{\IV} _{\varthetab;\phi_1 }
&=& \frac{a_{k}}{2\sqrt{n}} \mathbf{G}_{k}^{\betab} \mathbf
{L}_{k}^{\betab,\Lamb_{\Vb}} (\Vb^{-1/2})^{\otimes2}\sum_{i=1}^{n}
\frac{d_{i}^2({\bolds\theta}, \Vb)}{\sigma^2}
\operatorname{vec}(\Ub_{i}({\bolds\theta}, \Vb)\Ub_{i}\pr({\bolds\theta
},\Vb))\\
&=& \frac{a_{k}}{2\sqrt{n}\sigma^2} \mathbf{G}_{k}^{\betab} \mathbf
{L}_{k}^{\betab,\Lamb_{\Vb}} (\Vb^{-1})^{\otimes2} \\
& &{}
\times\sum_{i=1}^{n} \operatorname{vec}\bigl((\Xb_{i}-{\bolds\theta})(\Xb
_{i}-{\bolds\theta})\pr- \bigl(D_{k}(g_{1})/k\bigr) \Sigb\bigr).
\end{eqnarray*}
Hence, using a root-$n$ consistent estimator
$\hat{\varthetab}:= (\hat{{\bolds\theta}}\pr, \hat{\sigma
}^{2}, (\dvecrond\hat{\Lamb}_{\Vb})\pr, (\vecop\hat{\betab
})\pr)\pr$ and letting $\hat{\Sigb}:= \hat{\sigma}^{2} \hat
{\betab} \hat{\Lamb}_{\Vb}\hat{\betab}\pr$, Slutsky's lemma yields
\begin{eqnarray*}
\Deltab^{\III} _{\hat\varthetab;\phi_1 }
&=& \frac{a_{k}}{2\sqrt{n}\hat\sigma^2} \Mb_{k}^{\hat{\Lamb
}_{\Vb}}\mathbf{H}_{k}(\hat{\betab}\pr)^{\otimes2}({\hats{\Vb
}}{}^{-1})^{\otimes2}\\
& &{}
\times\sum_{i=1}^{n} \operatorname{vec}\bigl((\Xb_{i}-\hat{{\bolds\theta}}
)(\Xb_{i}-\hat{{\bolds\theta}})\pr- \bigl({D}_{k}(g_{1})/k\bigr)\hat\Sigb
\bigr)\nonumber\\
&=& \frac{a_{k}}{2\sqrt{n}\hat\sigma^2} \Mb_{k}^{\hat{\Lamb
}_{\Vb}}\mathbf{H}_{k}(\hat{\betab}\pr)^{\otimes2}({\hats{\Vb
}}{}^{-1})^{\otimes2} \\
& &{}
\times\Biggl\{\sum_{i=1}^{n} \operatorname{vec}\bigl((\Xb_{i}-{\bolds\theta
})(\Xb_{i}-{\bolds\theta})\pr-\bigl({D}_{k}(g_{1})/k\bigr) \Sigb
\bigr) \\
&&\hspace*{17.8pt}{} - n \operatorname{vec}\bigl((\bar{\Xb} - {\bolds\theta})(\hat{{\bolds
\theta}} - {\bolds\theta})\pr\bigr) - n \operatorname{vec}\bigl((\hat{{\bolds
\theta}} - {\bolds\theta})(\bar{\Xb} - {\bolds\theta})\pr\bigr) \\
&&\hspace*{17.8pt}{} + n \operatorname{vec}\bigl((\hat{{\bolds\theta}} - {\bolds\theta})(\hat
{{\bolds\theta}} - {\bolds\theta})\pr\bigr) - n \bigl({D}_{k}(g_{1})/k\bigr)
\operatorname{vec}(\hat\Sigb- \Sigb) \Biggr\} \\
&=& \Deltab^{\III} _{\varthetab;\phi_1 } - \frac{a_{k}
D_{k}(g_{1})}{2k \sigma^2} \Mb_{k}^{{\Lamb}_{\Vb}}\mathbf{H}_{k}({\betab
}\pr)^{\otimes2}({{\Vb}}^{-1})^{\otimes2} n^{1/2}
\operatorname{vec}(\hat\Sigb- \Sigb)\\
&&{} +o_\mathrm{P}(1),
\end{eqnarray*}
and, similarly,
\begin{eqnarray*}
\Deltab^{\IV} _{\hat\varthetab;\phi_1 }
&=& \frac{a_{k}}{2\sqrt{n}\hat\sigma^2} \mathbf{G}_{k}^{\hat{\betab
}} \mathbf{L}_{k}^{\hat{\betab},\hat{\Lamb}_{\Vb}}({\hats{\Vb }}{}^{-1})^{\otimes2} \\
& &{}
\times\sum_{i=1}^{n} \operatorname{vec}\bigl((\Xb_{i}-\hat{{\bolds\theta}}
)(\Xb_{i}-\hat{{\bolds\theta}})\pr-\bigl({D}_{k}(g_{1})/k\bigr)\hat\Sigb
\bigr)\\
&=& \Deltab^{\IV} _{\varthetab;\phi_1 } - \frac
{a_{k}D_{k}(g_{1})}{2k\sigma^2} \mathbf{G}_{k}^{{\betab}} \mathbf
{L}_{k}^{{\betab},{\Lamb_{\Vb}}}({{\Vb}}^{-1})^{\otimes2} n^{1/2}
\operatorname{vec}(\hat\Sigb- \Sigb)\\
&&{} +o_\mathrm{P}(1).
\end{eqnarray*}
Writing $\hat\Sigb- \Sigb=(\hat\sigma^2-\sigma^2)\hats{\Vb}
+\sigma^2 (\hats{\Vb}-\Vb)$, applying Slutsky's lemma again, and
using~(\ref{lll1}),~(\ref{lll2}) and the fact that $\mathbf{K}_{k}\vecop
(\mathbf{A})= \vecop(\mathbf{A}\pr)$, we obtain
%
%e8.14 ###
%
\begin{eqnarray}\label{astrois}\quad
\Deltab^{\III} _{\hat\varthetab;\phi_1 }
&=& \Deltab^{\III} _{\varthetab;\phi_1 } - \frac{a_{k}
{D}_{k}(g_{1})}{2k} \Mb_{k}^{{\Lamb}_{\Vb}}\mathbf{H}_{k}({\betab}\pr
)^{\otimes2}({{\Vb}}^{-1})^{\otimes2}
n^{1/2} \operatorname{vec}(\hats{\Vb}- \Vb)\nonumber\\
&&{} +o_\mathrm{P}(1) \nonumber\\[-8pt]\\[-8pt]
&=& \Deltab^{\III} _{\varthetab;\phi_1 } - \frac
{a_{k}{D}_{k}(g_{1})}{4k} \Mb_{k}^{{\Lamb}_{\Vb}}\mathbf{H}_{k}({\betab
}\pr)^{\otimes2}({{\Vb}}^{-1})^{\otimes2}
[\mathbf{I}_{k^2}+\mathbf{K}_k]\nonumber\\
& &\hspace*{37.1pt}{}
\times n^{1/2} \operatorname{vec}(\hat
\Vb-
\Vb)+o_\mathrm{P}(1)\nonumber
\end{eqnarray}
and
%
%e8.15 ###
%
\begin{eqnarray}\label{asquatre}
\Deltab^{\IV} _{\hat\varthetab;\phi_1 }
&=& \Deltab^{\IV} _{\varthetab;\phi_1 } - \frac{a_{k}
{D}_{k}(g_{1})}{2k} \mathbf{G}_{k}^{{\betab}} \mathbf{L}_{k}^{{\betab},
{\Lamb_{\Vb}}}({{\Vb}}^{-1})^{\otimes2}
n^{1/2} \operatorname{vec}(\hats{\Vb}- \Vb) +o_\mathrm{P}(1) \nonumber\hspace*{-35pt}\\
&=& \Deltab^{\IV} _{\varthetab;\phi_1 } - \frac
{a_{k}{D}_{k}(g_{1})}{4k} \mathbf{G}_{k}^{{\betab}} \mathbf
{L}_{k}^{{\betab
},{\Lamb_{\Vb}}}({{\Vb}}^{-1})^{\otimes2} [\mathbf{I}_{k^2}+\mathbf{K}_k] n^{1/2} \operatorname{vec}(\hat
\Vb- \Vb)\hspace*{-35pt}\\
&&{}+o_\mathrm{P}(1).\nonumber\hspace*{-35pt}
\end{eqnarray}
Now, \citet{KN93} showed that
\[
n^{1/2}
\pmatrix{\dvec(\hat{\Lamb}_{\Vb}- \Lamb_{\Vb}) \cr
\vecop(\hat{\betab}- \betab)}
=
n^{1/2} \Xib_{\betab, \Lamb_{\Vb}} \vecop( \hats{\Vb}- \Vb)
+ o_\mathrm{P}(1),
\]
where $\Xib_{\betab, \Lamb_{\Vb}}$ was defined in~(\ref{ttttt}).
Similar computations as in the proof of Proposition~\ref{LAN} then yield
%
%e8.16 ###
%
\begin{eqnarray}\label{deltmeth}
&&
n^{1/2} \vecop(\hats{\Vb}- \Vb)
\nonumber\\
&&\qquad=
n^{1/2}
(\Xib_{\betab, \Lamb_{\Vb}}\pr\Xib_{\betab, \Lamb_{\Vb}})^{-1}
\Xib_{\betab, \Lamb_{\Vb}}\pr
\pmatrix{\dvec(\hat{\Lamb}_{\Vb}- \Lamb_{\Vb}) \cr
\vecop(\hat{\betab}- \betab)
}+ o_\mathrm{P}(1)
\nonumber\\[-8pt]\\[-8pt]
&&\qquad=
(\mathbf{L}_{k}^{{\betab},{\Lamb_{\Vb}}})\pr(\mathbf{G}_{k}^{\betab})\pr
n^{1/2} \vecop(\hat{\betab}- \betab) \nonumber\\
& &\qquad\quad{}
+ \betab^{\otimes2} \mathbf{H}_{k} \pr n^{1/2} \dvec( \hat{\Lamb
}_{\Vb}- \Lamb_{\Vb}) + o_\mathrm{P}(1).\nonumber
\end{eqnarray}
The result for $\Deltab^{\III} _{\hat\varthetab;\phi_1}$ then
follows by plugging~(\ref{deltmeth}) into~(\ref{astrois}) and using
the facts that $\mathbf{H}_{k}({\betab}\pr)^{\otimes2}({{\Vb
}}{}^{-1})^{\otimes2}
(\mathbf{L}_{k}^{{\betab},{\Lamb_{\Vb}}})\pr=\mathbf{0}$ and $n^{1/2}
\dvec( \hat{\Lamb}_{\Vb}- \Lamb_{\Vb})=n^{1/2}
(\Mb_{k}^{{\Lamb}_{\Vb}} )\pr\dvecrond( \hat{\Lamb
}_{\Vb}- \Lamb_{\Vb})+o_\mathrm{P}(1)$ as $\ny$ (the latter is a
direct consequence of the definition of $\Mb_{k}^{{\Lamb}_{\Vb}}$
and the delta method). As for the result for $\Deltab^{\IV} _{\hat
\varthetab;\phi_1}$, it follows similarly by plugging (\ref
{deltmeth}) into~(\ref{asquatre}) by noting that $\mathbf
{G}_{k}^{{\betab
}} \mathbf{L}_{k}^{{\betab},{\Lamb_{\Vb}}}({{\Vb}}^{-1})^{\otimes2}
[\mathbf{I}_{k^2}+\mathbf{K}_k]\betab^{\otimes2} \mathbf{H}_{k} \pr
=\mathbf{0}$.
\end{pf*}
\begin{pf*}{Proof of Lemma~\ref{alihoprimeprime}}
Throughout fix $\varthetab=({\bolds\theta}\pr,\sigma^2,(\dvecrond
\Lamb)\pr,(\vecop\betab)\pr)\pr\in\mathcal{H}_{0;q}^{\Lamb
\prime\prime}$ and $g_1\in\mathcal{F}_a$, and define $\tilde{\Vb
}:= \hat{\betab}_{\mathrm{Tyler}} \tilde{\Lamb}_{\Vb} \hat
{\betab}_{\mathrm{Tyler}} ^{ \prime}$. Since $\mathbf{K}_{k} \vecop
(\mathbf{A})= \vecop(\mathbf{A}\pr)$ and
$\mathbf{c}_{p,q} \pr
\mathbf{H}_{k} \hat{\betab}_{\mathrm{Tyler}}^{\prime\otimes2}
({\tilde{\Vb}}{}^{1/2})^{\otimes2} (\vecop\mathbf{I}_{k}) =\mathbf{0}$,
we obtain, from~(\ref{HPres}),
%
%e8.17 ###
%
\begin{eqnarray}
\label{HPasymplin}
{\utT}{}_{K}^{(n)} &=&
\biggl(\frac{nk(k+2)}{\mathcal{J}_k(K)} \biggr)^{1/2}
( a_{p,q} (\tilde{\Lamb}_\Vb))^{-1/2}
\mathbf{c}_{p,q} \pr\mathbf{H}_{k} \hat{\betab}_{\mathrm
{Tyler}}^{\prime\otimes2} ({\tilde{\Vb}}{}^{1/2})^{\otimes2}
\mathbf{J}_{k}^{\perp} \vecop\bigl( {\utSb}{}_{\hat\varthetab;K}^{(n)}\bigr)
\nonumber\hspace*{-30pt}\\
&=& \biggl(\frac{nk(k+2)}{\mathcal{J}_k(K)} \biggr)^{1/2}
( a_{p,q} (\tilde{\Lamb}_\Vb))^{-1/2}
\mathbf{c}_{p,q} \pr\mathbf{H}_{k} \hat{\betab}_{\mathrm
{Tyler}}^{\prime\otimes2} ({\tilde{\Vb}}{}^{1/2})^{\otimes2}
\mathbf{J}_{k}^{\perp} \vecop\bigl( {\utSb}{}_{\varthetab; K}^{(n)}\bigr)
\nonumber\hspace*{-30pt}\\[-8pt]\\[-8pt]
& &{}
- \biggl(\frac{{\mathcal J}_{k}^{2}(K,g_{1})}{4 k(k+2)\mathcal
{J}_k(K)} \biggr)^{1/2}
( a_{p,q} (\tilde{\Lamb}_\Vb))^{-1/2}
\mathbf{c}_{p,q} \pr\mathbf{H}_{k} \hat{\betab}_{\mathrm
{Tyler}}^{\prime\otimes2} ({\tilde{\Vb}}{}^{1/2})^{\otimes2}
\nonumber\hspace*{-30pt}\\
& &\hspace*{11pt}{}
\times({\Vb}^{-1/2})^{\otimes2} n^{1/2}\vecop(\tilde{\mathbf{V}}-
\Vb) + o_\mathrm{P}(1)\nonumber\hspace*{-30pt}
\end{eqnarray}
as $\ny$, under $\mathrm{P}^{(n)}_{\varthetab; g_{1}}$.

We now show that the second term in~(\ref{HPasymplin}) is $o_\mathrm
{P}(1)$ as $\ny$, under $\mathrm{P}^{(n)}_{\varthetab; g_{1}}$.
Since $n^{1/2}\vecop(\tilde{\mathbf{V}}- \Vb)$ is $O_\mathrm{P}(1)$%
%under ${
, Slutsky's lemma yields
\begin{eqnarray*}
&&
( a_{p,q} (\tilde{\Lamb}_\Vb))^{-1/2}
\mathbf{c}_{p,q} \pr\mathbf{H}_{k} \hat{\betab}_{\mathrm
{Tyler}}^{\prime\otimes2} ({\tilde{\Vb}}{}^{1/2})^{\otimes2}
({\Vb}^{-1/2})^{\otimes2} n^{1/2}\vecop(\tilde{\mathbf{V}}- \Vb)\\
&&\qquad
= ( a_{p,q} ({\Lamb}_\Vb))^{-1/2}
\mathbf{c}_{p,q} \pr\mathbf{H}_{k} \hat{\betab}_{\mathrm
{Tyler}}^{\prime\otimes2}\hspace*{1pt} n^{1/2}\vecop(\tilde{\mathbf{V}}- \Vb) +
o_\mathrm{P}(1).
\end{eqnarray*}
By construction of the estimator $\tilde{\Lamb}_{\Vb}$, $\mathbf
{c}_{p,q} \pr\mathbf{H}_{k} \hat{\betab}_{\mathrm{Tyler}}^{\prime
\otimes2} (\vecop\tilde{\mathbf{V}})=0$, so that we have to show that
$n^{1/2}\mathbf{c}_{p,q} \pr\mathbf{H}_{k} \vecop(\hat{\betab}_{\mathrm
{Tyler}}\pr\Vb\hat{\betab}_{\mathrm{Tyler}})$ is $o_\mathrm{P}(1)$.
We only do so for $\varthetab$ values such that
$
\lambda_{1; \Vb}= \cdots= \lambda_{q; \Vb}=: \lambda_{1}^{*} >
\lambda_{2}^{*}:= \lambda_{q+1; \Vb}= \cdots= \lambda_{k; \Vb},
$ which is the most difficult case (extension to the general case is
straightforward, although notationally more tricky).
Note that the fact that $\varthetab\in\mathcal{H}_{0;q}^{\Lamb
\prime\prime}$ then implies that
%
%e8.18 ###
%
\begin{equation}\label{ffff}
-pq\lambda_{1}^{*}+(1-p)(k-q)\lambda_{2}^{*}=0.
\end{equation}
Partition $\mathbf{E}:= \betab\pr\hat{\betab}_{\mathrm
{Tyler}}$ into
%
%e8.19 ###
%
\begin{equation} \label{fsdff}
\mathbf{E}=
\pmatrix{\mathbf{E}_{11} & \mathbf{E}_{12} \cr
\mathbf{E}_{21} & \mathbf{E}_{22}},
\end{equation}
where\vspace*{1pt} $\mathbf{E}_{11}$ is $q \times q$ and $\mathbf{E}_{22}$ is $(k-q)
\times(k-q)$. As shown in Anderson [(\citeyear{A63}), page 129] $n^{1/2}(\mathbf
{E}_{11}\mathbf{E}_{11}\pr- \mathbf{I}_{q})= o_\mathrm
{P}(1)=n^{1/2}(\mathbf{E}_{22}\mathbf{E}_{22}\pr- \mathbf{I}_{k-q})$
and $n^{1/2}\mathbf{E}_{12}=O_\mathrm{P}(1)=n^{1/2}\mathbf{E}_{21}'$ as
$\ny$, under $\mathrm{P}^{(n)}_{\varthetab; g_{1}}$ [actually, \citet{A63} proves this
only for $\mathbf{E}=\betab\pr\betab_{\Sb}$ and under Gaussian
densities, but his proof readily extends to the present situation].
Hence, still as $\ny$, under $\mathrm{P}^{(n)}_{\varthetab; g_{1}}$,
%
%e8.20 ###
%
\begin{eqnarray}
&& n^{1/2} \mathbf{c}_{p,q} \pr\mathbf{H}_{k} \vecop(\hat{\betab}_{\mathrm
{Tyler}}\pr\Vb\hat{\betab}_{\mathrm{Tyler}}) \nonumber\\
&&\qquad
= -p \{ n^{1/2} \lambda_{1}^{*} \operatorname{tr}(\mathbf{E}_{11}\pr
\mathbf{E}_{11}) + n^{1/2} \lambda_{2}^{*} \operatorname{tr}(\mathbf
{E}_{21}\pr\mathbf{E}_{21}) \} \nonumber\\
&&\qquad\quad{}
+ (1-p) \{ n^{1/2} \lambda_{1}^{*} \operatorname{tr}(\mathbf{E}_{12}\pr
\mathbf{E}_{12}) + n^{1/2} \lambda_{2}^{*} \operatorname{tr}(\mathbf
{E}_{22}\pr\mathbf{E}_{22}) \} \\
&&\qquad
= -p \{ n^{1/2} \lambda_{1}^{*} \operatorname{tr}(\mathbf{I}_{q}) \}
+ (1-p) \{ n^{1/2} \lambda_{2}^{*} \operatorname{tr}(\mathbf{I}_{k-q})
\} + o_\mathrm{P}(1)\nonumber\\
&&\qquad
= o_\mathrm{P}(1); \nonumber
\end{eqnarray}
see~(\ref{ffff}).
We conclude that
the second term in~(\ref{HPasymplin}) is $o_\mathrm{P}(1)$, so that
\begin{eqnarray*}
{\utT}{}_{K}^{(n)}
& = &
\biggl(\frac{nk(k+2)}{\mathcal{J}_k(K)} \biggr)^{1/2}
( a_{p,q} (\tilde{\Lamb}_\Vb))^{-1/2}
\mathbf{c}_{p,q} \pr\\
&&{} \times\mathbf{H}_{k} \hat{\betab}_{\mathrm
{Tyler}}^{\prime\otimes2} ({\tilde{\Vb}}{}^{1/2})^{\otimes2}
\mathbf{J}_{k}^{\perp} \vecop\bigl( {\utSb}{}_{\varthetab; K}^{(n)}\bigr)
\\
&&{} + o_\mathrm{P}(1)
\\
&=&
\biggl(\frac{nk(k+2)}{\mathcal{J}_k(K)} \biggr)^{1/2}
( a_{p,q} (\tilde{\Lamb}_\Vb))^{-1/2}
\mathbf{c}_{p,q} \pr\\
&&{} \times\mathbf{H}_{k} \mathbf{E}^{\prime\otimes2} (\betab
\pr)^{\otimes2}
({\tilde{\Vb}}{}^{1/2})^{\otimes2} \mathbf{J}_{k}^{\perp} \vecop\bigl(
{\utSb}{}_{\varthetab; K}^{(n)}\bigr)\\
&&{} + o_\mathrm{P}(1).
\end{eqnarray*}

Since $n^{1/2} \mathbf{J}_{k}^{\perp} \vecop( {\utSb}{}_{\varthetab;
K}^{(n)})$ is $O_\mathrm{P}(1)$ under $\mathrm{P}_{\varthetab;
g_{1}}^{(n)}$, Slutsky's lemma entails
%
%e8.21 ###
%
\begin{eqnarray}\qquad\hspace*{5pt} {\utT}{}_{K}^{(n)}& = &
\biggl(\frac{nk(k+2)}{\mathcal{J}_k(K)} \biggr)^{1/2}
( a_{p,q} ({\Lamb}_\Vb))^{-1/2}
\mathbf{c}_{p,q} \pr\nonumber\\
&&{} \times\mathbf{H}_{k} (\operatorname{diag}(\mathbf
{E}_{11}\pr, \mathbf{E}_{22}\pr))^{\otimes2} (\betab\pr)^{\otimes2}
\nonumber\\
& &{}
\times({{\Vb}}{}^{1/2})^{\otimes2} \mathbf{J}_{k}^{\perp}
\vecop\bigl( {\utSb}{}_{\varthetab; K}^{(n)}\bigr)+
o_\mathrm{P}(1)\nonumber\\[-8pt]\\[-8pt]
& = & \biggl(\frac{nk(k+2)}{\mathcal{J}_k(K)} \biggr)^{1/2}
( a_{p,q} ({\Lamb}_\Vb))^{-1/2}
\mathbf{c}_{p,q} \pr\nonumber\\
&&{} \times\mathbf{H}_{k} (\operatorname{diag}(\mathbf
{E}_{11}\pr, \mathbf{E}_{22}\pr))^{\otimes2} (\betab\pr)^{\otimes2}
\nonumber\\
& &{}
\times({{\Vb}}{}^{1/2})^{\otimes2} \vecop\bigl( {\utSb
}{}_{\varthetab; K}^{(n)}\bigr)+ o_\mathrm{P}(1),\nonumber
\end{eqnarray}
where we used the facts that ${\utSb}{}_{\varthetab; K}^{(n)}$ is
$O_\mathrm{P}(1)$ and that
\begin{eqnarray*}
&&
n^{1/2} \mathbf{c}_{p,q} \pr\mathbf{H}_{k} (\operatorname{diag}(\mathbf
{E}_{11}\pr
, \mathbf{E}_{22}\pr))^{\otimes2} (\betab\pr)^{\otimes2}({{\Vb
}}^{1/2})^{\otimes2} (\vecop\mathbf{I}_{k}) \\
&&\qquad = n^{1/2} \{ -p \lambda_{1}^{*} \operatorname{tr}(\mathbf
{E}_{11}\mathbf{E}_{11}^\prime) + (1-p) \lambda_{2}^{*} \operatorname
{tr}(\mathbf{E}_{22}\mathbf{E}_{22}^\prime) \}
\\
&&\qquad = n^{1/2} \{ - p \lambda_{1}^{*} \operatorname{tr}(\mathbf{I}_{q}) +
(1-p) \lambda_{2}^{*} \operatorname{tr}(\mathbf{I}_{k-q}) \}
+o_\mathrm{P}(1)\\
&&\qquad= o_\mathrm{P}(1).
\end{eqnarray*}
Then, putting [with the same partitioning as in~(\ref{fsdff})]
\[
\betab\pr{\Vb}^{1/2} {\utSb}{}_{\varthetab; K}^{(n)}{\Vb
}^{1/2}\betab=:\mathbf{D}_{\varthetab; K}^{(n)}=:
\pmatrix{\bigl(\mathbf{D}_{\varthetab; K}^{(n)}\bigr)_{11} & \bigl(\mathbf
{D}_{\varthetab; K}^{(n)}\bigr)_{12} \vspace*{2pt}\cr\bigl(\mathbf{D}_{\varthetab;
K}^{(n)}\bigr)_{21} & \bigl(\mathbf{D}_{\varthetab; K}^{(n)}\bigr)_{22}},
\]
the asymptotic properties of ${\utSb}{}_{\varthetab; K}^{(n)}$ and
$\mathbf{E}_{jj}$, $j=1,2$ imply that
\begin{eqnarray*}
{\utT}{}_{K}^{(n)} & =& \biggl(\frac{nk(k+2)}{\mathcal{J}_k(K)}
\biggr)^{1/2}
( a_{p,q} ({\Lamb}_{\Vb}))^{-1/2} \bigl\{ -p \operatorname{tr}\bigl(\mathbf
{E}_{11}\pr\bigl(\mathbf{D}_{\varthetab; K}^{(n)}\bigr)_{11} \mathbf{E}_{11}\bigr)
\\
& &\hspace*{145.6pt}{}
+ (1-p) \operatorname{tr}\bigl(\mathbf{E}_{22}\pr\bigl(\mathbf{D}_{\varthetab;
K}^{(n)}\bigr)_{22} \mathbf{E}_{22}\bigr) \bigr\} +o_\mathrm{P}(1) \\[-4pt]
& =& \biggl(\frac{nk(k+2)}{\mathcal{J}_k(K)} \biggr)^{1/2}
( a_{p,q} ({\Lamb}_{\Vb}))^{-1/2} \bigl\{ -p
\operatorname{tr}\bigl(\bigl(\mathbf{D}_{\varthetab; K}^{(n)}\bigr)_{11} \bigr) \\
& &\hspace*{145pt}{}
+ (1-p) \operatorname{tr}\bigl(\bigl(\mathbf{D}_{\varthetab; K}^{(n)}\bigr)_{22}\bigr) \bigr\}
+o_\mathrm{P}(1)
\\&=& \biggl(\frac{nk(k+2)}{\mathcal{J}_k(K)} \biggr)^{1/2}
( a_{p,q} ({\Lamb}_{\Vb}))^{-1/2}
\mathbf{c}_{p,q} \pr\mathbf{H}_{k} (\betab\pr)^{\otimes2} ({{\Vb
}}^{1/2})^{\otimes2} \vecop\bigl( {\utSb}{}_{\varthetab; K}^{(n)}\bigr)
 +
o_\mathrm{P}(1)
\\
&=&
{\utT}{}_{\varthetab; K}^{(n)}+ o_\mathrm{P}(1)
\end{eqnarray*}
as $\ny$, under $\mathrm{P}^{(n)}_{\varthetab; g_{1}}$, which
establishes the result.
\end{pf*}
\begin{pf*}{Proof of Proposition~\ref{ranktestbeta}} Fix $\varthetab
_0\in
{\mathcal H}_{0;1}^{\betab\prime}$ and $g_1\in\mathcal{F}_a$. We
have already shown in Section~\ref{gsjdlr} that
$
{\utQ}{}^{(n)}_{K}- {\utQ}{}_{\varthetab_0,K}^{(n)}=o_\mathrm{P}(1)
$
as $\ny$ under $\mathrm{P}^{(n)}_{\varthetab_0;g_1}$. Proposition \ref
{Hajek}(i) then yields
%
%e8.22 ###
%
\begin{eqnarray}\label{repres2}\qquad
{\utQ}{}^{(n)}_{K}
&=&
\Deltab^{\IV\prime} _{\varthetab_0 ;K,g_1 }
[(\Gamb^{\IV}_{\varthetab_0 ; K})^{-} - \mathbf{P}_{k}^{\betab_{}^{0}}
( (\mathbf{P}_{k}^{\betab_{}^{0}})\pr\Gamb
^{\IV}_{\varthetab_0 ; K} \mathbf{P}_{k}^{\betab_{}^{0}} )^{-}
(\mathbf{P}_{k}^{\betab_{}^{0}})\pr
]
%({\Gamb}^{\IV}_{\hat{\varthetab};K})^{\perp}
\Deltab^{\IV} _{\varthetab_0 ;K,g_1 }
\nonumber\\[-8pt]\\[-8pt]
&&{} +o_\mathrm{P}(1),\nonumber
\end{eqnarray}
still as $\ny$ under $\mathrm{P}^{(n)}_{\varthetab_0;g_1}$. Now, since
\begin{eqnarray*}
&&
\Gamb^{\IV}_{\varthetab_0;K}
[(\Gamb^{\IV}_{\varthetab_0 ; K})^{-}- \mathbf{P}_{k}^{\betab
_{}^{0}} ( (\mathbf{P}_{k}^{\betab_{}^{0}})\pr\Gamb^{\IV
}_{\varthetab_0 ; K} \mathbf{P}_{k}^{\betab_{}^{0}} )^{-}(\mathbf
{P}_{k}^{\betab_{}^{0}})\pr
] \\
& &\qquad
=\tfrac{1}{2}\mathbf{G}_{k}^{\betab_0}\operatorname{diag}\bigl(\mathbf
{I}_{k-1},\mathbf{0}_{(k-2)(k-1)/2 \times(k-2)(k-1)/2}\bigr)
(\mathbf{G}_{k}^{\betab_0} )\pr
\end{eqnarray*}
is idempotent\vspace*{-4pt} with rank $(k-1)$ [compare with~(\ref{degreefreedom})],
it follows that ${\utQ}{}_{K}^{(n)}$ is asymptotically chi-square\vspace*{2pt} with
$(k-1)$ degrees of freedom
under $\mathrm{P}^{(n)}_{\varthetab_0;g_1}$, which establishes the
null-hypothesis part of (i). For local alternatives, we restrict to
those parameter values $\varthetab_0\in{\mathcal H}_{0}^{\betab}$
for which we have ULAN. From contiguity,~(\ref{repres2}), also holds
under alternatives\vadjust{\goodbreak} of the form $\mathrm{P}^{(n)}_{\varthetab
_0+n^{-1/2}\taub^{(n)};g_1}$. Le Cam's third lemma then implies
that ${\utQ}{}_{K}^{(n)}$, under $ \mathrm{P}^{(n)}_{\varthetab
_0+n^{-1/2}\taub^{(n)};g_1}$, is asymptotically noncentral chi-square,
still with $(k-1)$ degrees of freedom, but
with noncentrality parameter
\begin{eqnarray*}
% \label{limrank2}
&&\lim_{\ny}
\bigl\{ \bigl(\taub^{\IV(n)}\bigr)\pr\\
&&\qquad\hspace*{4pt}{}\times \bigl[ {\Gamb}^{\IV}_{\varthetab
_0;K,g_{1}} [(\Gamb^{\IV}_{\varthetab_0 ; K})^{-}- \mathbf
{P}_{k}^{\betab_{}^{0}} ( (\mathbf{P}_{k}^{\betab_{}^{0}})\pr\Gamb
^{\IV}_{\varthetab_0 ; K} \mathbf{P}_{k}^{\betab_{}^{0}} )^{-}(\mathbf
{P}_{k}^{\betab_{}^{0}})\pr
] {\Gamb}^{\IV}_{\varthetab_0;K,g_{1}}\bigr] \taub^{\IV(n)}
\bigr\}.
\end{eqnarray*}
Evaluation of this limit %in~(\ref{limrank2})
completes part (i) of the proof.

As for parts\vspace*{-4pt} (ii) and (iii), the fact that ${\utphi}{}^{(n)}_{\betab;
K} $ has asymptotic level $\alpha$ directly
follows from the asymptotic null distribution just established and the
classical Helly--Bray theorem, while asymptotic optimality under
$K_{f_1}$ scores is a consequence of the asymptotic equivalence, under
density $f_1$, of ${\utQ}{}^{(n)}_{K_{f_1}}$ and the optimal
parametric test statistic for density $f_1$.
\end{pf*}
%
%%%%%%%%%%%%%%%%%%%%%%%%%%%%%%%%%%%%%%%%%%%%%%%%%%%%%%%%%%%%%%%%%%%%%%%%%%%%%%%%%%%%%%%%%%%%%%%%%%%%%%%%%%%%%%%%%%%%%%%%%%%%%%%%
%%%%%%%%%%%%%%%%%%%%%%%%%%%%%%%%%%%%%%%%%%%%%%%%%%%%%%%%%%%%%%%%%%%%%%%%%%%%%%%%%%%%%%%%%%%%%%%%%%%%%%%%%%%%%%%%%%%%%%%%%%%%%%%%%
%%%%%%%%%%%%%%%%%%%%%%%%%%%%%%%%%%%%%%%%%%%%%%%%%%%%%%%%%%%%%%%%%%%%%%%%%%%%%%%%%%%%%%%%%%%%%%%%%%%%%%%%%%%%%%%%%%%%%%%%%%%%%%%%
\begin{pf*}{Proof of Proposition~\ref{ranktestlambda}}
Fix $\varthetab_0\in{\mathcal H}_{0;q}^{\Lamb\prime\prime}$
and $g_1\in\mathcal{F}_a$. It directly follows from
Lemma~\ref{alihoprimeprime} and Proposition~\ref{Hajek} that
\begin{eqnarray*}
{\utT}{}_{ K}^{(n)}&=&
%&=& {\utT}{}_{\varthetab_0 ; K}\n+o_\mathrm{P}(1) \\
%&=&
( \operatorname{grad}\pr h(\dvecrond\Lamb_{\Vb}^{0}) (\Gamb
^{\III}_{{\bolds\vartheta_0}; K})^{-1}
\operatorname{grad} h(\dvecrond\Lamb_{\Vb}^{0}) )^{-1/2} \\
%%[1mm]
& &{}
\times\operatorname{grad} \pr h(\dvecrond\Lamb_{\Vb}^{0}) (\Gamb
^{\III}_{{\bolds\vartheta_0}; K})^{-1} \Deltab^{\III} _{{\bolds
\vartheta_0}; K, g_{1}} +o_\mathrm{P}(1)
\end{eqnarray*}
as $\ny$, under $\mathrm{P}^{(n)}_{\varthetab_0;g_1}$, hence
also---provided that $\varthetab_0\in{\mathcal H}_{0}^{\Lamb
}$---under the contiguous sequences %of alternatives
$\mathrm{P}^{(n)}_{\varthetab_0+n^{-1/2}\taub^{(n)};g_1}$. Parts (i)
and (ii) result from the fact that ${\utDelta}{}_{\varthetab_0 ;
K,g_{1}}^{\III}$ is asymptotically normal with mean zero under $\mathrm
{P}^{(n)}_{\varthetab_0 ; g_{1}}$ and mean\looseness=-1
\[
\lim_{\ny}
\bigl\{ {\mathcal{J}_k(K,g_{1})}/{\bigl(k(k+2)\bigr)}\mathbf{D}_{k}(\Lamb_{\Vb})
\taub^{\III(n)} \bigr\}
\]\looseness=0
under\vspace*{1pt} $\mathrm{P}^{(n)}_{\varthetab_0 + n^{-1/2} \taub^{(n)}; g_{1}}$
(Le Cam's third lemma; again, for $\varthetab_0\in{\mathcal
H}_{0}^{\Lamb}$), and with covariance matrix $ {\mathcal
{J}_k(K)}/{(k(k+2))} \mathbf{D}_{k}(\Lamb_{\Vb})$ under both.
Parts (iii) and (iv) follow as in the previous proof.
\end{pf*}
\end{appendix}

\section*{Acknowledgments}
The authors very gratefully acknowledge the extremely careful and
insightful editorial handling of this unusually long and technical
paper. The original version received very detailed and constructive
comments from two anonymous referees and a (no less anonymous)
Associate Editor. Their remarks greatly helped improving the
exposition.

\printaddresses

\end{document}